\documentclass[11pt]{article}
\usepackage{amsmath,amsthm}
\usepackage{amsbsy,amssymb}
\usepackage{amsfonts}
\usepackage{graphicx}
\usepackage{epstopdf}
\usepackage[lined,boxed,linesnumbered]{algorithm2e}
\ifpdf
  \DeclareGraphicsExtensions{.eps,.pdf,.png,.jpg}
\else
  \DeclareGraphicsExtensions{.eps}
\fi
\usepackage{booktabs} 
\usepackage{multirow} 
\usepackage{pdflscape}
\usepackage{longtable}
\usepackage{etoolbox,ragged2e}
\makeatletter
\patchcmd{\LT@makecaption}
   {\@tempboxa{#1{#2: }#3}}
   {\@tempboxa{#1{\textsc{#2\ }}\itshape #3}}{}{}
\patchcmd{\LT@makecaption}
   {#1{#2: }#3}
   {\Centering \footnotesize #1\textsc{#2}\par{\itshape #3}}{}{}
\patchcmd{\LT@makecaption}
   {\hbox to\hsize{\hfil\box\@tempboxa\hfil}}
   {\Centering #1\textsc{#2}\par{\itshape #3}}{}{}
\makeatother

\usepackage{graphicx}
\newcommand{\bzero}{{0}}
\newcommand{\ba}{{a}}
\newcommand{\bb}{{b}}
\newcommand{\bc}{{c}}

\newcommand{\bg}{{g}}

\newcommand{\bo}{{o}}

\newcommand{\bt}{{t}}

\newcommand{\bv}{{v}}

\newcommand{\bx}{{x}}
\newcommand{\by}{{y}}

\newcommand{\bA}{{A}}
\newcommand{\bB}{{B}}
\newcommand{\bC}{{C}}
\newcommand{\bD}{{D}}
\newcommand{\bE}{{E}}

\newcommand{\bI}{{I}}
\newcommand{\bJ}{{J}}
\newcommand{\bK}{{K}}
\newcommand{\bL}{{L}}
\newcommand{\bM}{{M}}
\newcommand{\bN}{{N}}
\newcommand{\bO}{{O}}
\newcommand{\bP}{{P}}
\newcommand{\bQ}{{Q}}
\newcommand{\bR}{{R}}
\newcommand{\bS}{{S}}

\newcommand{\bU}{{U}}
\newcommand{\bV}{{V}}
\newcommand{\bW}{{W}}
\newcommand{\bX}{{X}}

\newcommand{\bZ}{{Z}}

\newcommand{\bTheta}{{\Theta}}
\newcommand{\bLambda}{{\Lambda}}

\newcommand{\bSigma}{{\Sigma}}

\def\diag{\mathop{\rm diag}\nolimits}
\def\tr{\mathop{\rm tr}\nolimits}
\def\rank{\mathop{\rm rank}\nolimits}
\def\argmin{\mathop{\rm argmin}\nolimits}
\def\argmax{\mathop{\rm argmax}\nolimits}
\def\dom{\mathop{\rm dom}\nolimits}
\def\dist{\mathop{\rm dist}\nolimits}
\def\corr{\mathop{\rm corr}\nolimits}
\def\vect{\mathop{\rm vec}\nolimits}

\makeatletter 
\newenvironment{varsubequations}[1]
 {
  \addtocounter{equation}{-1}
  \begin{subequations}
  \renewcommand{\theparentequation}{#1-}
  \def\@currentlabel{#1}
 }
 {
  \end{subequations}\ignorespacesafterend
 }
\makeatother

\renewcommand{\S}{Section~}
\newcommand{\Pf}{Proof~}

 \def\citep{\cite}
\def\citet{\cite}

\usepackage{natbib}

\usepackage{color}
\usepackage{url}
\usepackage[margin=1in]{geometry}

\theoremstyle{plain}
\newtheorem{theorem}{Theorem}[section]
\newtheorem{corollary}{Corollary}[section]
\newtheorem{proposition}{Proposition}[section]
\newtheorem{lemma}{Lemma}[section]

\theoremstyle{definition}
\newtheorem{example}{Example}[section]
\newtheorem{definition}{Definition}[section]

\theoremstyle{remark}
\newtheorem{remark}{Remark}[section]

\title{Orthogonal Trace-Sum Maximization: Applications, Local Algorithms, and Global Optimality
}

\author{
	Joong-Ho Won\\
	Department of Statistics\\ 
	Seoul National University\\
	wonj@stats.snu.ac.kr\\
	\and 
	Hua Zhou\\
	Department of Biostatistics\\
	University of California, Los Angeles\\
	huazhou@ucla.edu
	\and
	Kenneth Lange\\
	Departments of Computational Medicine, Human Genetics, and Statistics\\
	University of California, Los Angeles\\
	klange@ucla.edu
}
\date{}

\usepackage{amsopn}

\begin{document}

\maketitle

\begin{abstract}
This paper studies the problem of maximizing the sum of traces of matrix quadratic forms on a product of Stiefel manifolds. 
This orthogonal trace-sum maximization (OTSM) problem generalizes many interesting problems such as generalized canonical correlation analysis (CCA), Procrustes analysis, and cryo-electron microscopy of the Nobel prize fame. 
For these applications finding global solutions is highly desirable but it has been unclear how to find even a stationary point, let alone testing its global optimality.
Through a close inspection of Ky Fan's classical result (1949) on the variational formulation of the sum of largest eigenvalues of a symmetric matrix, 
and a semidefinite programming (SDP) relaxation of the latter,
we first provide a simple method to certify  global optimality of a given stationary point of OTSM. 
This method only requires testing whether a symmetric matrix is positive semidefinite.
A by-product of this analysis is an unexpected strong duality between Shapiro-Botha (1988) and Zhang-Singer (2017).
After showing that a popular algorithm for generalized CCA and Procrustes analysis may generate oscillating iterates, we propose a simple fix that provably guarantees convergence to a stationary point. 
The combination of our algorithm and certificate reveals novel global optima of various instances of OTSM.
 \end{abstract}

\allowdisplaybreaks
\section{Introduction}
\subsection{Orthogonal trace-sum maximization}
Given $\bS_{ij}=\bS_{ji}^T \in \mathbb{R}^{d_i\times d_j}$ for $i,j = 1,\dotsc,m$, 
and $r\le\min_{i=1,\dotsc,m} d_i$, we are interested in solving the following optimization problem: 
\begin{equation}\tag{OTSM}\label{eqn:tracemax}
	\text{maximize} ~~ \frac{1}{2}\sum_{i, j=1}^m \tr(\bO_i^T \bS_{ij} \bO_j) ~~
	\text{subject~to} ~~ \bO_i \in \mathcal{O}_{d_i,r},~ i=1,\dotsc,m
	,
\end{equation}
where 
$\mathcal{O}_{d,r} = \{O \in \mathbb{R}^{d\times r}: \bO^T\bO=\bI_r\}$ is the Stiefel manifold of (partially) orthogonal matrices \citep{Boothby86DiffGeometryBook}; $\bI_r$ denotes the identity matrix of order $r$. 
In the sequel, we call \eqref{eqn:tracemax} the \emph{orthogonal trace-sum maximization} problem. 
OTSM arises in many interesting settings, as follows.

\paragraph{Canonical correlation analysis}
Canonical correlation analysis 
\ifx\citealt\undefined 
	(CCA, \cite{hotelling1936relations})
\else
	\citep[CCA,][]{hotelling1936relations} 
\fi
seeks directions to maximize the correlation between two sets of $n$ observations of variables of possibly different dimensions, $\bA_1\in\mathbb{R}^{n\times d_1}$ and $\bA_2\in\mathbb{R}^{n\times d_2}$:
\[
	\text{maximize} ~~ \corr(\bA_1\bt_1,\bA_2\bt_2) ~~
	\text{subject to} ~~ \bt_i^T\bt_i=1, ~ i=1,2,
\]
where $\bt_i\in\mathbb{R}^{d_i}$ are the optimization variables, and $\corr(\cdot,\cdot)$ denotes the Pearson correlation coefficient between two sample vectors. 
Generalizations of CCA (i) handle more than two sets of variables $\bA_1,\dotsc,\bA_m$ (\mbox{$m\ge 2$}), and (ii) seek partial rotation matrices (as opposed to vectors) of $\bA_i$'s to achieve maximal agreement. 
The popular MAXDIFF and MAXBET criteria \citep{vandegeer1984linear,ten1988generalized,hanafi2006analysis,liu2015maximization} solve
\begin{align*}
	\text{maximize} ~~ \sum_{i<j}\tr(\bO_i^T\bA_i^T\bA_j\bO_j) ~~
	\text{subject to} ~~ \bO_i \in \mathcal{O}_{d_i,r}, ~ i=1,\dotsc,m;
	\tag{MAXDIFF}
	\\
	\text{maximize} ~ \frac{1}{2}\sum_{i,j=1}^m\tr(\bO_i^T\bA_i^T\bA_j\bO_j) ~
	\text{subject to} ~ \bO_i \in \mathcal{O}_{d_i,r}, ~ i=1,\dotsc,m.
	\tag{MAXBET}
\end{align*}
Both MAXDIFF and MAXBET are instances of \eqref{eqn:tracemax} with $\bS_{ij}=\bA_i^T\bA_j$ (if MAXDIFF, $\bS_{ii}=\bzero$), for $i, j=1,\dotsc, m$.
It is worth noting that when $d_1=\dotsb=d_m = d= r$, i.e., the fully orthogonal case, MAXDIFF coincides with MAXBET up to an additive constant.

\paragraph{Procrustes analysis and little Grothendieck problem}
If 
the variables are fully orthogonal 
and there exist $\bS_{11}, \bS_{12}, \dotsc,\bS_{mm} \in \mathbb{R}^{d\times d}$ such that the symmetric $md\times md$ block matrix $\tilde{\bS}=(\bS_{ij})_{i,j=1}^m$ 
is positive semidefinite (denoted $\bS \succeq \bzero$), 
then \eqref{eqn:tracemax} reduces to the little Grothendieck problem over the orthogonal group \citep{BandeiraKennedySinger16LittleGrothendieck}, 
which arises in generalized Procrustes analysis \citep{gower1975generalized,TenBerge1977orthogonal,goodall1991procrustes}.
Given a collection of $n$ landmarks of $d$-dimensional images $\bA_i \in \mathbb{R}^{n\times d}$, $i=1,\dotsc,m$,
the goal is to find orthogonal matrices 
that minimize the pairwise discrepancy
\begin{equation}\label{eqn:procrustes}
	f(\bO_1,\ldots,\bO_m)  =  \frac{1}{2} \sum_{i,j=1}^m \|\bA_i \bO_i - \bA_j\bO_j \|_{\text{F}}^2
	= -\sum_{i<j}\tr(\bO_i^T\bA_i^T\bA_j\bO_j) + \text{const.}
\end{equation} 
subject to the constraints that $\bO_i \in \mathcal{O}_{d,d}$ for all $i$,
where $\|\cdot\|_{\mathrm{F}}$ is the Frobenius norm.
This problem is a special case of \eqref{eqn:tracemax} with $\bS_{ij} = \bA_i^T\bA_j$ for $i,j=1,\dotsc,m$. 
Clearly, $\tilde{\bS} = [\bA_1, \dotsc, \bA_m]^T[\bA_1, \dotsc, \bA_m] \succeq \bzero$.
When $m=2$, problem \eqref{eqn:procrustes} reduces to ordinary (partial) Procrustes analysis \citep[Chapter 7]{DrydenMardia16Book}. 

\paragraph{Cryo-EM and orthogonal least squares}
Another instance of \eqref{eqn:tracemax} 
involving fully orthogonal matrices is
the least squares regression problem 
that minimizes the squared Frobenius norm of the difference between a given $n\times d$ matrix $\bA_{K+1}$ and linear combination $\sum_{i=1}^K \bA_i \bO_i$ of given $n \times d$ matrices $\bA_i$ with $\bO_i \in \mathcal{O}_{d,d}$, $i=1,\dotsc, K$.
This least-squares problem has a direct application in single-particle reconstruction with cryo-electron microscopy (cryo-EM) celebrated by the 2017 Nobel Prize in Chemistry. 
Then we can equivalently minimize
\begin{align}\label{eqn:aug-LScriterion}
	\frac{1}{2}\Big\| \bA_{K+1}(-\bO_{K+1}) -\sum_{i=1}^K \bA_i \bO_i \Big\|_{\text{F}}^2 
	&= \sum_{i < j} \tr(\bO_i^T\bA_i^T \bA_j \bO_j ) + 
	\text{const.}
\end{align}
subject to the orthogonality constraints on $\bO_1, \ldots, \bO_{K+1}$. 
Any minimizer $(\tilde{\bO}_1, \ldots, \allowbreak \tilde{\bO}_{K+1})$ of \eqref{eqn:aug-LScriterion} supplies a minimizer $(- \tilde{\bO}_{1} \tilde{\bO}_{K+1}^T, \ldots, - \tilde{\bO}_{K} \tilde{\bO}_{K+1}^T)$ of 
the original problem.
This is a special case of \eqref{eqn:tracemax} with $\bS_{ij}=-\bA_i^T\bA_j$. 
In cryo-EM, reconstruction of the three-dimensional (3D) map of a particle involves estimating viewing directions of its 2D projections. Retrieval of the orthogonal matrices representing the orientations is posed as the above least squares problem \citep{vainshtein1986,vanheel1987,wang2013orientation,bendory2020}.

\subsection{Global solutions of orthogonal trace-sum maximization}
Each instance of \eqref{eqn:tracemax} above can be posed as a maximum likelihood estimation problem under an appropriate model. 
Finding its global solution is highly desirable for correct inference.
While it attains a maximum because each $\mathcal{O}_{d_i,r}$ is compact and the objective function is continuous in $\mathbb{R}^{d_1\times r}\times\dotsb\times\mathbb{R}^{d_m\times r}$, 
\eqref{eqn:tracemax} is a nonconvex optimization problem since the constraint set 
$\mathcal{O}_{d_1,r}\times\dotsb\times\mathcal{O}_{d_m,r}$ 
is nonconvex.
Except for the special case of $m=2$ in which an analytic global maximizer can be found using the singular value decomposition (SVD) 
\citep{vandegeer1984linear,goodall1991procrustes}, 
we generally have to resort to iterative methods. 
The nonconvexity of the problem makes it difficult to test global optimality of a candidate (local) solution.

To add further difficulties, the global solution to \eqref{eqn:tracemax} is not unique. 
If $(\bO_1^\star,\dotsc,\bO_m^\star)$ is a solution to \eqref{eqn:tracemax}, then for any $\bR\in\mathcal{O}_{r,r}$, $(\bO_1^\star \bR, \dotsc,\bO_m^\star \bR)$ is also a solution.

 \subsection{Contributions}
The contributions of this paper are:
(i) providing a simple certificate that guarantees the \emph{global} optimality of
a \emph{local} stationary point of problem \eqref{eqn:tracemax}
(\S\ref{sec:global}); 
and
(ii) showing that a standard algorithm for generalized CCA and Procrustes analysis may exhibit oscillation, and 
proposing an efficient proximal block relaxation algorithm with a convergence guarantee to a stationary point 
(\S\ref{sec:blockascent}).
Our certificate and duality results are developed in close analogy to the classical result by Ky Fan \citep{fan1949theorem} on the variational formulation of the sum of largest eigenvalues of a symmetric matrix (\S\ref{sec:eigen}). 
(In the Supplementary Material, we also establish a duality between problem \eqref{eqn:tracemax} and another eigenvalue optimization problem.
As a special case, a strong duality between two separately known results in the literature \citep{shapiro1988dual,zhang2017disentangling} is shown.)
The certificate only requires testing positive semidefiniteness of a symmetric matrix constructed from a stationary point and data. Therefore it is simple to verify global optimality.
The convergence theory for the proposed algorithm proves that the \emph{whole} sequence $\{\bO^{k}=(\bO_1^{k},\dotsc,\bO_m^{k})\}$ of iterates converges to a stationary point at least at a sublinear rate ---
this result is stronger than convergence of the objective value sequence or convergence of a subsequence of $\{\bO^{k}\}$.
To our knowledge, there has been no convergence result of this stronger kind for the related problems. 
Some numerical results of the proposed algorithm combined with the certificate are presented in \S\ref{sec:numerical}.

 \section{Preliminary: the $m=1$ case and Ky Fan theorem}\label{sec:eigen}
As a preparation for what follows, we review the classical results on variational formulations of the sum of $r$ largest eigenvalues of a symmetric matrix.
For a matrix $\bS$ in the vector space of $d\times d$ symmetric real matrices (denoted $\mathbb{S}^{d}$), 
let $\lambda_i(\bS)$ be the $i$th largest eigenvalue of $\bS$.
Then it is well known
\begin{equation}\label{eqn:kyfan}
	\sum_{i=1}^r \lambda_i(\bS) = \max_{\bO\in\mathcal{O}_{d,r}} \tr(\bO^T\bS\bO)
= 
\max_{\bU\in\mathbb{S}^d} \left\{ \tr(\bS\bU): \bzero \preceq \bU \preceq \bI_d, \; \tr(\bU)=r \right\}. 	
\end{equation}
The first equality is the celebrated Ky Fan theorem \citep{fan1949theorem}, where the involved nonconvex optimization problem over a Stiefel manifold is a special case of \eqref{eqn:tracemax} for $m=1$.
The second equality is due to \citet{overton1993optimality,hiriart1995sensitivity}, which states that there always is a tight convex relaxation of Ky Fan's nonconvex problem. It is also well known that the dual of this convex semidefinite programming (SDP) problem is an SDP
\begin{equation}\label{eqn:eigendual}
\begin{array}{ll}
\text{minimize} & rz + \tr(\bM)  \\
\text{subject to} 
  & z\bI_d+\bM-\bL = \bS, 
    ~    
  \bM \succeq \bzero, 
    ~
    \bL \succeq \bzero
\end{array}
\end{equation}
for variables $\bM, \bL \in \mathbb{S}^d$ and $z\in\mathbb{R}$
\citep{nesterov1994interior,alizadeh1995interior,pataki1998rank}.

Let us examine the relation between stationary points of these optimization problems closely.
For Ky Fan's nonconvex problem, 
the Lagrangian is 
\[
\mathcal{L}(\bO,\bLambda) = -\tr(\bO^T\bS\bO) + \frac{1}{2}\tr\big[ \bLambda(\bO^T\bO - \bI_r)\big],
\]
by rewriting the constraint $\bO\in\mathcal{O}_{d,r}$ as an equality constraints $\bO^T\bO=\bI_r$.
The Lagrange multiplier matrix $\bLambda$ is symmetric due to the symmetry of the corresponding constraint. 
Point $\bO\in\mathcal{O}_{d,r}$ is a 
stationary point if 
the directional derivative of $\mathcal{L}$
with respect to $\bW\in\mathbb{R}^{d\times r}$
\[
d_{\bW}\mathcal{L} = -\tr\big[(\bS\bO)^T \bW\big] + \tr\big[ (\bO\bLambda)^T \bW\big]  
\]
vanishes for any $\bW$, i.e.,
if $\bO$ satisfies the necessary condition for first-order local optimality.
This is equivalent to the existence of a symmetric matrix $\bLambda$ satisfying
\begin{equation}\label{eqn:eigen1storder}
	\bO\bLambda = \bS\bO.
\end{equation}
Further, using the constraint $\bO^T\bO=\bI_r$, we have a representation $\bLambda=\bO^T\bS\bO \in \mathbb{S}^r$.
The Karush-Kuhn-Tucker (KKT) optimality conditions assert that a locally maximal point is also stationary \citep[Theorem 12.1]{nocedal2006numerical}.

The second-order 
necessary condition for a local maximum
is
\begin{equation}\label{eqn:eigen2ndorder}
d_{\bW}^2\mathcal{L} = \tr(\bLambda\bW^T\bW) - \tr(\bW^T\bS\bW) \ge 0
\end{equation}
for all $\bW \in \mathbb{R}^{d\times r}$ such that $\bW^T\bO + \bO^T\bW = \bzero$
\citep[Theorem 12.5]{nocedal2006numerical};
the set of such $\bW$ is the tangent space of the Stiefel manifold $\mathcal{O}_{d,r}$ at $\bO$.
For the convex problem (either the primal or dual), the 
KKT conditions
are
\begin{equation}\label{eqn:eigenKKT}
\begin{split}
	&\bzero \preceq \bU \preceq \bI_d, \quad \tr(\bU)=r, 
	\quad
	\bM \succeq \bzero, 
	\quad
	\bM + z\bI_d - \bL = \bS, \\
	&\tr(\bL \bU ) = 0, 
	\quad
	\tr[\bM (\bI_d-\bU)] = 0, 
	\quad
	\bL \succeq \bzero.
\end{split}
\end{equation}

Assume that $\bar{\bO}$ is a locally maximizer.
Let $\bar{\bLambda}=\bar{\bO}^T\bS\bar{\bO}$, $\bar{\bU}=\bar{\bO}\bar{\bO}^T$ and $\bar{\bM}=\bar{\bO}(\bar{\bLambda}-\bar{z}\bI_r)\bar{\bO}^T$ for $\bar{z}=\lambda_{\min}(\bar{\bLambda})$, the smallest eigenvalue of $\bar{\Lambda}$. 
Since $\bar{\bO}$ is a stationary point,
it is easy to verify that $(\bar{\bU},\bar{\bM},\bar{\bL}=\bar{\bM}+\bar{z}\bI_d-\bS)$ satisfies all the KKT conditions in \eqref{eqn:eigenKKT} but $\bar{\bL} \succeq \bzero$.\footnote{The construction of matrices $\bar{\bM}$ and $\bar{\bL}$ is inspired by \citet[Theorem 3.3]{pataki1998rank}, in which optimality conditions of the dual SDP \eqref{eqn:eigendual} is given. In particular, if $\bar{\bO}$ consists of the orthonormal eigenvectors of $\bS$ corresponding the the $r$ largest eigenvalues, it follows that $\lambda_{r}(\bS) \ge \bar{z} \ge \lambda_{r+1}(\bS)$, $\bar{\bM}=\bar{\bO}\diag(\lambda_1(\bS) - \bar{z}, \dotsc, \lambda_r(\bS) - \bar{z})\bar{\bO}^T$ and $\bar{\bL}=\bar{\bO}^{\perp}\diag(\bar{z} - \lambda_{r+1}(\bS), \dotsc, \bar{z} - \lambda_d(\bS))\bar{\bO}^{\perp T}$.}
If $r<d$,
let $\bar{\bO}^{\perp}\in\mathcal{O}_{d-r,r}$ 
consist of an orthonormal columns that span the null space of $\bar{\bO}$,
and
choose $\bW=\bar{\bO}^\perp\bK$ for $\bK \in \mathbb{R}^{(d-r)\times r}$ so that it satisfies $\bW^T\bar{\bO} + \bar{\bO}^T\bW = \bzero$.
Then
after some algebra,
\[
	\tr(\bar{\bLambda}\bW^T\bW) - \tr(\bW^T\bS\bW) 
	= 
	\tr(\bK\bar{\bLambda}\bK^T)
	- \tr(\bK^T\bar{\bO}^{\perp T}\bS\bar{\bO}^\perp\bK)
	\ge 0
	.
\]
Since $\bK$ is arbitrary, choose $\bK=\beta v_{\min}^T$ for $\beta\in\mathbb{R}^{d-r}$ and 
$v_{\min} \in \mathbb{R}^r$, where the latter is the unit eigenvector of $\bar{\bLambda}$ associated with $\bar{z}$.
So,
\[
	\tr(\bK\bar{\bLambda}\bK^T)
	- \tr(\bK^T\bar{\bO}^{\perp T}\bS\bar{\bO}^\perp\bK)
	= \bar{z}\beta\beta^T - \beta^T \bar{\bO}^{\perp T}\bS\bar{\bO}^\perp \beta \ge 0,
	\quad
	\forall \beta \in \mathbb{R}^{d-r}.
\]
This implies $\bar{z}\bI_{d-r}- \bar{\bO}^{\perp T}\bS\bar{\bO^\perp} \succeq \bzero$.
On the other hand, as any $y\in\mathbb{R}^d$ can be written as $y=\bar{\bO}\alpha + \bar{\bO}^\perp\beta$ for $\alpha\in\mathbb{R}^r$ and $\beta\in\mathbb{R}^{d-r}$, 
it can be easily seen that
$$
    y^T(\bar{\bM}+\bar{z}\bI_d-\bS)y = \beta^T(\bar{z}\bI_{d-r}-\bar{\bO}^{\perp T}\bS\bar{\bO}^\perp)\beta \ge 0
$$
Thus 
$\bar{\bL}=\bar{\bM}+\bar{z}\bI_d-\bS \succeq \bzero$ and $\bar{\bO}$ satisfies all the KKT conditions.
(If $r=d$, it is immediate that $\bar{\bM}=\bS - \bar{z}\bI_d$, or $\bar{\bL}=\bzero$.)
This implies that all the local maxima of Ky Fan's nonconvex optimization problem are global maxima.

Conversely, if a 
stationary point $\bar{\bO}$ satisfies $\bar{\bL} \triangleq \bar{\bM}+\bar{z}\bI_d-\bS \succeq \bzero$ for $\bar{\bLambda}$, $\bar{\bU}$, and $\bar{\bM}$ constructed as in the beginning of the previous paragraph, then it is also  globally optimal, which obviously is 
locally maximal.
In other words, $\bar{\bL} \succeq \bzero \implies$ (global optimum) $\implies$ (local maximum) $\implies \bar{\bL} \succeq \bzero$, hence $\bar{\bL} \succeq \bzero$ is \emph{necesssary and sufficient} for a stationary point to be globally optimal.

The above analysis of 
Ky Fan's problem \eqref{eqn:kyfan}
sheds light on \eqref{eqn:tracemax} in three ways. 
First, there can be a tight convex relaxation to \eqref{eqn:tracemax}.
Second, by analyzing the KKT conditions of the convex relaxation, we may be able to certify a 
stationary point of \eqref{eqn:tracemax} to be globally optimal.
Third, the dual of the convex relaxation may have to do with a sum of the eigenvalues of a block matrix constructed from $\bS_{ij}$'s.
In the sequel, we carry out an analysis of the OTSM problem inspired by the Ky Fan problem. The added complexity due to $m > 1$ reveals both similarities and differences between the two problems.

 \section{Certificate of global optimality}\label{sec:global}

\subsection{Local optimality conditions}\label{sec:global:local}
\subsubsection{First-order conditions}
Rewriting the constraints $\bO_i\in\mathcal{O}_{d_i,r}$ as equality constraints $\bO_i^T\bO_i=\bI_r$, 
the Lagrangian of \eqref{eqn:tracemax} is 
\[
\mathcal{L}(\bO_1,\dotsc,\bO_m,\bLambda_1,\dotsc,\bLambda_m) = -\frac{1}{2}\sum_{i,j=1}^m\tr(\bO_i^T\bS_{ij}\bO_j) + \frac{1}{2}\sum_{i=1}^m \tr\big[ \bLambda_i (\bO_i^T\bO_i - \bI_r)\big],
\]
where the Lagrange multiplier matrices $\bLambda_i$ are symmetric due to the symmetry of the corresponding constraints.
In parallel with \S\ref{sec:eigen}, 
a point $\bO=(\bO_1,\dotsc,\bO_m)$ is a stationary point of problem \eqref{eqn:tracemax} if 
the directional derivative of $\mathcal{L}$
with respect to $\bW=(\bW_1,\dotsc,\bW_m)\in\mathbb{R}^{d_1\times r}\times\dotsb\times\mathbb{R}^{d_m\times r}$
\begin{equation}\label{eqn:first}
d_{\bW}\mathcal{L} = -\sum_{i=1}^m\sum_{j=1}^m\tr\big[(\bS_{ij}\bO_j)^T \bW_i\big] + \sum_{i=1}^m \tr\big[ (\bO_i\bLambda_i)^T \bW_i\big]  
\end{equation}
vanishes for any $\bW$.
A local maximum satisfies condition \eqref{eqn:first}, which is equivalent to
\begin{equation}\label{eqn:firstorder}
	\bO_i\bLambda_i = \sum_{j=1}^m\bS_{ij}\bO_j, \quad i=1,\dotsc,m,
\end{equation}
resembling the first-order condition \eqref{eqn:eigen1storder} of the Ky Fan problem.

Using the constraint $\bO_i^T\bO_i=\bI_r$, we further have
a representation for $\bLambda_i$:
\begin{equation}\label{eqn:lagrange}
	\bLambda_i = \bO_i^T\big(\sum_{j=1}^m\bS_{ij}\bO_j\big) 
	= \big(\sum_{j=1}^m\bS_{ij}\bO_j\big)^T\bO_i.
\end{equation}
The second equality follows from the symmetry of the Lagrange multiplier. 
Substituting this quantity in equation \eqref{eqn:firstorder}, we obtain
\begin{align}\label{eqn:firstorder2}
	\bO_i\big(\sum_{j=1}^m\bS_{ij}\bO_j\big)^T\bO_i = \sum_{j=1}^m \bS_{ij}\bO_j, 
	\quad i=1,\dotsc,m.
\end{align}

\subsubsection{Second-order condition}
The second-order 
necessary condition for local maximality of \eqref{eqn:tracemax} is 
\begin{equation}\label{eqn:second}
d_{\bW}^2\mathcal{L} = \sum_{i=1}^m \tr(\bLambda_i\bW_i^T\bW_i) - \tr(\bW^T\tilde{\bS}\bW) \ge 0
\end{equation}
for all $\bW=[\bW_1^T, \dotsc, \bW_m^T]^T$,
such that $\bW_i$ is a tangent vector of $\mathcal{O}_{d_i,r}$ at $\bO_i$, i.e.,
\begin{equation}\label{eqn:tangent}
\bW_i^T\bO_i + \bO_i^T\bW_i = \bzero, \quad i=1,\dotsc,m
,
\end{equation}
where
$\tilde{\bS} = (\bS_{ij})$ is the symmetric $D\times D$ block matrix whose $(i,j)$ block is $S_{ij}\in\mathbb{R}^{d_i\times d_j}$, $i, j=1,\dotsc,m$.
Again see the resemblance of condition \eqref{eqn:second} to the second-order condition \eqref{eqn:eigen2ndorder} for Ky Fan's problem.

Unfortunately, in \eqref{eqn:tracemax} we do not have the luxury of all the 
local maxima
being globally optimal. We revisit this issue after studying a potentially tight convex relaxation to the problem in 
the next subsection.
As a partial result, the following characterization of the Lagrange multipliers associated with a 
local maximum
can be deduced from \citet[pp.\;489-490]{ten1988generalized}: 
\begin{proposition}\label{prop:positivity}
	If $(\bO_1,\dotsc,\bO_m)\in\mathcal{O}_{d_1,r}\times \cdots \times\mathcal{O}_{d_m,r}$ is a 
	local maximizer of
	\eqref{eqn:tracemax},
	then $\bLambda_i$ as defined in equation \eqref{eqn:lagrange} is positive semidefinite, for 
	$i=1,\dotsc,m$.
\end{proposition}
\noindent We provide an alternative proof based on the Lagrangian in Appendix \ref{sec:proofs}.

\subsection{Semidefinite programming relaxation}\label{sec:global:sdp}

By introducing an appropriate matrix variable and constraints, we can obtain an upper bound of the optimal value of \eqref{eqn:tracemax} by that of an SDP relaxation. 
Besides providing tight bounds, the SDP formulation paves the way toward certifying the global optimality of a local solution.
If $D=\sum_{i=1}^m d_i$, then we can define a $D\times D$ matrix
\begin{equation}\label{eqn:blockU}
	\bU \triangleq \frac{1}{m} \bO \bO^T, \quad
	\bO = [\bO_1^T, \cdots, \bO_m^T]^T \in \mathbb{R}^{D\times r}
	,
\end{equation}
so that 
$\sum_{i<j}\tr(\bO_i^T\bS_{ij}\bO_j)$,
the objective function of \eqref{eqn:tracemax}, is equal to
$\frac{m}{2}\tr(\tilde{\bS}\bU)$,
where
$\tilde{\bS}=(\bS_{ij})$.
We can express \eqref{eqn:tracemax} in terms of the matrix $\bU$
by imposing appropriate constraints. The proof of the following proposition is in the Appendix.
\begin{proposition}\label{prop:sdp}
Problem \eqref{eqn:tracemax} is equivalent to the optimization problem 
	\begin{equation}\label{eqn:sdprank}
	\begin{array}{ll}
	\text{maximize} & ({m}/{2})\tr(\tilde{\bS}\bU) \\
	\text{subject to} &
		\bU \succeq \bzero, ~
		\rank(\bU)=r, 
	    ~
	    m\bU_{ii} \preceq \bI_{d_i}, ~
		\tr(m\bU_{ii})=r, 
		~
		i=1,\dotsc,m,
	\end{array}
	\end{equation}
where the optimization variable is a symmetric $D\times D$ matrix $\bU$; $\bU_{ii}$ denotes the $i$th diagonal block of $\bU$ whose size is $d_i\times d_i$,
and $\bA \preceq \bB$ denotes the L{\"o}wner order, i.e., $\bB-\bA \succeq \bzero$.
\end{proposition}

By dropping the rank constraint from problem \eqref{eqn:sdprank} we obtain a convex, SDP relaxation of \eqref{eqn:tracemax}:
\begin{equation}\label{eqn:sdp_primal}\tag{P-SDP}
	\begin{array}{ll}
	\text{maximize}  & ({m}/{2})\tr(\tilde{\bS}\bU) \\
	\text{subject to} & 
		\bU \succeq \bzero, 
	    ~
	    m\bU_{ii} \preceq \bI_{d_i}, ~
		\tr(m\bU_{ii}) = r, \quad i=1,\dotsc,m.
	\end{array}
\end{equation}
This relaxation is tight if the solution $\bU^\star$ has rank $r$. The solution to \eqref{eqn:tracemax} is recovered by the decomposition \eqref{eqn:blockU}.
The dual of \eqref{eqn:sdp_primal} is easily seen to be the following SDP
\begin{equation}\label{eqn:sdp_dual}\tag{D-SDP}
\begin{array}{ll}
\text{minimize} & \sum_{i=1}^m [rz_i + \tr(\bM_i) ] \\
\text{subject to} 
  & \bZ+\bM-\bL = \tilde{\bS}, \quad
   \bL \succeq \bzero, 
  \quad
  \bM_i \succeq \bzero, \quad i=1,\dotsc,m, 
\end{array}
\end{equation}
where $\bZ=\diag(mz_1\bI_{d_1},\dotsc,mz_m\bI_{d_m})$ and $\bM=\diag(m\bM_1,\dotsc,m\bM_m)$.
The optimization variables are $\bL\in\mathbb{S}^{D}$, $\bM_i\in\mathbb{S}^{d_i}$, $z_i\in\mathbb{R}$, $i=1,\dotsc,m$.
Strong duality between \eqref{eqn:sdp_primal} and \eqref{eqn:sdp_dual} holds (e.g., Slater's condition is satisfied). 

A rank-$r$ solution to the SDP relaxation \eqref{eqn:sdp_primal}, if it exists, yields a globally optimal solution to the original problem \eqref{eqn:tracemax}. 
However, solving these convex programs is computationally challenging even with modern convex optimization solvers 
due to their lifted dimensions. 
Moreover, if the optimal SDP solution $\bU$ has rank greater than $r$, the factor $\bO$ in \eqref{eqn:blockU} is infeasible to the original problem \eqref{eqn:tracemax}.

Thus it is natural to ask when  the candidate rank-$r$ solution \eqref{eqn:blockU} to \eqref{eqn:sdp_primal} constructed from a stationary point $(\bO_1,\dotsc,\bO_m)$ of \eqref{eqn:tracemax} becomes actually an optimal solution. 
If this is the case, then the local solution globally solves \eqref{eqn:tracemax}. 
We explore this path in the next subsection.

\subsection{Certifying global optimality of a stationary point}\label{sec:global:certificate}
The KKT conditions 
for \eqref{eqn:sdp_primal} and \eqref{eqn:sdp_dual} are 
\begin{varsubequations}{KKT}
\begin{align}
	\bU &\succeq \bzero  \label{eqn:KKT:Ulower} \\
	m\bU_{ii} &\preceq \bI_{d_i},~ i=1,\dotsc,m \label{eqn:KKT:Uupper} \\
	\tr(m\bU_{ii})&=r \label{eqn:KKT:Utrace} \\
	\bM_i &\succeq \bzero, ~ i=1,\dotsc,m \label{eqn:KKT:Mlower} \\
	\bZ + \bM - \bL &= \tilde{\bS}
	\label{eqn:KKT:S} 
	\\
	\tr( \bL \bU ) &= 0 \label{eqn:KKT:Lslack} \\
	\tr( \bM_i(\bI_{d_i}-m\bU_{ii})) &= 0,~ i=1,\dotsc,m \label{eqn:KKT:Mslack} \\
	\bL &\succeq \bzero \label{eqn:KKT:Llower} 
	,
\end{align}
\end{varsubequations}
where $\bZ=\diag(mz_1\bI_{d_1},\dotsc,mz_m\bI_{d_m})$ for $z_1,\dotsc,z_m\in\mathbb{R}$ and $\bM\allowbreak=\allowbreak\diag\allowbreak(m\bM_1,\allowbreak\dotsc,\allowbreak m\bM_m)$.
If any tuple $(\bU,\bZ,\bM,\bL)$ satisfies conditions \eqref{eqn:KKT:Ulower}--\eqref{eqn:KKT:Llower}, then
$\bU$ is an optimal solution to \eqref{eqn:sdp_primal}
and 
$(\bZ,\bM,\bL)$ is
optimal for \eqref{eqn:sdp_dual} 
\citep{vandenberghe1996semidefinite}.

Now suppose $\bar{\bO}=(\bar{\bO}_1,\dotsc,\bar{\bO}_m)$ is a 
stationary point of \eqref{eqn:tracemax}. 
Recalling equation \eqref{eqn:lagrange}, let the associated Lagrange multipliers be $\bar{\bLambda}_i = \sum_{j=1}^m \bar{\bO}_i^T\bS_{ij}\bar{\bO}_j$. 
We can find the quantities that satisfy the KKT conditions above 
in a similar, but not completely obvious, manner to \S \ref{sec:eigen}.
The matrix
\begin{equation}\label{eqn:Ubar}
	\bar{\bU} \triangleq \frac{1}{m}\bar{\bO}\bar{\bO}^T
\end{equation}
clearly satisfies \eqref{eqn:KKT:Ulower}, \eqref{eqn:KKT:Uupper}, and \eqref{eqn:KKT:Utrace}.
Now let
$\tau_i$ be the smallest eigenvalue of the symmetric matrix $\bar{\bLambda}_i$. 
Then
\begin{equation}\label{eqn:Mbar}
	\bar{\bM}_i \triangleq \frac{1}{m}\bar{\bO}_i\bar{\bLambda}_i\bar{\bO}_i^T-\bar{z}_i\bar{\bO}_i\bar{\bO}_i^T = \bar{\bO}_i\left(\frac{1}{m}\bar{\bLambda}_i-\bar{z}_i\bI_r\right)\bar{\bO}_i^T
\end{equation}
satisfies \eqref{eqn:KKT:Mlower}
for any 
$\bar{z}_i \le {\tau_i}/{m}$.
If we define
block diagonal matrices
\begin{equation}\label{eqn:zbar}
	\bar{\bZ} = \diag(m\bar{z}_1\bI_{d_1},\dotsc,m\bar{z}_m\bI_{d_m}), 
	\quad
	\bar{\bM} = \diag(m\bar{\bM}_1,\dotsc,m\bar{\bM}_m), 
\end{equation}
then,
\begin{align*}
	\tr\big[ (\bar{\bM}+\bar{\bZ}) \bar{\bU} \big]
	&= \sum_{i=1}^m \tr\big[ (\bar{\bM}_i+\bar{z}_i\bI_r)m\bar{\bU}_{ii} \big]
	\\
	&= \frac{1}{m}\sum_{i=1}^m \tr\big[ (\bar{\bO}_i\bar{\bLambda}_i\bar{\bO}_i^T-m\bar{z}_i\bar{\bO}_i\bar{\bO}_i^T+m\bar{z}_i\bI_r) \bar{\bO}_i\bar{\bO}_i^T \big] \\
	&= \frac{1}{m}\sum_{i=1}^m\tr(\bar{\bLambda}_i)
	= \frac{1}{m}\sum_{i=1}^m\sum_{j=1}^m\tr(\bar{\bO}_i^T\bS_{ij}\bar{\bO}_j)
	= \tr(\tilde{\bS}\bar{\bU}) 
	.
\end{align*}
This satisfies \eqref{eqn:KKT:S} and \eqref{eqn:KKT:Lslack}
for $\bar{\bL}\triangleq\bar{\bM}+\bar{\bZ}-\tilde{\bS}$.
Finally, 
\begin{align*}
	\tr(m\bar{\bM}_i) &= \tr\big[\bar{\bO}_i(\bar{\bLambda}_i-m\bar{z}_i\bI_r)\bar{\bO}_i^T\big] = \tr(\bar{\bLambda}_i-m\bar{z}_i\bI_r), \\
	\tr\big[ (m\bar{\bM}_i)(m\bar{\bU}_{ii}) \big]
	&= \tr\big[ (\bar{\bO}_i(\bar{\bLambda}_i-m\bar{z}_i\bI_r)\bar{\bO}_i^T)(\bar{\bO}_i\bar{\bO}_i^T) \big]
	= \tr(\bar{\bLambda}_i-m\bar{z}_i\bI_r),
\end{align*}
thus \eqref{eqn:KKT:Mslack} is satisfied.
In short, the choices of variables \eqref{eqn:Ubar}, \eqref{eqn:Mbar}, and \eqref{eqn:zbar} satisfy all the KKT conditions except \eqref{eqn:KKT:Llower} for any $\bar{z}_i\le{\tau_i}/{m}$.

To satisfy this final KKT condition, let 
$\bar{\bL}(\bar{z}_1,\dotsc,\bar{z}_m)=\bar{\bM}+\bar{\bZ}-\tilde{\bS}$
and observe that 
$\bar{\bM}+\bar{\bZ} = \diag\big[\bar{\bO}_1\bar{\bLambda}_1\bar{\bO}_1^T+m\bar{z}_1\bar{\bO}_1^{\perp}\bar{\bO}_1^{\perp T}, \dotsc, \bar{\bO}_m\bar{\bLambda}_m\bar{\bO}_m^T+m\bar{z}_m\bar{\bO}_m^{\perp} \bar{\bO}_m^{\perp T}\big]$,
where $\bO_i^{\perp}\in\mathcal{O}_{d_i,d_i-r}$ satisfies $\bO_i^{\perp}\bO_i^{\perp T} = \bI_{d_i} - \bO_i\bO_i^T$.
If the linear matrix inequality 
(LMI)
\begin{equation}\label{eqn:simplelmi}
	\bar{\bL}(\bar{z}_1,\dotsc,\bar{z}_m) \succeq \bzero, \quad
	\bar{z}_i \le \tau_i/m, \quad i=1,\dotsc,m
\end{equation}
has a feasible point $(\bar{z}_1^\star,\dotsc,\bar{z}_m^\star)$, then this is a certificate that $(\bar{\bO}_1,\dotsc,\bar{\bO}_m)$ is a global maximizer of \eqref{eqn:tracemax}.
While in general LMIs are solved by interior-point methods \citep{nesterov1994interior}, for LMI \eqref{eqn:simplelmi} it is unnecessary. 
Since 
$\bar{\bO}_i^{\perp}\bar{\bO}_i^{\perp T}$
is positive semidefinite, $\bar{\bL}$ is monotone (in Loewner order) in the scalars $\bar{z}_1,\dotsc,\bar{z}_m$, i.e.,
\[
	\bar{\bL}(z_1,\dotsc,z_m) \succeq \bar{\bL}(w_1,\dotsc,w_m)
\]
whenever $z_i\ge w_i$ for $i=1,\dotsc,m$.
Thus it is sufficient to check the positive semidefiniteness at values $\bar{z}_i=\tau_i/m$. If it holds, we have found a feasible point. If not, the LMI is infeasible.
We state this result as the following theorem.
\begin{theorem}\label{thm:certificate}
	Suppose 
	$\bar{\bO}=(\bar{\bO}_1,\dotsc,\bar{\bO}_m)$ is a 
	stationary point of \eqref{eqn:tracemax}.
	Let $\bar{\bLambda}_i=\sum_{j=1}^m \bar{\bO}_i^T\bS_{ij}\bar{\bO}_j$, 
	and 
	$\tau_i$ be the smallest eigenvalue of $\bar{\bLambda}_i$
	for $i=1,\dotsc,m$.
	Then 
	$\bar{\bO}$
	is a global optimum of \eqref{eqn:tracemax} if
\begin{align*}\label{eqn:certificate}\tag{CERT}
	\bL^\star =
    \diag\left(
	\bar{\bO}_1\bar{\bLambda}_1\bar{\bO}_1^T+\tau_1\bar{\bO}_1^{\perp}\bar{\bO}_1^{\perp T}, 
	\dotsc,
	\bar{\bO}_m\bar{\bLambda}_m\bar{\bO}_m^T+\tau_m\bar{\bO}_m^{\perp}\bar{\bO}_m^{\perp T}
	\right)
	- \tilde{\bS}
	\succeq \bzero
	.
\end{align*}
\end{theorem}
\begin{remark}\label{rem:smallercertificate}
Let $\bar{\bO} = {m}^{-1/2}[\bar{\bO}_1^T, \cdots, \bar{\bO}_m^T]^T$. It is easy to see that $\bar{\bO}\in\mathcal{O}_{D,r}$ and $\bL^\star\bar{\bO}=\bzero$
using the 
stationarity
condition \eqref{eqn:firstorder}.
Therefore it suffices to test $(\bar{\bO}^\perp)^T\bL^\star\bar{\bO}^\perp \succeq \bzero$, 
where $\bar{\bO}^\perp\in \mathcal{O}_{D,D-r}$ fills out $\bar{\bO}$ to a fully orthogonal matrix.
This matrix is $(D-r)\times(D-r)$
and may be easier to handle than the $D\times D$ matrix $\bL^\star$.
\end{remark}

\begin{example}\label{ex:hard}
To see the significance of Theorem \ref{thm:certificate}, consider the example examined by 
\ifx\citealt\undefined 
	Ten Berge \cite[p. 270]{TenBerge1977orthogonal}
\else
	\citet[p. 270]{TenBerge1977orthogonal}
\fi
in the context of the MAXDIFF problem ($\bS_{ii}=\bzero$, $i=1,\dotsc,m$). 
Let 
$m=3$, $d_1=d_2=d_3=d$,  
and set $\bS_{12}=-\bI_d$, $\bS_{13}=\bI_d$, $\bS_{23}=\bI_d$.
Using Theorem \ref{thm:certificate},
it is easy to see that any triple of (partially) orthogonal matrices $(\bar{\bO}_1,\bar{\bO}_2,\bar{\bO}_3)\in\mathcal{O}_{d,r}\times\mathcal{O}_{d,r}\times\mathcal{O}_{d,r}$ such that $\bar{\bO}_3=\bar{\bO}_1+\bar{\bO}_2$ satisfies \eqref{eqn:certificate}
for any $r\le d$.
Indeed, for choices $\bLambda_i=\bI_r$ and $\tau_i=1$ for all $i$,
\[
	\bL^\star=\begin{bmatrix} \bI_d & \bI_d & -\bI_d \\
							\bI_d & \bI_d & -\bI_d \\
							-\bI_d & -\bI_d & \bI_d \end{bmatrix}
	= \begin{bmatrix} \bI_d \\ \bI_d \\ -\bI_d \end{bmatrix}
	  \begin{bmatrix} \bI_d & \bI_d & -\bI_d \end{bmatrix}
	\succeq \bzero
	.
\]
Specifically, for $d=3$ and $r=2$, the triple
\begin{equation}\label{eqn:hardoptimum}
	\bar{\bO}_1=\begin{bmatrix} 1 & 0 \\ 0 & 1 \\ 0 & 0 \end{bmatrix},
	\quad
	\bar{\bO}_2=\begin{bmatrix} -1/2 & \sqrt{3}/2 \\ -\sqrt{3}/2 & -1/2 \\ 0 & 0 \end{bmatrix},
	\quad
	\bar{\bO}_3=\begin{bmatrix} 1/2 & \sqrt{3}/2 \\ -\sqrt{3}/2 & 1/2 \\ 0 & 0 \end{bmatrix}
\end{equation}
is a global maximizer.
The global maximum value is $3$.
\end{example}

Unlike the Ky Fan problem of the previous section, condition \eqref{eqn:certificate} is hardly necessary for global optimality.
This is a key difference of \eqref{eqn:tracemax} from the Ky Fan problem.
To see this, let $\bar{\bO}=(\bar{\bO}_1,\dotsc,\bar{\bO}_m)$ be a 
local maximizer,
and
recall
$\bI_{d_i} = \bar{\bO}_i\bar{\bO}_i^T + \bar{\bO}_i^\perp\bar{\bO}_i^{\perp T}$.
It follows from 
the tangency condition \eqref{eqn:tangent} that
\begin{align*}
\tr(\bar{\bLambda}_i\bW_i^T\bW_i) 
& = \tr(\bW_i\bar{\bLambda}_i\bW_i^T\bar{\bO}_i\bar{\bO}_i^T ) +
\tr(\bW_i\bar{\bLambda}_i\bW_i^T\bar{\bO}_i^\perp\bar{\bO}_i^{\perp T} ) \\
& = \tr(\bar{\bO}_i^T\bW_i\bar{\bLambda}_i\bW_i^T\bar{\bO}_i ) +
\tr(\bar{\bLambda}_i\bW_i^T\bar{\bO}_i^\perp\bar{\bO}_i^{\perp T}\bW_i ) \\
& = \tr(\bW_i^T\bar{\bO}_i\bar{\bLambda}_i\bar{\bO}_i^T\bW_i ) +
\tr(\bar{\bLambda}_i\bW_i^T\bar{\bO}_i^\perp\bar{\bO}_i^{\perp T}\bW_i ) 
\ge
\tau_i \tr(\bW_i^T\bar{\bO}_i^\perp\bar{\bO}_i^{\perp ^T}\bW_i )
.
\end{align*}
The last inequality is due to
$\bW_i^T\bar{\bO}_i^\perp(\bar{\bO}_i^{\perp})^T\bW_i \succeq \bzero$ 
and $\bar{\bLambda}_i \succeq \tau_i\bI_r$
\ifx\citealt\undefined 
	(see, e.g., \cite[pp. 482--483]{lange2013}). 
\else
	\citep[see, e.g.,][pp. 482--483]{lange2013}. 
\fi
Thus the local maximality condition \eqref{eqn:second} is sufficiently satisfied if
\begin{align*}
\sum_{i=1}^m
\tr(\bW_i^T\bar{\bO}_i\bar{\bLambda}_i\bar{\bO}_i^T\bW_i ) 
+\sum_{i=1}^m \tau_i 
\tr\big(\bW_i^T\bar{\bO}_i^{\perp}\bar{\bO}_i^{\perp T}\bW_i \big) -\tr(\bW^T\tilde{\bS}\bW)
\ge 0
\end{align*}
for all 
$\bW=[\bW_1^T,\dotsc,\bW_m^T]^T$. 
This will be the case if $\bL^\star \succeq \bzero$, but not only if.

Nevertheless, there are a few special cases that condition \eqref{eqn:certificate} is also necessary for global optimality: 
\begin{corollary}\label{cor:nec_maxdiff}
    For the MAXDIFF problem with $m=2$, i.e., $\bS_{11} = \bzero$ and $\bS_{22} = \bzero$, 
    if a point $(\bar{\bO}_1, \bar{\bO}_2)$ is a global maximizer, then
    condition \eqref{eqn:certificate} is satisfied.
    This is true for any $r \le \min\{d_1, d_2\}$. 
\end{corollary}
\begin{corollary}\label{cor:nec_hanafi}
    For $m=2$ and $r=1$, if a point $(\bar{\bO}_1, \bar{\bO}_2)$ is a global maximizer of \eqref{eqn:tracemax}, 
    then condition \eqref{eqn:certificate} is satisfied.
\end{corollary}
\begin{corollary}\label{cor:nec_data}
    If $\bS_{ij}$ has rank less than or equal to $r$ with singular value decomposition $\bS_{ij} = \bV_i\bSigma_{ij}\bV_j^T$ for $\bV_i \in \mathcal{O}_{d_i,r}$, $\bV_j \in \mathcal{O}_{d_j,r}$, and $\bSigma_{ij}$ is an $r\times r$ nonnegative diagonal matrix for $i, j = 1, \dotsc, m$, then
    $(\bV_1, \dotsc, \bV_m)$ solves \eqref{eqn:tracemax} globally.
\end{corollary}
\noindent The last corollary provides a data-only sufficient condition for global optimality, which does not require computing a stationary point. 
Its hypothesis is satisfied, e.g., if $\bA_i$'s share left singular vectors in the MAXDIFF or MAXBET problems: $\bA_i = \bR \bSigma_i \bV_i^T$ for some $\bR \in \mathcal{O}_{r,r}$.
Corollary \ref{cor:nec_hanafi} is due to \citet[Result 4]{hanafi2003global}, about which we discuss in the next subsection.
Proofs of the other corollaries are provided in Appendix.

We also point out a special case in which local optimality implies condition \eqref{eqn:certificate}. 
For $m=2$ involving fully orthogonal matrices (i.e., $d_1=d_2=r$; note in this case \eqref{eqn:tracemax} coincides with both MAXDIFF and MAXBET), 
it is well-known that 
local maximality
also implies global optimality \citep{fischer1965,TenBerge1977orthogonal}. Theorem \ref{thm:certificate} recovers this result. To see this, 
for a 
locally maximal point
$\bar{\bO}=(\bar{\bO}_1, \bar{\bO}_2)$, observe that $\bar{\bLambda}_1=\bar{\bO}_1^T\bS_{12}\bar{\bO}_2=\bar{\bO}_2^T\bS_{21}\bar{\bO}_1=\bar{\bLambda}_2$ by symmetry. 
Further, $\bar{\bLambda}=\bar{\bLambda}_1=\bar{\bLambda}_2$ is positive semidefinite by Proposition \ref{prop:positivity}. 
Then,
\begin{align*}
    \bL^{\star} = \begin{bmatrix} \bar{\bO}_1\bar{\bLambda}\bar{\bO}_1^T  & -\bS_{12} \\ -\bS_{12}^T & \bar{\bO}_2\bar{\bLambda}\bar{\bO}_2^T \end{bmatrix}
    = 
    \begin{bmatrix} \bar{\bO}_1 \\ -\bar{\bO}_2 \end{bmatrix}
    \bar{\bLambda}
    \begin{bmatrix} \bar{\bO}_1  \\ -\bar{\bO}_2 \end{bmatrix}^T
    \succeq \bzero
    ,
\end{align*}
since $\bS_{12} = \bO_1\bLambda\bO_2^T = \bS_{21}^T$.
Combined with Corollary \ref{cor:nec_maxdiff}, we see the circle
$\bL^{\star} \succeq \bzero \implies$ (global optimum) $\implies$ (local optimum) $\implies$ $\bL^{\star} \succeq \bzero$ for $m=2$ and $r=d_1=d_2$, like the Ky Fan problem.

Table \ref{tbl:parallelism} summarizes the results discussed so far  and their parallelism with those in \S \ref{sec:eigen}.

\subsection{Other certificates}

Sufficient conditions for global optimality of problem \eqref{eqn:tracemax} appear understudied. 
Here we discuss three such conditions collected from the generalized CCA and Procrustes analysis literature.

\ifx\citealt\undefined 
	Ten Berge \cite{TenBerge1977orthogonal} 
\else
	\citet{TenBerge1977orthogonal} 
\fi
shows that if $d_1=\dotsb=d_m=r$ (fully orthogonal) and $\bar{\bLambda}_{ij} \triangleq \bar{\bO}_i^T\bS_{ij}\bar{\bO}_j$ is symmetric and positive semidefinite for all $i<j$ for global optimality, then $(\bar{\bO}_1, \dotsc, \bar{\bO}_m)$ is a global solution.
This sufficient condition is excessively strong, and in fact 
he uses Example \ref{ex:hard} to show that the condition is hardly met: no orthogonal matrices $(\bar{\bO}_1, \bar{\bO}_2, \bar{\bO}_3)$ exist such that $-\bar{\bO}_1^T\bar{\bO}_2$, $\bar{\bO}_1^T\bar{\bO}_3$, and $\bar{\bO}_2^T\bar{\bO}_3$ are simultaneously symmetric and positive semidefinite.
Nevertheless, Ten Berge's sufficient condition is implied by Theorem \ref{thm:certificate}. To see this,  observe that $\bar{\bLambda}_i = \sum_{j=1}^m \bar{\bLambda}_{ij}$ is symmetric (and positive semidefinite), satisfying 
the stationarity condition \eqref{eqn:firstorder}.
Also since $\bS_{ij} = \bar{\bO}_i\bar{\bLambda}_{ij}\bar{\bO}_j^T$,
\begin{align*}
    &\bL^{\star} = \diag(\bar{\bO}_1\bar{\bLambda}_1\bar{\bO}_1, \dotsc, \bar{\bO}_m\bar{\bLambda}_m\bar{\bO}_m) - \tilde{\bS}
    \\
    &= 
    \begin{bmatrix} \bar{\bO}_1 & & \\ & \ddots & \\ & & \bar{\bO}_m \end{bmatrix}
    \begin{bmatrix} \sum_{j\neq 1}\bar{\bLambda}_{1j} & -\bar{\bLambda}_{12} & \dotsb & -\bar{\bLambda}_{1m} \\
    \vdots &  & \ddots & \vdots \\
    -\bar{\bLambda}_{m1} & \dotsb & -\bar{\bLambda}_{m,m-1} & \sum_{j\neq m}\bar{\bLambda}_{mj} 
    \end{bmatrix}
    \begin{bmatrix} \bar{\bO}_1^T & & \\ & \ddots & \\ & & \bar{\bO}_m^T \end{bmatrix}
    .
\end{align*}
It is easy to check that the middle block matrix is positive semidefinite, so is $\bL^{\star}$.

\ifx\citealt\undefined 
	Hanafi and Ten Berge \cite{hanafi2003global} 
\else
	\citet[Result 1]{hanafi2003global} 
\fi
show that if unit vectors $\bo_i \in \mathbb{R}^{d_i}$, $i=1, \dotsc, m$, satisfy $\tilde{\bS}\bo = \diag(\lambda_1, \dotsc, \lambda_m)\bo$ with $\bo = (\bo_1^T, \dotsc, \bo_m^T)^T$ and $\diag(\lambda_1\bI_{d_1}, \dotsc, \lambda_m\bI_{d_m}) - \tilde{\bS} \succeq \bzero$, then $(\bo_1, \dotsc, \bo_m)$ maximizes $\frac{1}{2}\sum_{i,j=1}^m\bo_i^T\bS_{ij}\bo_j$ globally.
Obviously this is a special case of Theorem \ref{thm:certificate} for $r=1$. Corollary \ref{cor:nec_hanafi} follows from this result.

If a stationary point $(\bar{\bO}_1, \dotsc, \bar{\bO}_m)$ satisfies the second-order condition \eqref{eqn:second} for \emph{all}  $\bW=[\bW_1^T, \dotsc, \bW_m^T]^T \in \mathbb{R}^{D\times r}$, i.e., each $\bW_i$ does not necessarily observe the tangency condition \eqref{eqn:tangent} at $\bar{\bO}_i$, then this is obviously sufficient for the point to be globally optimal.
\ifx\citealt\undefined 
	Liu et al. \cite[Theorem 2.4]{hanafi2003global} 
\else
	\citet[Theorem 2.4]{liu2015maximization} 
\fi
describe this condition in a matrix form: 
\begin{equation}\label{eqn:liu}
\mathcal{L}^{\star}
\triangleq
\diag(\bK_{d_1, r}^T (\bI_{d_1} \otimes \bar{\bLambda}_1) \bK_{d_1, r}, \ldots, \bK_{d_m, r}^T (\bI_{d_m} \otimes \bar{\bLambda}_m) \bK_{d_m, r})
-
\mathcal{S}
\succeq \bzero,
\end{equation}
where $\mathcal{S}=(\bI_r \otimes \bS_{ij})$, and $\bK_{mn}$ is the commutation matrix such that $\bK_{mn} \vect{\bA} = \vect{\bA^T}$ for $\bA \in \mathbb{R}^{m\times n}$;
$\vect(\cdot)$ is the usual vectorization operator, and 
$\otimes$ denotes the Kronecker product.
Condition \eqref{eqn:liu} can be related to Theorem \ref{thm:certificate} by the similarity transform of $\mathcal{L}^{\star}$
with $\mathcal{K}=\diag(\bK_{d_1,r}, \dotsc, \bK_{d_m,r})$: 
\begin{align*}
\mathcal{K}
\mathcal{L}^\star 
\mathcal{K}^T
&=  
\diag(
    \bI_{d_1} \otimes \bar{\bLambda}_1,
    \ldots,
    \bI_{d_m} \otimes \bar{\bLambda}_m
)
- 
\mathcal{K}\mathcal{S}\mathcal{K}^T
\\
&= 
\diag(
    \bI_{d_1} \otimes \bar{\bLambda}_1,
    \ldots,
    \bI_{d_m} \otimes \bar{\bLambda}_m
)
-  
\tilde{\mathcal{S}}
\;\succeq\; \bzero,
\quad
\tilde{\mathcal{S}} = (\bS_{ij} \otimes \bI_r).
\end{align*}
The second equality uses the fact $\bK_{d_i,r}(\bI_r \otimes \bS_{ij}) \bK_{d_j,r}^T = \bS_{ij} \otimes \bI_r$. 
Observe the resemblance of the last line to condition \eqref{eqn:certificate}.
When $r=1$, these two conditions actually coincide, hence also with the Hanafi-Ten Berge condition.
For $r>1$, besides the expenses of constructing a larger matrix than $\bL^{\star}$ ($rD\times rD$ vs. $D \times D$), condition \eqref{eqn:liu} is usually stronger than Theorem \ref{thm:certificate}. 
For instance, in Example \ref{ex:portwine} of \S\ref{sec:numerical} stationary points obtained by Algorithm \ref{alg:blockrelax} in \S\ref{sec:blockascent} satisfy condition \eqref{eqn:certificate} for all possible values of $r$, but for those points $\mathcal{L}^{\star} \succeq \bzero$ only when $r=1$.

\begin{table*}[h!]
\caption{Comparison between Ky Fan's problem 
and the orthogonal trace-sum maximization problem 
\eqref{eqn:tracemax}. Set $\Gamma_2$ refers to the set of 
locally maximal points.}
\label{tbl:parallelism}
\begin{center}
\begin{small}
\begin{sc}
\begin{tabular}{@{}p{3.75cm}p{2.75cm}p{5.5cm}@{}}
\toprule
& Ky Fan & 
OTSM \\
\midrule
Domain  & $\bar{\bO}\in\mathcal{O}_{d,r}$          &  $\bar{\bO}=(\bar{\bO}_1,\dotsc,\bar{\bO}_m)\in\mathcal{O}_{d_1,r}\times\dotsb\times\mathcal{O}_{d_m,r}$               \\
Lagrange  multiplier(s) & $\bar{\bLambda}=\bar{\bO}^T\bS\bar{\bO}=\bar{\bLambda}^T$ & $\bar{\bLambda}_i=\sum_{j=1}^m\bar{\bO}_i^T\bS\bar{\bO}_j=\bar{\bLambda}_i^T$,  $i=1,\dotsc, m$ \\
Lifting matrix  &   $\bar{\bU}=\bar{\bO}\bar{\bO}^T$       &  $\bar{\bU}=\frac{1}{m}\bar{\bO}\bar{\bO}^T$                 \\
Cutoff matrix & $\bar{z}=\lambda_{\min}(\bar{\bLambda})$       & 
$\bar{\bZ}=\diag([\tau_i\bI_{d_i}]_{i=1}^m)$, $\tau_i=\lambda_{\min}(\bLambda_i)$    \\ 
Nonneg. part of $\bar{\bLambda}$
& $\bar{\bM}=\bar{\bO}(\bar{\bLambda}-\bar{z}\bI_r)\bar{\bO}^T$ 
& 
$\bar{\bM}=\diag([\bar{\bO}_i(\bar{\bLambda}_i-\tau_i\bI_r)\bar{\bO}_i^T]_{i=1}^m)$ \\
Nonpos. part of $\bar{\bLambda}$
& $\bar{\bL} = \bar{\bM} + \bar{z}\bI_d - \bS$
& $\bar{\bL} = \bar{\bM} + \bar{\bZ} - \tilde{\bS}$ \\
Certificate 
matrix
& $\bar{\bL}\succeq\bzero$ ($\forall\bar{\bO}\in\Gamma_2$)
& $\bar{\bL}\succeq\bzero$ ($\exists\bar{\bO}\in\Gamma_2$)  \\
\bottomrule
\end{tabular}
\end{sc}
\end{small}
\end{center}
\end{table*}
 
\section{Proximal block relaxation algorithm}\label{sec:blockascent}

In order to apply Theorem \ref{thm:certificate} to verify condition \eqref{eqn:certificate}, an algorithm that generates iterates 
converging to a stationary point 
is needed. Although problem \eqref{eqn:tracemax} has been studied in the generalized CCA and Procrustes analysis context for a long time, algorithms that possess this desired property appear rare. 
In this section, we propose such an algorithm.

\subsection{Oscillation of the standard algorithm}\label{sec:blockascent:oscillation}

We first point out a flaw in the  block ascent algorithm 
widely employed in both the generalized CCA  \citep{ten1984orthogonal,ten1988generalized,hanafi2006analysis}
and the Procrustes analysis contexts \citep{TenBerge1977orthogonal,goodall1991procrustes,DrydenMardia16Book}.
This algorithm cyclically updates each orthogonal matrix $\bO_i$ with other blocks $\bO_j$, $j\neq i$, held fixed.  To update the $i$th block in the $k+1$st cycle, let
$\bO^{\text{prev}}=(\bO_1^{k+1},\dotsc,\bO_{i-1}^{k+1},\bO_i^k, \bO_{i+1}^k,\dotsc,\bO_{m}^k)$ then maximize $\tr\big[\bO_i^T(\sum_{j=1}^m\bS_{ij}\bO_j^{\text{prev}})\big]$. 
This block update scheme is natural since the domain $\mathcal{O}_{d_1,r}\times\dotsb\times\mathcal{O}_{d_m,r}$ has a product structure.
Furthermore, each maximization is explicit: let us invoke 
the von Neuman-Fan inequality 
\[
	\tr(\bA^T\bB) \le \sum_{l} \sigma_l(\bA)\sigma_l(\bB),
\]
which holds for any two matrices $\bA$ and $\bB$ of the same dimensions with the $l$th largest singular values $\sigma_l(\bA)$ and $\sigma_l(\bB)$, respectively; equality is attained when $\bA$ and $\bB$ share a simultaneous ordered SVD
\ifx\citealt\undefined 
	(see, e.g., \cite{lange2013}).
\else
	\citep[see, e.g.,][pp. 482-483]{lange2013}.
\fi
Thus if $\bB = \sum_{j=1}^m \bS_{ij}\bO_j^{\text{prev}}$ has a SVD of $\bP_i\bD_i\bQ_i^T$, where $\bD_i$ is $r\times r$ nonnegative diagonal, then the optimal choice of $\bA = \bO_i$ is $\bP_i\bQ_i^T$. The latter matrix is orthogonal.
This method can be considered a linearized version of the alternating variable algorithm of \citet[Algorithm 4.1]{liu2015maximization}.

However, convergence of this standard algorithm is not guaranteed. To be precise, let $\bO^k=(\bO_1^{k},\dotsc,\bO_m^k)$ be the $k$th iterate after $k$ cycles of the algorithm. 
While it can be shown that the sequence of objective values $\{f(\bO^k)\}$ converges 
\citep{ten1988generalized},
it cannot be said that the iterates $\{\bO^k\}$ themselves converge. 
The main reason is that the map $\bB=\sum_{j\neq i}\bS_{ij}\bO_j \mapsto \bP_i\bQ_i^T$ is \emph{set-valued}. 
If $\bB$ is rank deficient, 
any orthonormal basis of the null space of $\bB^T$ (\textit{resp.} $\bB$)  can be chosen as left (\textit{resp.} right) singular vectors corresponding to the zero singular value; the product $\bP_i\bQ_i^T$ may not be unique 
\citep[Proposition 7]{absil2012projection}.
Each update of the $i$th block may place it too far from its previous location.
To see the potential peril of this update scheme,
let us revisit Example \ref{ex:hard}. 
Suppose the algorithm is initialized with $\bO^0=(\bI,\bJ,\bI)$, where
\[
	\bI=\begin{bmatrix} 1 & 0 \\ 0 & 1 \\ 0 & 0 \end{bmatrix}\in\mathcal{O}_{3,2}
	\quad\text{and}\quad
	\bJ=\begin{bmatrix} 0 & 1 \\ 1 & 0 \\ 0 & 0 \end{bmatrix}\in\mathcal{O}_{3,2}
	.
\]
Both $\bI-\bJ$ and $\bI+\bJ$ have rank 1, and $-\bJ\in\argmax_{\bO_i\in\mathcal{O}_{3,2}}\tr\big[\bO_i^T(\bI-\bJ)\big]$, $\bI\in\argmax_{\bO_i\in\mathcal{O}_{3,2}}\allowbreak\tr\big[\bO_i^T(\bI+\bJ)\big]$.
Taking these particular values as the outputs of an instance of the above set-valued map, we have the following sequence of $\bO^k$:
\[
	(\bI,\bJ,\bI) \to (-\bJ,\bI,-\bJ) \to (-\bI,-\bJ,-\bI) 
	\to (\bJ,-\bI,\bJ) \to (\bI,\bJ,\bI) \to \cdots.
\]
Thus the standard algorithm 
\emph{oscillates} while the objective  does not change from a suboptimal value of $1$ (recall the globally optimal value is $3$).

\subsection{Proximal regularization}
We propose a simple modification of the standard algorithm that leads to a convergent algorithm.
Define a bivariate function $f_i:\mathcal{O}_{d_i,r}\times(\mathcal{O}_{d_1,r}\times\dotsb\times\mathcal{O}_{d_m,r})\to\mathbb{R}$ as $f(\bO_i, \bTheta)= \tr\big[\bO_i^T(\sum_{j=1}^m\bS_{ij}\bTheta_j)\big]$ for $i=1,\dotsc,m$, where $\bTheta=(\bTheta_1, \dotsc, \bTheta_m)$.
Then the objective function of \eqref{eqn:tracemax} can be denoted by 
$f(\bO) = \frac{1}{2}\sum_{i=1}^m f_i(\bO_i, \bO)$,
where $\bO=(\bO_1,\dotsc,\bO_m)$.
For the update of the $i$th block in the $k+1$st cycle, 
we 
consider a spherical quadratic approximation of $f_i$ at $\bO_i^k$, i.e., 
\begin{equation}\label{eqn:proximal}
		f_i(\bO_i^k,\bO^{\text{prev}}) + \tr[\nabla_1 f_i(\bO_i^k, \bO^{\text{prev}})^T (\bO_i - \bO_i^{k})] - \frac{1}{2\alpha}\|\bO_i-\bO_i^k\|_{\rm F}^2
		+ \sum_{j\neq i}f_j(\bO_j^{\text{prev}}, \bO^{\text{prev}})
		,
\end{equation}
where $\bO^{\text{prev}}=(\bO_1^{k+1},\dotsc,\bO_{i-1}^{k+1},\bO_i^k, \bO_{i+1}^k,\dotsc,\bO_{m}^k)$,
and $\nabla_1 f_i(\bO_i, \bTheta)=\sum_{j=1}^m \bS_{ij}\bTheta_j$ is the derivative of $f_i$ in the first variable;
$\alpha>0$ is a given constant.
We then maximize this approximation
with respect to $\bO_i$, with the other coordinates held fixed.
This partial maximization is also explicit, as
it can be easily seen that maximizing the objective \eqref{eqn:proximal} is equivalent to maximizing $\tr\big[\bO_i^T(\sum_{j=1}^m\bS_{ij}\bO_j^{\text{prev}}+\alpha^{-1}\bO_i^k)\big]$. 
Therefore we can employ 
the von Neuman-Fan inequality 
to $\bA=\bO_i$ and
$\bB = \sum_{j=1}^m \bS_{ij}\bO_j^{\text{prev}}+{\alpha}^{-1}\bO_i^k$.
If $\bB$
has an SVD of $\bP_i\bD_i\bQ_i^T$, 
then the optimal choice of $\bA$ is again $\bP_i\bQ_i^T$, 
which is orthogonal.
This fact suggests Algorithm \ref{alg:blockrelax}, 
which includes that standard algorithm as a special case ($\alpha=+\infty$).
The quadratic regularization term in objective \eqref{eqn:proximal} keeps the update $\bO_i^{k+1}$ in the proximity of its previous value $\bO_i^k$, and the $\alpha$ moderates the degree of attraction.
Algorithm \ref{alg:blockrelax} is also an instance of the minorization-maximization (MM) algorithm 
\ifx\citealt\undefined 
	(see, e.g., \cite{lange2016mm}): 
\else
	\citep[see, e.g.,][]{lange2016mm}: 
\fi
at each update, the surrogate function defined on $(\mathcal{O}_{d_1,r}\times\dotsb\times\mathcal{O}_{d_m,r})\times(\mathcal{O}_{d_1,r}\times\dotsb\times\mathcal{O}_{d_m,r})$
\[
	g(\bO \mid \bTheta) = f(\bTheta) + \sum_{i=1}^m\tr[\nabla_1 f_i(\bTheta_i, \bTheta)^T (\bO_i - \bTheta_i)] - \frac{1}{2\alpha}\sum_{i=1}^m\|\bO_i-\bO_i^k\|_{\rm F}^2
\]
minorizes the objective function $f(\bO)$ 
at 
$\bTheta=(\bO_1^{k+1}, \allowbreak \dotsc, \allowbreak \bO_{i-1}^{k+1}, \allowbreak \bO_i^k, \allowbreak \bO_{i+1}^k, \allowbreak \dotsc, \allowbreak \bO_m^k)$
for a certain range of $\alpha$
and is partially maximized.
As a consequence of being an MM algorithm,
each update monotonically improves the objective function $f$.
Due to the the compactness of each $\mathcal{O}_{d_i,r}$, actually 
the sequence of objective values $\{f(\bO^k)\}$ converges.

In the next subsection we proceed to show that the sequence of iterates $\{\bO^k\}$ converges to a 
stationary point, in contrast to the standard algorithm.
In particular, in the example of the previous subsection, with any finite $\alpha>0$
the maximizers of $\tr\big[\bO_i^T(-\bJ+\bI+\alpha^{-1}\bI)\big]$, $\tr\big[\bO_i^T(-\bI+\bI+\alpha^{-1}\bJ)\big]$, and $\tr\big[\bO_i^T(\bI+\bJ+\alpha^{-1}\bI)\big]$ in $\mathcal{O}_{3,2}$ are uniquely determined by $\bI$, $\bJ$, and $\bI$. 
This yields $\bO^0=(\bI,\bJ,\bI)=\bO^1=\bO^2=\dotsb$ in Algorithm \ref{alg:blockrelax}. In fact the point $(\bI,\bJ,\bI)$ is a 
stationary point.
(Convergence to a global maximizer, e.g., \eqref{eqn:hardoptimum}, requires a good initial point. We discuss this in \S\ref{sec:numerical} with another global solution.)
Note, however,
the map in Lines 5 and 6 of Algorithm \ref{alg:blockrelax} is nevertheless set-valued, since there is no guarantee of full rank of $\bB$.

\begin{algorithm}[h!]
\caption{Proximal block relaxation algorithm for solving \eqref{eqn:tracemax}}
\label{alg:blockrelax}
\begin{tabbing}
{\small~~1:}\enspace Init\=ialize \= $\bO_1,\dotsc,\bO_m$; Set $\alpha\in(0, 1/\max_{i=1,\dotsc,m}\|\bS_{ii}\|_2)$\\
{\small~~2:}\enspace For $k=1,2,\dotsc$ \\
{\small~~3:}\enspace\> For $i=1,\dotsc,m$ \\
{\small~~4:}\enspace\>\> Set $\bB = \sum_{j=1}^m\bS_{ij}\bO_j+{\alpha}^{-1}\bO_i$ \\
{\small~~5:}\enspace\>\> Compute SVD of $\bB$ as $\bP_i\bD_i\bQ_i^T$ \\
{\small~~6:}\enspace\>\> Set $\bO_i = \bP_i\bQ_i^T$ \\
{\small~~7:}\enspace\> End For \\
{\small~~8:}\enspace\> If there is no progress, then break \\
{\small~~9:}\enspace End For \\
{\small~10:}\enspace Return $(\bO_1,\dotsc,\bO_m)$
\end{tabbing}
\end{algorithm}

 \subsection{Global convergence}\label{sec:blockascent:convergence}

Algorithm \ref{alg:blockrelax} with $\alpha>0$ converges 
despite of the non-uniqueness of the map in Lines 5 and 6:
\begin{theorem}\label{thm:convergence}
	The sequence $\{(\bO_1^{k},\dotsc,\bO_m^{k})\}$ generated by Algorithm \ref{alg:blockrelax} converges to a 
	stationary point 
	of \eqref{eqn:tracemax} for 
	$\alpha \in (0, 1/\max_{i=1,\dotsc,m}\|\bS_{ii}\|_2)$;
	$\|\cdot\|_2$ denotes the spectral norm.
	Furthermore, the rate of convergence is at least sublinear.
\end{theorem}
This result is stronger than typical global convergence results that all the limit points are stationary \citep{zangwill1969nonlinear,lange2013}, or that the gradient vanishes \citep{nocedal2006numerical,absil2007genrtr}.
Theorem \ref{thm:convergence}
can be shown using Theorems 1 and 2 in 
\ifx\citealt\undefined 
	Xu and Yin \cite{xu2017globally} 
\else
	\citet{xu2017globally} 
\fi
by noting that Algorithm \ref{alg:blockrelax} falls into their ``deterministic block prox-linear'' class of algorithms and problem \eqref{eqn:tracemax} possesses the Kurdyka-\L{}ojasiewicz property \citep{attouch2010proximal}.
In the Supplementary Material, 
we provide a simpler proof utilizing the closedness \citep{zangwill1969nonlinear} of the map in Lines 5 and 6, and the geometry of the product of Stiefel manifolds. 
\begin{remark}
    In case $\bS_{ii}=\bzero$ for $i=1,\dotsc,m$, e.g., the MAXDIFF problem, the $\alpha$ can be chosen as an arbitrary positive constant.
\end{remark}
 \section{Numerical experiments}\label{sec:numerical}
\subsection{Setup}\label{sec:numerical:setup}
In this section we test Algorithm \ref{alg:blockrelax} equipped with the certificates of global optimality and suboptimality discussed in \S\ref{sec:global} with both synthetic and real-world examples.
Algorithm \ref{alg:blockrelax} was implemented in the Julia programming language \citep{bezanson2017julia} and run on a standard laptop computer (Macbook Pro, i5@2.4GHz, 16GB RAM).
We set the proximity constant $\alpha=1000$ and terminated the algorithm if the mean change ${m}^{-1}\sum_{i=1}^m\|\bO_i^k-\bO_i^{k-1}\|_{\rm F}$ was less than $10^{-8}$ and the relative change of the objective function was less than $10^{-10}$,
or a maximum iteration of 50000 was reached.
For comparison, we also tested the generic Riemanian trust-region method \citep{absil2007genrtr} implemented in the Manopt MATLAB toolbox \citep{boumal2014manopt}. The maximum number of outer iterations of this method was set to 10000.
For both methods, four initialization strategies were considered:
\begin{enumerate}
    \item (``eye") The $i$th block $\bO_i^0$ of the initial point takes the first $r$ columns of $\bI_{d_i}$. 
    \item (``tb") 
Take the eigenvectors corresponding to the $r$ largest eigenvalues of
the data matrix 
$\tilde{\bS}=(\bS_{ij})$
to form a $D\times r$ orthogonal matrix $\tilde{\bV}$. Split $\tilde{\bV}$ into $m$ blocks so that $\tilde{\bV}=[\tilde{\bV}_1^T,\dotsc,\tilde{\bV}_m^T]^T$ where $\tilde{\bV}_i\in\mathbb{R}^{d_i\times r}$, $i=1,\dotsc,m$. Project each block $\tilde{\bV}_i$ to the Stiefel manifold $\mathcal{O}_{d_i,r}$ to obtain $\bO_i^0$. 
    \item (``sb") Replace the diagonal blocks $\bS_{ii}$ of 
    $\tilde{\bS}$ by $-\sum_{j=1}^m(\bS_{ij}\bS_{ij}^T)^{1/2}$, where $\bM^{1/2}$ denotes the matrix square root of the positive semidefinite matrix $\bM$. 
    Take the eigenvectors corresponding to the $r$ largest eigenvalues 
    of the resulting negative semidefinite matrix
    to form a $D\times r$ orthogonal matrix $\tilde{\bV}$. 
    Proceed as strategy ``tb.''
    \item (``lww1") Set $\bO_1^0$ to the top $r$ eigenvectors of $\bS_{11}$. Then set $\bO_k^0 = \bU_k \bQ_k$, where $\bU_k$ is the 
top $r$ eigenvectors of $\bS_{kk}$ and $\bQ_k$ is the Q factor in the QR decomposition of $\bU_k \sum_{j<k} \bS_{kj}\bO_j^0$, $k=2, \dotsc, m$.
\end{enumerate}
The initial point of strategy ``tb'' coincides with that gives the second upper bound of the orthogonal Procrustes problem \citep[p. 273]{TenBerge1977orthogonal}, and also with the starting point strategy 2 of \citet[p. 1495]{liu2015maximization} for the MAXBET problem.
Strategy ``sb'' extends \citet[p. 380]{shapiro1988dual} for the orthogonal Procrustes problem; 
see Lemma \ref{lem:shapirobotha} in Appendix,
and also the Supplementary Material.
Strategy ``lww1'' is the starting point strategy 1 by \citet[p. 1494]{liu2015maximization}. 

\subsection{Small examples}\label{sec:numerical:small}
\begin{example}[CCA of Port Wine Data]\label{ex:portwine}
We consider generalized CCA of the subset of the data from sensory evaluation of port wines 
analyzed by 
\ifx\citealt\undefined 
	Hanafi and Kiers \cite[Table 2]{hanafi2006analysis}.
\else
	\citet[Table 2]{hanafi2006analysis}.
\fi
The goal is to capture the agreement between $m=4$ assessors 
in the assessment of the appearance of $n= 8$ port wines. 
Note the dimensions are disparate: $d_1=4$, $d_2=3$, $d_3=4$, and $d_4=3$.
The MAXDIFF criterion was tested for all possible $r=1, 2, 3$. 
The results are summarized in Table \ref{tbl:portwine_maxdiff}.
Algorithm \ref{alg:blockrelax} achieved global optimum for all $r$ and for all initial point strategies. A similar phenomenon occurred with MAXBET,
whose results are provided in the Supplementary Material,
except for $r=3$ with strategies ``eye'' and ``lww1.''
    On the other hand, Manopt occasionally 
    converged to a stationary point confirmed to be not locally maximal (violating the conclusion of Proposition \ref{prop:positivity}). 
    In addition, 
    Algorithm \ref{alg:blockrelax} was orders of magnitudes faster than Manopt.
    In all cases certified to be globally optimal, the smallest eigenvalues of $\bL^{\star}$ in condition \eqref{eqn:certificate} (denoted $\lambda_{\min}(\bL^{\star})$) were numerically zero up to the fourteenth digit after the decimal point, whereas those of $\mathcal{L}^{\star}$ (denoted $\lambda_{\min}(\mathcal{L}^{\star})$) in condition \eqref{eqn:liu} were often definitely negative.
\end{example}

\begin{example}
\label{ex:hard2}
    We revisit Example \ref{ex:hard} for $d=3$ and $r=1, 2, 3$.
    The results are summarized in Table \ref{tbl:tenberge}.
    Strategies ``eye'' and ``lww1'' gave suboptimal stationary points as initial points, and both algorithms could not make a  progress.
    Strategy ``sb'' yielded global optima for both $r=1, 2$ 
    For $r=3$, no strategy could certify global optimality. 
    For $r=2$, while both ``sb'' and ``tb'' were successful,  Algorithm \ref{alg:blockrelax} took the full 50000 iterations to achieve the same accuracy as Manopt, which in this case took 22 outer iterations.
    The stationary points reached from the two initial points were all quite different from each other, and also from the analytic solution \eqref{eqn:hardoptimum}. 
    The error $\|\bar{\bO}_1+\bar{\bO}_2-\bar{\bO}_3\|_{\infty}$ was between $8.877\times 10^{-7}$ and $6.146\times 10^{-6}$. 
    Together with the smallest eigenvalue of $\bL^\star$ computed being $-6.674\times 10^{-6}$, this relatively large error reflects the hardness of this problem illustrated in 
    Example \ref{ex:hard}.
    This difficulty was also experienced with an extra run of the commercial interior-point method solver MOSEK \citep{mosek} to solve the convex relaxtion \eqref{eqn:sdp_primal}. 
    While the optimal objective value was $3$ up to the eighth digit after the decimal point,  MOSEK   failed to obtain a rank-two solution.
    With this exception, Algorithm \ref{alg:blockrelax} terminated in a fraction of time for Manopt.
\end{example}

\paragraph{Additional examples}
    In the Supplementary Material, Examples 5.1 and 5.2 of \citet{liu2015maximization} are considered under both MAXDIFF and MAXBET criteria, and new global optima are found.

\subsection{Simulation studies}\label{sec:numerical:large}
\begin{figure}[t!]
\begin{center}
\begin{tabular}{ll}
\includegraphics[height=0.18\textheight]{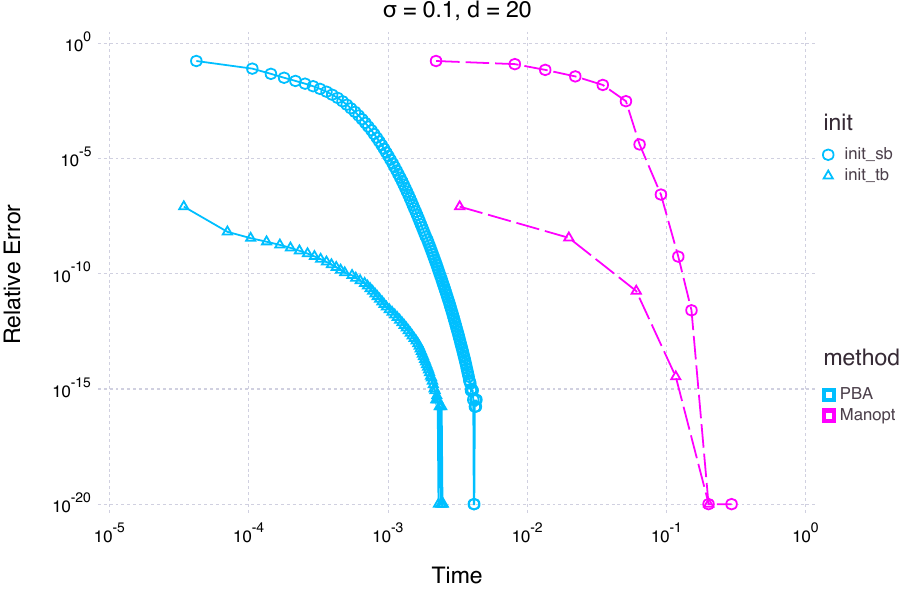}
&
\includegraphics[height=0.18\textheight]{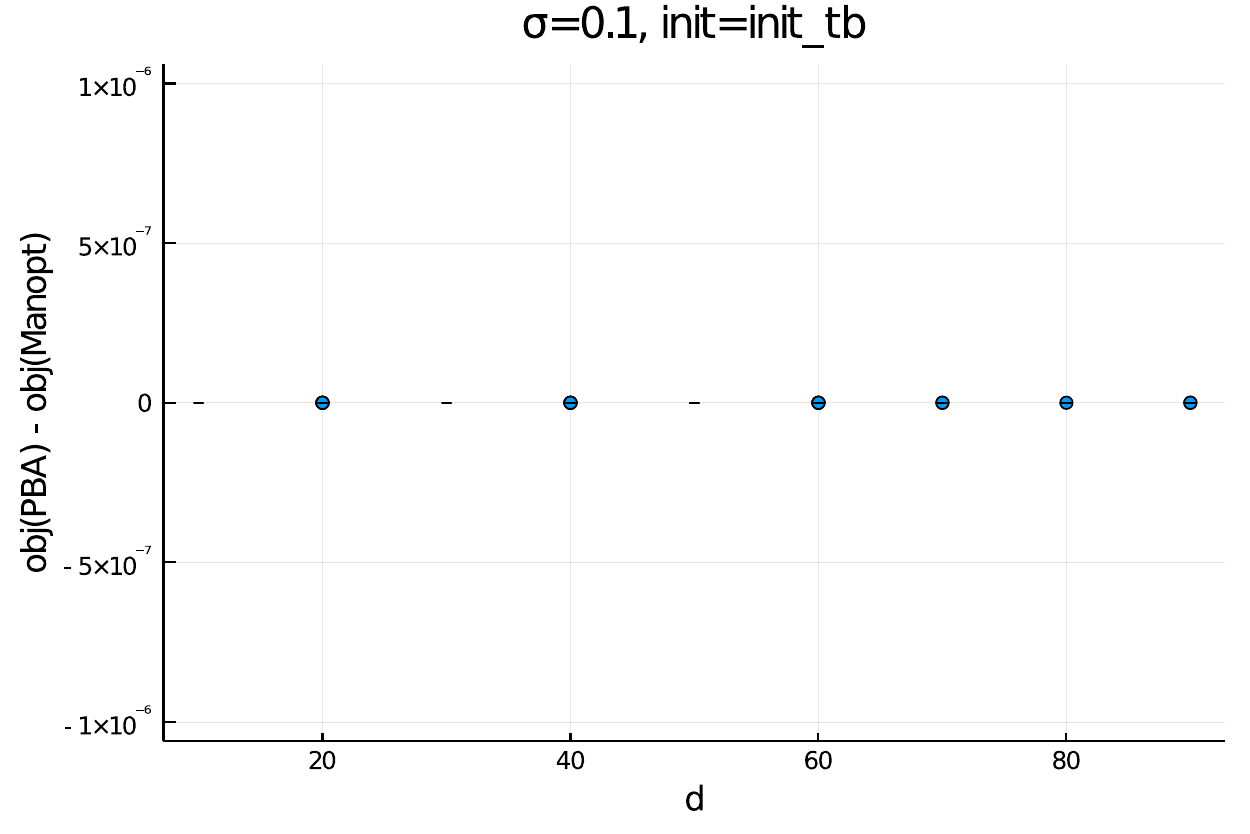}
\\
\includegraphics[height=0.18\textheight]{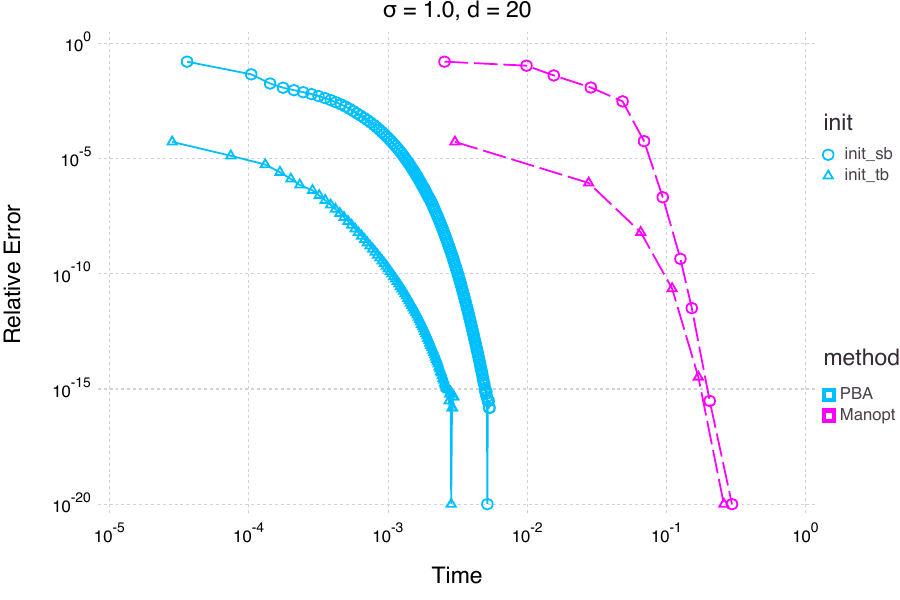}
&
\includegraphics[height=0.18\textheight]{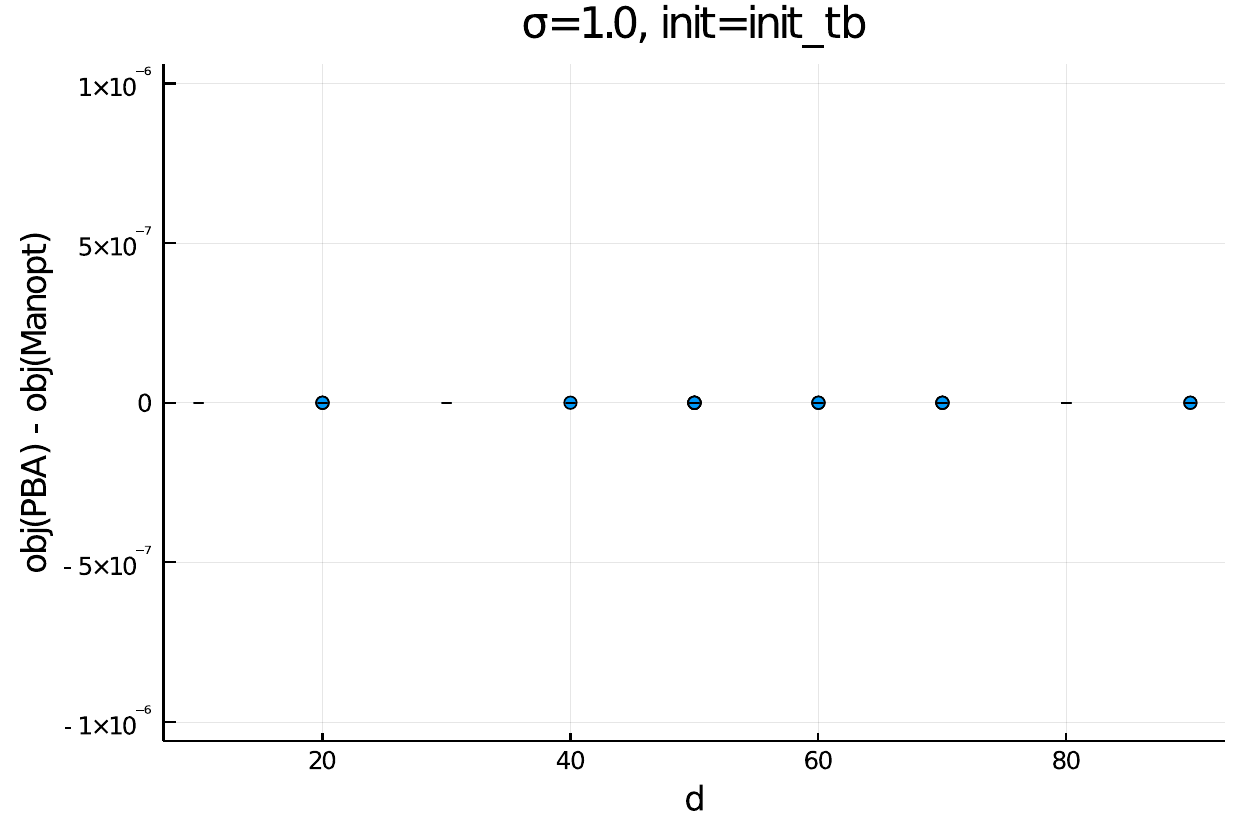}
\\
\includegraphics[height=0.18\textheight]{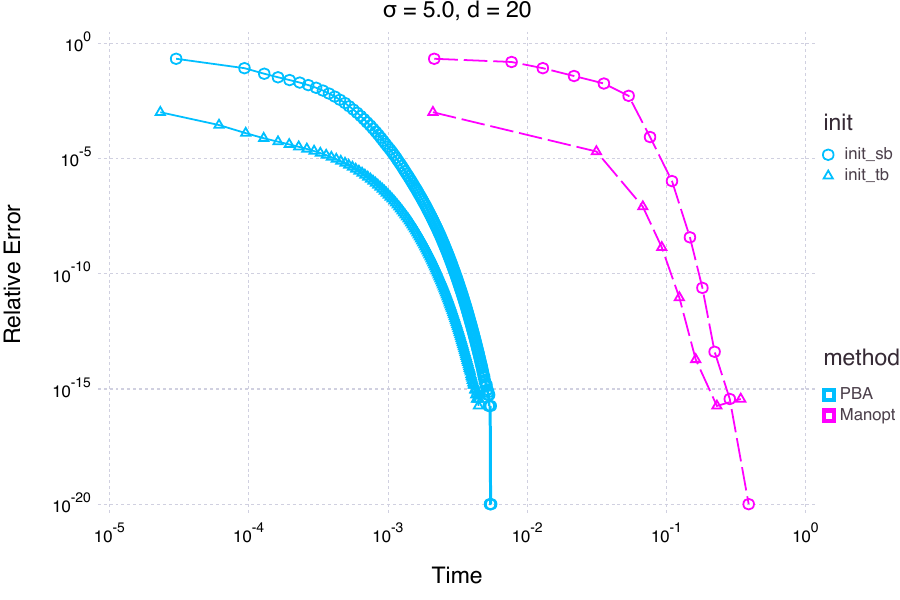}
&
\includegraphics[height=0.18\textheight]{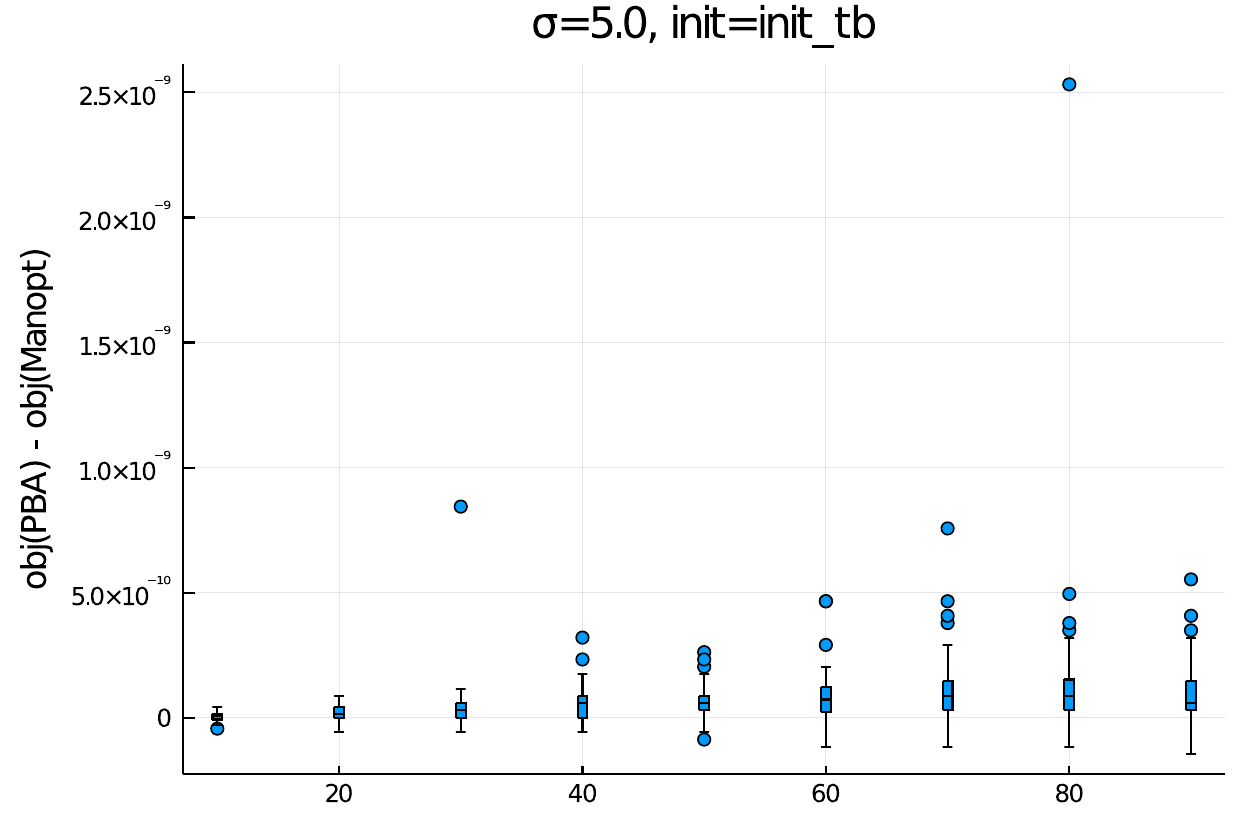}
\\
\includegraphics[height=0.18\textheight]{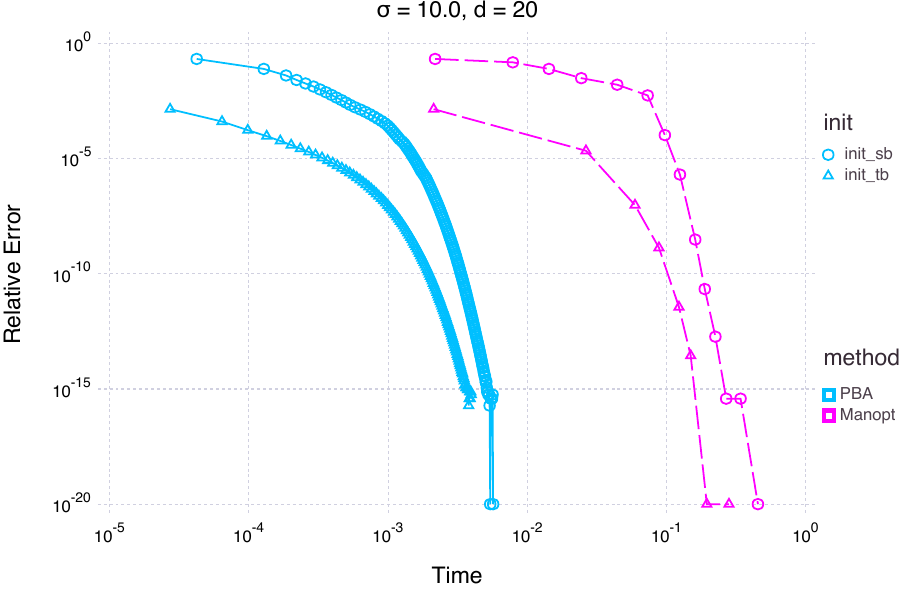}
&
\includegraphics[height=0.18\textheight]{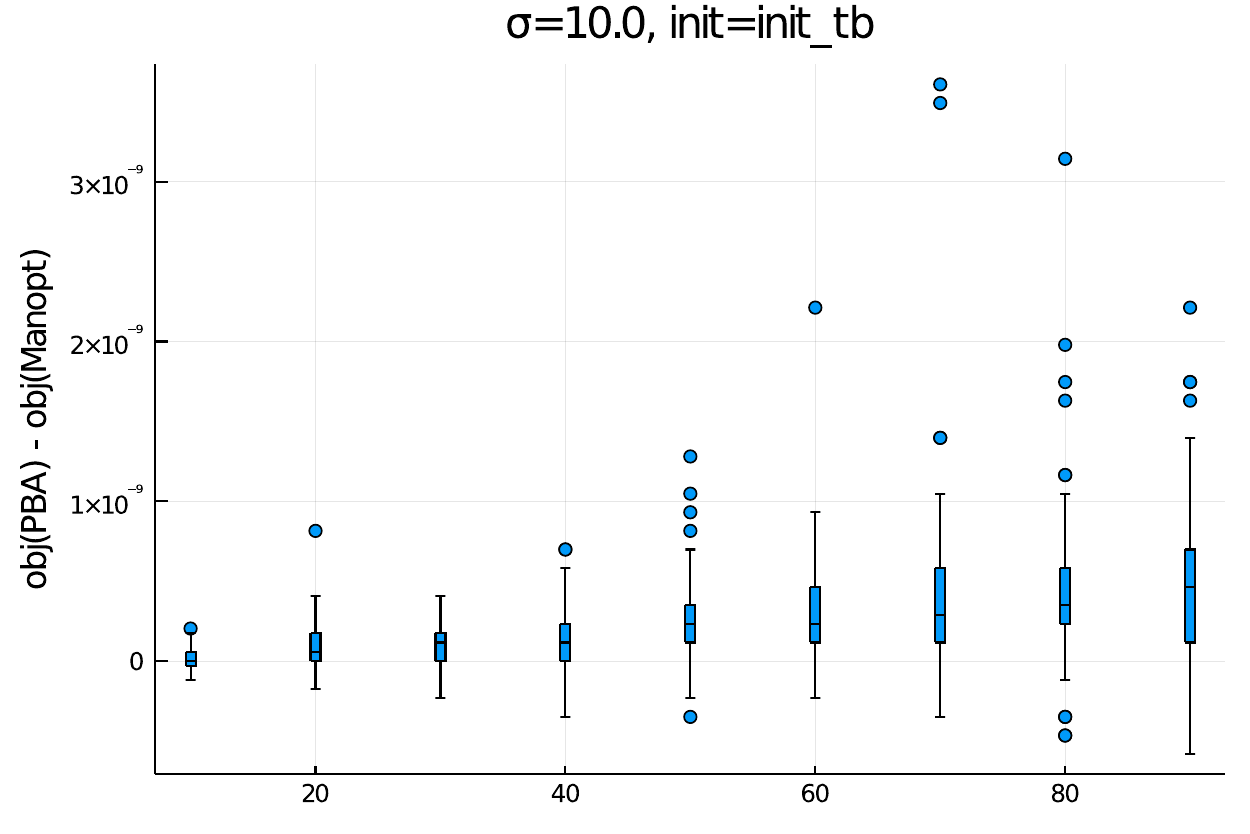}
\end{tabular}
\end{center}
\caption{Simulation studies. Left: relative error vs. wall clock time for each method and initialization strategy. Right: objective value difference between Algorithm \ref{alg:blockrelax} and Manopt at convergence.}\label{fig:summary_r3m5}
\end{figure}

Following the orthogonal Procrustes analysis model, we generated $n$ sets of $d$ dimensional landmarks from the standard normal distribution independently, and randomly rotated them by $m$ orthogonal matrices of size $d\times r$. For this set of $n\times d$ matrices, 
normal error with variance $\sigma^2$ to obtain $\bA_i$ was added, $i=1,\dotsc,m$. 
The data matrix 
$\tilde{\bS}=(\bS_{ij})$ was constructed 
with $\bS_{ij}=\bA_i^T\bA_j$, $i\neq j$, and $\bS_{ii}=\bzero$.
Values of $m=5$, $n=100$, and 
$d \in \{10, 20, \dotsc, 100\}$ were used.  
The noise levels considered were $\sigma\in\{0.1, 1.0, 5.0, 10.0\}$. 
The rank $r$ was set to $3$. 
Initial value strategies ``sb'' and ``tb'' were used for both Algorithm \ref{alg:blockrelax} and Manopt, as they showed good performance in the small examples. 
One hundred samples of random sets were generated for each combination of simulation parameters.
In Fig. \ref{fig:summary_r3m5}, error-vs.-time curves are plotted for typical instances.
Algorithm \ref{alg:blockrelax} were more than an order of magnitudes faster than Manopt. (
For $d=100$,
Manopt did not terminate for more than three days, hence the results were omitted.) 
There were little difference between the two initial value strategies,
hence the differences of the final objective values between Algorithm \ref{alg:blockrelax} and Manopt are plotted in Fig. \ref{fig:summary_r3m5} for strategy ``tb'' only.
The final objective values of the two algorithms agreed in most cases, while Algorithm \ref{alg:blockrelax} tended to give larger values.
The proportions of certified global optima are reported in Tables \ref{tbl:freq_r3} for $d$ that are multiples of tens. 
Not surprisingly, when the noise level was low both Algorithm \ref{alg:blockrelax} and Manopt almost always solved \eqref{eqn:tracemax} globally.
Even if $\sigma$ was as large as $10.0$, the success rate was between 8 to 24\%.
 \subsection{Real-world examples}\label{sec:numerical:real}
\begin{example}[Cryo-EM]\label{ex:cryo}
	\textit{Ab initio} modeling for the single-particle reconstruction (SPR) problem in cryo-EM refers to the procedure of obtaining a preliminary 3D map of the particle in the ice from 2D images by tomographic inversion. Since each cryo-EM image is a noisy projection of the particle with unknown orientation, reliable estimation of orientations from a collection of images is an important step in SPR.
	A popular approach is based on the common-lines property \citep{bendory2020}: the Fourier slice theorem implies that any pair of projection images possesses a pair of radial lines on which there Fourier transforms coincide.
	Once the common lines of all the pairs among $m$ projections, the orientations can be estimated via orthogonal least squares \citep{wang2013orientation}. For a pair of images $i$ and $j$, if the common line between images $i$ and $j$ appears in the direction of $\bc_{ij}=(\cos\theta_{ij}, \sin\theta_{ij}, 0)^T$ in image $i$ and in $\bc_{ji}=(\cos\theta_{ji}, \sin\theta_{ji}, 0)^T$ in image $j$, then the unknown 3D rotation matrices $\bO_i$ and $\bO_j$ in the special orthogonal group $SO(3)$ should approximately satisfy $\bO_i^T\bc_{ij} = \bO_{j}^T\bc_{ji}$. Thus for estimating these matrices for all pairs among the $m$ images, we may minimize
	\[
		\sum_{i<j}\|\bO_i^T\bc_{ij} - \bO_j^T\bc_{ji}\|_F^2
		,
	\]
	which is \eqref{eqn:tracemax} with $d_1=\dotsb=d_m=r=3$, $\bS_{ij} = \bc_{ij}\bc_{ji}^T$ for $i\neq j$ and $\bS_{ii}=\bzero$, $i=1,\dotsc,m$ (i.e., MAXDIFF), but the domain is $SO(3)\times \dotsb \times SO(3)$ instead of $\mathcal{O}_{3,3}\times \dotsb \times \mathcal{O}_{3,3}$.
	Algorithm \ref{alg:blockrelax} can be trivially modified for this setting, since the projection of $\bB=\bP\bD\bQ^T \in \mathbb{R}^{r\times r}$ (full SVD) onto $SO(r)$ is $\bP\diag(1, \dotsc, 1, -1)\bQ^T$, if the singular values of $B$ are sorted in descending order.

	We generated $m$ noisy projections of a ribosomal subunit  provided with the ASPIRE software for SPR\footnote{Available at \url{http://spr.math.princeton.edu/content/download-software}} that implements the orthogonal least squares method via SDP relaxation ($m=100, 500, 1000$). The orientations of the projections were distributed uniformly over $SO(3)$. White Gaussian noise was added to the clean projections to generate noisy images of size 65 by 65 with signal-to-noise ratios (SNR) $\infty$, 1, 1/2, 1/4, 1/8, and 1/16. Common-line pairs were detected with a $1^\circ$ resolution using the functionality of ASPIRE. Due to the presence of noise, the common-line detection rate deteriorates as SNR decreases. Orientations were estimated using two methods, i.e., SDP relaxation of ASPIRE (which utilizes the SDPLR solver\footnote{Available at \url{http://sburer.github.io/files/SDPLR-1.03-beta.zip}}) and Algorithm \ref{alg:blockrelax} initialized with ``sb.'' The mean squared error of the estimated rotation matrices were computed using the formula of \citep[Eq. (8.2)]{wang2013orientation}.
	Because we had difficulties in installing ASPIRE on the Macbook Pro laptop and common-line detection is computationally demanding, all the computations except for Algorithm \ref{alg:blockrelax} were conducted on a Linux workstation with two Intel Xeon E5-2680v2@2.80GHz CPUs (256GB RAM) and eight Nvidia GTX 1080 GPUs (8GB VRAM/GPU).

	The results are collected in Table \ref{tbl:cryoem}.
	Except for the extremely challenging case with a low number of measurements ($m=100$) and SNR (1/16), Algorithm \ref{alg:blockrelax} produced solutions of the same quality as ASPIRE/SDPLR in much shorter time (recall that ASPIRE was run on a much powerful workstation), which, in turn, are certified to be globally optimal except the case $m=1000$ and SNR=1/16.
	Hence these solutions cannot be further improved under the least squares regime.
	In case the two methods disagree, the solution computed by using an SDP relaxation of the orthogonal least squares method failed to be even first-order stationary, possibly due to the rounding procedure of the SDP solution to $SO(3)$, even though the resulting MSE was lower.
\end{example}

\begin{figure}[t!]
\begin{center}
$$
\begin{array}{cc}
\includegraphics[width=0.48\textwidth]{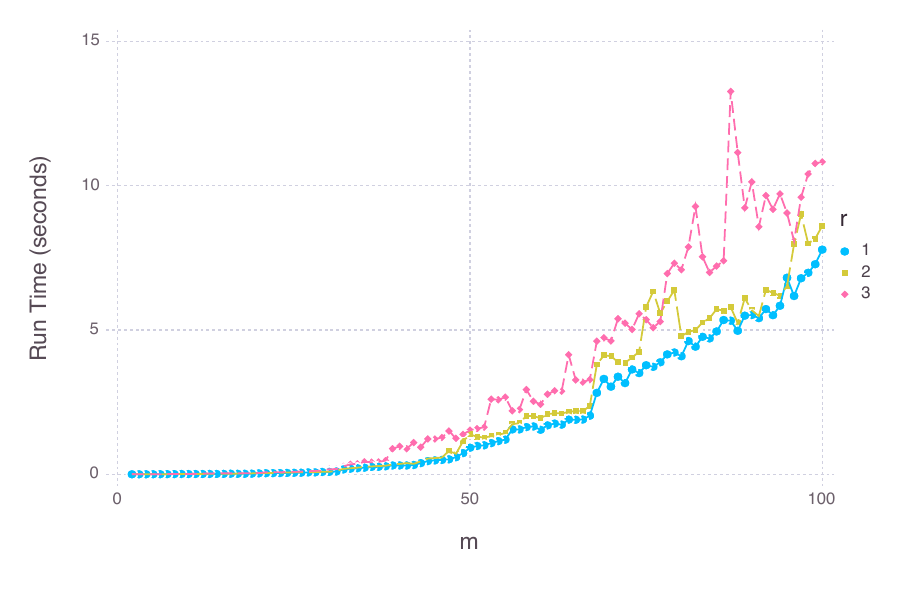} & \includegraphics[width=0.45\textwidth]{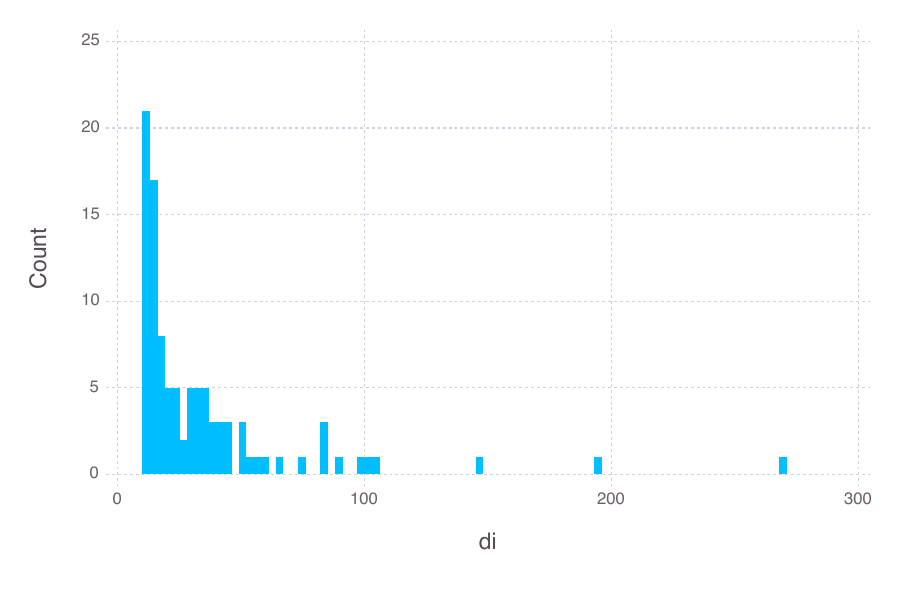}
\end{array}
$$
\end{center}
\vspace{-.4in}
\caption{Left panel displays the run times of the proximal block relaxation algorithm (Algorithm \ref{alg:blockrelax}) for finding the top $r$ canonical correlations among the first $m$ genes on chromosome 1 of $n=488,377$ UK Biobank samples. Right panel shows the distribution of $d_i$ (number of SNPs) in the first 100 genes.}\label{fig:ukb_gcca}
\end{figure}

\begin{example}[Generalized CCA]\label{ex:gCCA}
Our second real-world data example considers gene-level interaction analysis based on genotype data \citep{ZhaoJinXiong06Interaction,PengZhaoXue10GeneLevelInteraction}. Let $\bA_i \in \{0,1,2\}^{n \times d_i}$ be the genotype matrix of gene $i$, where $n$ is the number of individuals and $d_i$ is the number of single nucleotide polymorphisms (SNPs) in gene $i$. To test the interaction between $m$ genes, the maximal canonical correlations among $m$ genes, i.e., \eqref{eqn:tracemax}, were computed. To demonstrate the scalability of the proximal block relaxation algorithm (Algorithm \ref{alg:blockrelax}), we computed the top $r \in \{1,2,3\}$ generalized canonical correlations using the MAXDIFF criterion among the first $m \in \{2,\ldots,100\}$ genes on chromosome 1 of $n=488,377$ samples from the UK Biobank \citep{Sudlow15UKBioBank}. (In contrast, conventional analyses \citep{ZhaoJinXiong06Interaction,PengZhaoXue10GeneLevelInteraction} are restricted to $m=2$ and $r=1$.) 
The numbers of SNPs $d_i$ range from 10 to 271 with mean 34.33 in these genes. Figure \ref{fig:ukb_gcca} displays the run times, all under 15 seconds, of Algorithm \ref{alg:blockrelax} using the same convergence criteria as in Section \ref{sec:numerical:setup}, together with the histogram of $d_i$'s of the 100 genes. Among the 297 local solutions, 107 (36\%) of them were certified to be globally optimal using Theorem \ref{thm:certificate}. 
\end{example}
 \begin{table*}
\begin{center}
\caption{Port Wine Data, MAXDIFF. In ``classification,'' 
``not loc. opt.'' means that the iterate at termination is stationary but confirmed not to be locally maximal, as it violates Proposition \ref{prop:positivity}.
}\label{tbl:portwine_maxdiff}
\medskip
\footnotesize
\begin{tabular}{lllrlllll}
\toprule
$r$ & init & method & iter & time (sec) & obj & classification & $\lambda {\min}(\bL^{\star})$ & $\lambda {\min}(\mathcal{L}^{\star})$ \\
\midrule
\multirow{8}{*}{1} & \multirow{2}{*}{eye} & PBA & 10 & 0.0001923 & 209.8 & global opt. & -1.276e-14 & -1.276e-14 \\
 &  & Manopt & 9 & 0.07502 & 209.8 & global opt. & -2.482e-15 & -2.482e-15 \\
 & \multirow{2}{*}{sb} & PBA & 10 & 0.0002432 & 209.8 & global opt. & -2.101e-14 & -2.101e-14 \\
 &  & Manopt & 7 & 0.05024 & 209.8 & global opt. & \;6.234e-15 & \;6.234e-15 \\
 & \multirow{2}{*}{tb} & PBA & 9 & 0.0001739 & 209.8 & global opt. & -1.159e-14 & -1.159e-14 \\
 &  & Manopt & 5 & 0.03958 & 209.8 & global opt. & -1.951e-14 & -1.951e-14 \\
 & \multirow{2}{*}{lww1} & PBA & 10 & 0.0002291 & 209.8 & global opt. & \;2.599e-15 & \;2.599e-15 \\
 &  & Manopt & 9 & 0.05789 & 209.8 & global opt. & -1.024e-15 & -1.024e-15 \\ 
\cmidrule{2-9} 
\multirow{8}{*}{2} & \multirow{2}{*}{eye} & PBA & 10 & 0.0002719 & 271.2 & global opt. & -1.471e-14 & -79.61 \\
 &  & Manopt & 11 & 0.1084 & 271.2 & global opt. & -2.707e-15 & -79.61 \\
 & \multirow{2}{*}{sb} & PBA & 9 & 0.0002301 & 271.2 & global opt. & -6.812e-14 & -79.61 \\
 &  & Manopt & 7 & 0.08470 & 271.2 & global opt. & -4.465e-14 & -79.61 \\
 & \multirow{2}{*}{tb} & PBA & 9 & 0.0002238 & 271.2 & global opt. & -2.418e-14 & -79.61 \\
 &  & Manopt & 6 & 0.07996 & 271.2 & global opt. & -1.992e-14 & -79.61 \\
 & \multirow{2}{*}{lww1} & PBA & 10 & 0.0002762 & 271.2 & global opt. & -1.364e-13 & -79.61 \\
 &  & Manopt & 12 & 0.1024 & 271.2 & global opt. & -1.219e-14 & -79.61 \\
\cmidrule{2-9} 
\multirow{8}{*}{3} & \multirow{2}{*}{eye} & PBA & 13 & 0.0002659 & 284.1 & global opt. & -3.292e-14 & -106.1 \\
 &  & Manopt & 13 & 0.1814 & 280.7 & not loc. opt. & \text{--} & \text{--} \\
 & \multirow{2}{*}{sb} & PBA & 9 & 0.0001853 & 284.1 & global opt. & -1.370e-14 & -106.1 \\
 &  & Manopt & 6 & 0.1135 & 284.1 & global opt. & -5.021e-14 & -106.1 \\
 & \multirow{2}{*}{tb} & PBA & 10 & 0.0002105 & 284.1 & global opt. & -1.16e-13 & -106.1 \\
 &  & Manopt & 6 & 0.1078 & 284.1 & global opt. & -4.389e-15 & -106.1 \\
 & \multirow{2}{*}{lww1} & PBA & 11 & 0.0002893 & 284.1 & global opt. & \;4.392e-15 & -106.1 \\
 &  & Manopt & 12 & 0.1530 & 280.7 & not loc. opt. 
 & \text{--} & \text{--}
 \\
\bottomrule 
\end{tabular}
\end{center}
\end{table*}

\begin{table*}
\begin{center}
\caption{Example \ref{ex:hard2}. 
In ``classification,'' ``stationary'' means that the iterate at termination is stationary but its global optimality is not confirmed by using Theorem \ref{thm:certificate}.
}\label{tbl:tenberge}
\medskip
\footnotesize
\begin{tabular}{lllrlllll}
\toprule
$r$ & init & method & iter & time (sec) & obj & classification & $\lambda {\min}(\bL^{\star})$ & $\lambda {\min}(\mathcal{L}^{\star})$ \\
\midrule
\multirow{8}{*}{1} & \multirow{2}{*}{eye} & PBA & 2 & 0.0005041 & 1.000 & stationary & -1.000 & -1.000 \\
 &  & Manopt & 1 & 0.001597 & 1.000 & stationary & -1.000 & -1.000 \\
 & \multirow{2}{*}{sb} & PBA & 12 & 0.0001631 & 1.500 & global opt. & -1.608e-10 & -1.608e-10 \\
 &  & Manopt & 7 & 0.02612 & 1.500 & global opt. & -1.045e-7 & -1.045e-7 \\
 & \multirow{2}{*}{tb} & PBA & 2 & 0.0001112 & 1.000 & stationary & -1.000 & -1.000 \\
 &  & Manopt & 1 & 0.001039 & 1.000 & stationary & -1.000 & -1.000 \\
 & \multirow{2}{*}{lww1} & PBA & 2 & 0.0001182 & 1.000 & stationary & -1.000 & -1.000 \\
 &  & Manopt & 1 & 0.0009643 & 1.000 & stationary & -1.000 & -1.000 \\
\cmidrule{2-9} 
\multirow{8}{*}{2} & \multirow{2}{*}{eye} & PBA & 2 & 0.0001587 & 2.000 & stationary & -1.000 & -1.000 \\
 &  & Manopt & 1 & 0.001062 & 2.000 & stationary & -1.000 & -1.000 \\
 & \multirow{2}{*}{sb} & PBA & 50000 & 0.5982 & 3.000 & global opt. & -6.674e-6 & -6.674e-6 \\
 &  & Manopt & 22 & 0.2152 & 3.000 & global opt. & -1.241e-6 & -1.183e-6 \\
 & \multirow{2}{*}{tb} & PBA & 50000 & 0.4910 & 3.000 & global opt. & -6.674e-6 & -6.674e-6 \\
 &  & Manopt & 22 & 0.1925 & 3.000 & global opt. & -1.322e-6 & -1.322e-6 \\
 & \multirow{2}{*}{lww1} &  PBA & 2 & 0.0001049 & 2.000 & stationary & -1.000 & -1.000 \\
 &  & Manopt & 1 & 0.0009794 & 2.000 & stationary & -1.000 & -1.000 \\
\cmidrule{2-9} 
\multirow{8}{*}{3} & \multirow{2}{*}{eye} & PBA & 2 & 6.970e-5 & 3.000 & stationary & -1.000 & -1.000 \\
 &  & Manopt & 1 & 0.001210 & 3.000 & stationary & -1.000 & -1.000 \\
 & \multirow{2}{*}{sb} & PBA & 14 & 0.0002474 & 4.000 & stationary & -1.000 & -1.000 \\
 &  & Manopt & 8 & 0.03438 & 4.000 & stationary & -1.000 & -1.000 \\
 & \multirow{2}{*}{tb} & PBA & 11 & 0.0001370 & 4.000 & stationary & -1.000 & -1.000 \\
 &  & Manopt & 7 & 0.02888 & 4.000 & stationary & -1.000 & -1.000 \\
 & \multirow{2}{*}{lww1} & PBA & 2 & 5.438e-5 & 3.000 & stationary & -1.000 & -1.000 \\
 &  & Manopt & 1 & 0.0009759 & 3.000 & stationary & -1.000 & -1.000 \\
\bottomrule
\end{tabular}
\end{center}
\end{table*}

\begin{table}
\caption{Frequency of certified global optimality, large examples}\label{tbl:freq_r3}
\medskip
\footnotesize

\begin{minipage}{.48\textwidth}
\begin{tabular}{rclrrrr}
\toprule
\multirow{2}{*}[-2pt]{$d$} & \multirow{2}{*}[-2pt]{init} & \multirow{2}{*}[-2pt]{method} & \multicolumn{4}{c}{certification rate ($\sigma$)} \\
\cmidrule{4-7}
& & & $0.1$ & $1.0$ & $5.0$ & $10$ \\
\midrule
\multirow[t]{4}{*}{10} & \multirow[t]{2}{*}{sb} & \text{PBA} & 1.0 & .94 & .19 & .17 \\
 &  & \text{Manopt} & 1.0 & .94 & .19 & .17 \\
 & \multirow[t]{2}{*}{tb} & \text{PBA} & 1.0 & .94 & .19 & .17 \\
 &  & \text{Manopt} & 1.0 & .94 & .19 & .17 \\
\multirow[t]{4}{*}{20} & \multirow[t]{2}{*}{sb} & \text{PBA} & 1.0 & .88 & .14 & .17 \\
 &  & \text{Manopt} & 1.0 & .88 & .14 & .17 \\
 & \multirow[t]{2}{*}{tb} & \text{PBA} & 1.0 & .88 & .14 & .17 \\
 &  & \text{Manopt} & 1.0 & .88 & .14 & .17 \\
\multirow[t]{4}{*}{30} & \multirow[t]{2}{*}{sb} & \text{PBA} & 1.0 & .91 & .12 & .12 \\
 &  & \text{Manopt} & 1.0 & .91 & .12 & .12 \\
 & \multirow[t]{2}{*}{tb} & \text{PBA} & 1.0 & .91 & .12 & .12 \\
 &  & \text{Manopt} & 1.0 & .91 & .12 & .12 \\
\multirow[t]{4}{*}{40} & \multirow[t]{2}{*}{sb} & \text{PBA} & 1.0 & .86 & .10 & .19 \\
 &  & \text{Manopt} & 1.0 & .86 & .10 & .19 \\
 & \multirow[t]{2}{*}{tb} & \text{PBA} & 1.0 & .86 & .10 & .19 \\
 &  & \text{Manopt} & 1.0 & .86 & .10 & .19 \\
\multirow[t]{4}{*}{50} & \multirow[t]{2}{*}{sb} & \text{PBA} & 1.0 & .86 & .22 & .21 \\
 &  & \text{Manopt} & 1.0 & .86 & .22 & .21 \\
 & \multirow[t]{2}{*}{tb} & \text{PBA} & 1.0 & .86 & .22 & .21 \\
 &  & \text{Manopt} & 1.0 & .86 & .22 & .21 \\
\bottomrule
\end{tabular}
\end{minipage} 
\hfill
\begin{minipage}{.50\textwidth}
\begin{tabular}{rclrrrr}
\toprule
\multirow{2}{*}[-2pt]{$d$} & \multirow{2}{*}[-2pt]{init} & \multirow{2}{*}[-2pt]{method} & \multicolumn{4}{c}{certification rate ($\sigma$)} \\
\cmidrule{4-7}
& & & $0.1$ & $1.0$ & $5.0$ & $10$ \\
\midrule
\multirow[t]{4}{*}{60} & \multirow[t]{2}{*}{sb} & \text{PBA} & 1.0 & .78 & .21 & .20 \\
 &  & \text{Manopt} & 1.0 & .78 & .21 & .20 \\
 & \multirow[t]{2}{*}{tb} & \text{PBA} & 1.0 & .78 & .21 & .20 \\
 &  & \text{Manopt} & 1.0 & .78 & .21 & .20 \\
\multirow[t]{4}{*}{70} & \multirow[t]{2}{*}{sb} & \text{PBA} & 1.0 & .79 & .14 & .15 \\
 &  & \text{Manopt} & 1.0 & .79 & .14 & .15 \\
 & \multirow[t]{2}{*}{tb} & \text{PBA} & 1.0 & .79 & .14 & .15 \\
 &  & \text{Manopt} & 1.0 & .79 & .14 & .15 \\
\multirow[t]{4}{*}{80} & \multirow[t]{2}{*}{sb} & \text{PBA} & 1.0 & .91 & .18 & .14 \\
 &  & \text{Manopt} & 1.0 & .91 & .18 & .14 \\
 & \multirow[t]{2}{*}{tb} & \text{PBA} & 1.0 & .91 & .18 & .14 \\
 &  & \text{Manopt} & 1.0 & .91 & .18 & .14 \\
\multirow[t]{4}{*}{90} & \multirow[t]{2}{*}{sb} & \text{PBA} & 1.0 & .80 & .21 & .23 \\
 &  & \text{Manopt} & 1.0 & .80 & .21 & .23 \\
 & \multirow[t]{2}{*}{tb} & \text{PBA} & 1.0 & .80 & .21 & .23 \\
 &  & \text{Manopt} & 1.0 & .80 & .21 & .23 \\
\multirow[t]{4}{*}{100} & \multirow[t]{2}{*}{sb} & \text{PBA} & 1.0 & .76 & .23 & .21 \\
 &  & \text{Manopt} & - & - & - & - \\
 & \multirow[t]{2}{*}{tb} & \text{PBA} & 1.0 & .76 & .23 & .21 \\
 &  & \text{Manopt} & - & - & - & - \\
 \bottomrule
\end{tabular}
\end{minipage}
\end{table}

 \begin{table}
\caption{Example \ref{ex:cryo}. 
``CL rate'' refers to the common-line detection rate.
In ``classification,'' ``stationary'' means that the iterate at termination is stationary but its global optimality is not confirmed by using Theorem \ref{thm:certificate}; ``not stat.'' means that it does not satisfy the first-order local optimality condition \eqref{eqn:eigen1storder}.
}\label{tbl:cryoem}
\medskip
\footnotesize
\begin{center}
\begin{tabular}{rcrlrrrl}
\toprule
$m$ & SNR & CL rate & method & MSE & time (sec) & iters & classification
\\
\midrule
\multirow[t]{10}{*}{100} & \multirow[t]{2}{*}{$\infty$} & \multirow[t]{2}{*}{0.9990} & \text{PBA} & 0.0000 & 0.0582 & 14 &  global opt. \\
 &  &  & \text{SDPLR} & 0.0000 & 1.3764 & - &  global opt. \\
 & \multirow[t]{2}{*}{1} & \multirow[t]{2}{*}{0.9442} & \text{PBA} & 0.0010 & 0.0309 & 13 & global opt. \\
 &  &  & \text{SDPLR} & 0.0010 & 0.8879 & - & global opt. \\
 & \multirow[t]{2}{*}{1/2} & \multirow[t]{2}{*}{0.8275} & \text{PBA} & 0.0148 & 0.0392 & 13 & global opt. \\
 &  &  & \text{SDPLR} & 0.0148 & 0.8510 & - & global opt. \\
 & \multirow[t]{2}{*}{1/4} & \multirow[t]{2}{*}{0.6048} & \text{PBA} & 0.1253 & 0.0333 & 15 & global opt. \\
 &  &  & \text{SDPLR} & 0.1253 & 0.8089 & - & global opt. \\
 & \multirow[t]{2}{*}{1/8} & \multirow[t]{2}{*}{0.3628} & \text{PBA} & 0.6646 & 0.0635 & 29 & global opt. \\
 &  &  & \text{SDPLR} & 0.6646 & 0.9883 & - & global opt. \\
 & \multirow[t]{2}{*}{1/16} & \multirow[t]{2}{*}{0.1834} & \text{PBA} & 2.1363 & 0.2265 & 110 & stationary \\
 &  &  & \text{SDPLR} & 1.7991 & 1.3070 & - & not stat. \\
\multirow[t]{10}{*}{500} & \multirow[t]{2}{*}{$\infty$} & \multirow[t]{2}{*}{0.9994} & \text{PBA} & 0.0000 & 0.6996 & 14 &  global opt. \\
 &  &  & \text{SDPLR} & 0.0000 & 2.9501 & - &  global opt. \\
 & \multirow[t]{2}{*}{1} & \multirow[t]{2}{*}{0.8899} & \text{PBA} & 0.0038 & 0.6890 & 12 & global opt. \\
 &  &  & \text{SDPLR} & 0.0038 & 2.3761 & - & global opt. \\
 & \multirow[t]{2}{*}{1/2} & \multirow[t]{2}{*}{0.7250} & \text{PBA} & 0.0338 & 0.7755 & 14 & global opt. \\
 &  &  & \text{SDPLR} & 0.0338 & 3.1371 & - & global opt. \\
 & \multirow[t]{2}{*}{1/4} & \multirow[t]{2}{*}{0.4860} & \text{PBA} & 0.1959 & 0.9834 & 16 & global opt. \\
 &  &  & \text{SDPLR} & 0.1959 & 5.4611 & - & global opt. \\
 & \multirow[t]{2}{*}{1/8} & \multirow[t]{2}{*}{0.2678} & \text{PBA} & 0.7366 & 1.0938 & 20 & global opt. \\
 &  &  & \text{SDPLR} & 0.7366 & 13.1841 & - & global opt. \\
 & \multirow[t]{2}{*}{1/16} & \multirow[t]{2}{*}{0.1263} & \text{PBA} & 1.6252 & 1.5513 & 30 & global opt. \\
 &  &  & \text{SDPLR} & 1.6252 & 11.2194 & - & global opt. \\
\multirow[t]{10}{*}{1000} & \multirow[t]{2}{*}{$\infty$} & \multirow[t]{2}{*}{0.9994} & \text{PBA} & 0.0000 & 2.5524 & 13 &  global opt. \\
 &  &  & \text{SDPLR} & 0.0000 & 10.1951 & - &  global opt. \\
 & \multirow[t]{2}{*}{1} & \multirow[t]{2}{*}{0.9177} & \text{PBA} & 0.0017 & 2.5475 & 13 & global opt. \\
 &  &  & \text{SDPLR} & 0.0017 & 9.7137 & - & global opt. \\
 & \multirow[t]{2}{*}{1/2} & \multirow[t]{2}{*}{0.7889} & \text{PBA} & 0.0188 & 2.4491 & 13 & global opt. \\
 &  &  & \text{SDPLR} & 0.0188 & 19.7047 & - & global opt. \\
 & \multirow[t]{2}{*}{1/4} & \multirow[t]{2}{*}{0.5686} & \text{PBA} & 0.1297 & 2.6885 & 14 & global opt. \\
 &  &  & \text{SDPLR} & 0.1297 & 34.6575 & - & global opt. \\
 & \multirow[t]{2}{*}{1/8} & \multirow[t]{2}{*}{0.3365} & \text{PBA} & 0.5384 & 3.8236 & 20 & global opt. \\
 &  &  & \text{SDPLR} & 0.5384 & 67.5752 & - & global opt. \\
 & \multirow[t]{2}{*}{1/16} & \multirow[t]{2}{*}{0.1687} & \text{PBA} & 1.3403 & 6.6967 & 35 & stationary \\
 &  &  & \text{SDPLR} & 1.3403 & 102.6977 & - & stationary \\
\bottomrule
\end{tabular}
\end{center}
\end{table}
 \section{Conclusion}
We have presented an in-depth analysis of the orthogonal trace-sum maximization (OTSM) problem, which subsumes various linear and quadratic optimization problems on a product of Stiefel manifolds. In a close analogy with classical results on eigenvalue optimization, a fairly general condition for certifying global optimality of a stationary point of the problem is derived. A practical algorithm to reach a stationary point  with a global convergence guarantee is also proposed. We believe both are new to the literature. Numerical experiments show that the combination of our algorithm and certificate, with initial value strategies ``sb'' and ``tb'', can reveal global optima of various instances of OTSM. 
A further analysis on the probability of global optima of the algorithm, under some distributional assumption on the data, is warranted.

\appendix
\section{Technical proofs}\label{sec:proofs}

\subsection*{Proof of Proposition \ref{prop:positivity}}
	Let $\tau_i=\lambda_{\min}(\bLambda_i)$, the smallest eigenvalue of the symmetric matrix $\bLambda_i$, and $\bv_i\in\mathbb{R}^r$ be the associated unit eigenvector, i.e., $\bLambda_i\bv_i = \tau_i\bv_i$.
	If $\bO_i^{\perp}\in\mathcal{O}_{d_i-r,r}$ fills out $\bO_i$ to a fully orthogonal matrix,
	then for any $\bx_i \in \mathbb{R}^{r}$, $\bW_i=\bO_i^{\perp}\bx_i\bv_i^T$ and $\bW_j = \bzero$ for $j\neq i$ satisfy the tangency condition \eqref{eqn:tangent}. 
	(If $r=d_i$, set $\bO_i^{\perp}=\bzero$.)
	Then, $\tr(\bW^T\tilde{\bS}\bW)=\sum_{i=1}^m\tr(\bW_i^T\bS_{ii}\bW_i)=\sum_{i=1}^m\bx_i^T\bO_i^{\perp T}\bS_{ii}\bO_i^{\perp}\bx_i$, and 
\begin{equation*}
	\tr(\bLambda_i\bW_i^T\bW_i) = \tr(\bLambda_i \bv_i\bx_i^T\bO_i^{\perp T}\bO_i^{\perp}\bx_i\bv_i^T)
	= \tr(\bLambda_i \bv_i\bx_i^T\bx_i\bv_i^T)
	= \|\bx_i\|^2\bv_i^T\bLambda_i\bv_i = \|\bx_i\|^2\tau_i
	.
\end{equation*}
	Further, $\tr(\bLambda_j\bW_j^T\bW_j) = 0$ for $j\neq i$.
	Thus the second-order condition \eqref{eqn:second} entails
\begin{equation*}
	0 \le \sum_{i=1}^m \tr(\bLambda_i\bW_i^T\bW_i) - \tr(\bW^T\bS\bW) 	
	= \sum_{i=1}^m\bx_i^T(\tau_i\bI_r - \bO_i^{\perp T}\bS_{ii}\bO_i^{\perp})\bx_i
	.
\end{equation*}
	Thus $\diag(\tau_1\bI_r, \dotsc, \tau_m\bI_r) \succeq \diag(\bO_1^{\perp T}\bS_{11}\bO_1^{\perp}, \dotsc, \bO_m^{\perp T}\bS_{mm}\bO_m^{\perp})$. 
	Since $\bO_i^{\perp T}\bS_{ii}\bO_i^{\perp}$ is positive semidefinite
	(recall $\bS_{ii}\succeq\bzero$), it follows that
	$\tau_i \ge 0$, $i=1, \dotsc, m$.

\subsection*{Proof of Proposition \ref{prop:sdp}}\label{sec:proofs:sdp}
It suffices to show the constraints $\bO_i\in\mathcal{O}_{d_i,r}$, $i=1,\dotsc,m$, are equivalent to the constraints of problem \eqref{eqn:sdprank}.
From equation \eqref{eqn:blockU}, clearly the former implies the latter.
To show the opposite, first note that $\bU\succeq\bzero$ and $\rank(\bU)=r$ if and only if 
	$m\bU = \bE\bE^T$, 
	$\bE = [	\bE_1^T, \cdots, \bE_m^T ]^T \in \mathbb{R}^{D\times r}$,
for some $\bE_i\in\mathbb{R}^{d_i\times r}$, $i=1,\dotsc,m$.
Then $m\bU_{ii}=\bE_i\bE_i^T \preceq \bI_{d_i}$ and $\tr(m\bU_{ii})=\tr(\bE_i^T\bE_i)=r$ jointly imply that all $r$ singular values of $\bE_i$ are 1. That is, $\bE_i\in\mathcal{O}_{d_i,r}$.
 
\subsection*{Proof of Corollary \ref{cor:nec_maxdiff}}\label{sec:proofs:nec_maxdiff}
Suppose $\bS_{12}=\bS_{21}^T \in \mathbb{R}^{d_1\times d_2}$ has a singular value decomposition 
$\bS_{12}=\bU\bSigma\bV^T$ with $\bU\in\mathcal{O}_{d_1,d}$, $\bV\in\mathcal{O}_{d_2,d}$, and $\bSigma=\diag(\sigma_1, \dotsc, \sigma_d) \allowbreak \in \allowbreak \mathbb{S}^d$, where $d=\min\{d_1, d_2\}$ and $\sigma_1 \ge \dotsb \ge \sigma_d \ge 0$.
Since $\bO_1\bO_2^T \in \mathbb{R}^{d_1\times d_2}$ has $r$ unit singular values and the rest are zero, the von Neumann-Fan inequality entails
\begin{equation*}
    \tr(\bO_1^T\bS_{12}\bO_2) = \tr[(\bO_1\bO_2^T)^T\bS_{12}] \le \sum_{i=1}^r \sigma_i,
\end{equation*}
with equality if and only if $\bO_1 = \bU_1\bR$ (\textit{resp.} $\bO_2=\bV_1\bR$), where $\bR \in \mathcal{O}_{r,r}$ and $\bU_1$ (\textit{resp.} $\bU_2$) consists of the left (\textit{resp.} right) singular vectors of $\bS_{12}$ associated with $\sigma_1, \dotsc, \sigma_r$ in this order.
If we denote such orthogonal matrices by $\bar{\bO}_1$ and $\bar{\bO}_2$, then $\bar{\bO}_1^T\bU = [\bR\; \bzero] = \bar{\bO}_2^T\bV$
and $(\bar{\bO}_1, \bar{\bO}_2)$ is globally optimal stationary point.
The associate Lagrange multipliers are
    $\bar{\bLambda}_1 = \bar{\bO}_1^T\bS_{12}\bar{\bO}_2 = \bar{\bO}_2^T\bS_{21}\bar{\bO}_1 = \bar{\bLambda}_2$
    .
Thus
\begin{equation*}
    \bar{\bLambda} = \bar{\bLambda}_1 = \bar{\bLambda}_2 = \bar{\bO}_1^T\bU\bSigma\bV^T\bar{\bO}_2 = \bR^T\bSigma_1\bR = \bR^T\diag(\sigma_1, \dotsc, \sigma_r)\bR
\end{equation*}
and the diagonal blocks of the certificate matrix $\bL^{\star}$ are
\begin{equation*}
    \bL_{11}^{\star} = \bar{\bO}_1\bar{\bLambda}\bar{\bO}_1^T + \sigma_r\bar{\bO}_1^{\perp}\bar{\bO}_1^{\perp T} = \bU\tilde{\bSigma}\bU^T, 
    \quad
    \bL_{22}^{\star} = \bar{\bO}_2\bar{\bLambda}\bar{\bO}_2^T + \sigma_r\bar{\bO}_2^{\perp}\bar{\bO}_2^{\perp T} = \bV\tilde{\bSigma}\bV^T,
\end{equation*}
where 
$\tilde{\bSigma} = \diag(\bSigma_1, \bI_{d - r}) \succeq \bSigma$.
It follows that
\begin{align*}
    \bL^{\star} 
    = \begin{bmatrix} \bU\tilde{\bSigma}\bU^T & -\bU\bSigma\bV^T \\ -\bV\bSigma\bU^T & \bV\tilde{\bSigma}\bV^T \end{bmatrix}
    \succeq
    \begin{bmatrix} \bU\bSigma\bU^T & -\bU\bSigma\bV^T \\ -\bV\bSigma\bU^T & \bV\bSigma\bV^T \end{bmatrix}
    =
    \begin{bmatrix} \bU \\ -\bV \end{bmatrix}\bSigma
    \begin{bmatrix} \bU \\ -\bV \end{bmatrix}^T
    \succeq \bzero
    .
\end{align*}

\subsection*{Proof of Corollary \ref{cor:nec_data}}\label{sec:proofs:nec_data}
The following lemma extends Theorem 2 of \citet{shapiro1988dual} for $\bS_{ii} \succeq \bzero$, $i=1, \dotsc, m$, and is used to prove the claim of the corollary.
\begin{lemma}\label{lem:shapirobotha}
    If $\bS_{ij}=\bS_{ji}^T \in \mathbb{R}^{d_i\times d_j}$, $i, j = 1, \dotsc, m$, has a singular valude decomposition $\bS_{ij} = \bU_{ij}\bSigma_{ij}\bV_{ij}^T$ with $\bU_{ij}\in\mathcal{O}_{d_i,d}$, $\bV_{ij}\in\mathcal{O}_{d_j,d}$, and $\bSigma_{ij}=\diag(\sigma_1, \dotsc, \sigma_d)\in\mathbb{S}^d$, where $d=\min\{d_i, d_j\}$ and $\sigma_1 \ge \dotsb \ge \sigma_d \ge 0$,
    then for $\bX=\diag(\bX_1, \dotsc, \bX_m)$ with $\bX_i=\sum_{j=1}^m \bU_{ij}\bSigma_{ij}\bU_{ij}^T = \sum_{j=1}^m(\bS_{ij}\bS_{ij}^T)^{1/2}$, the symmetric matrix $\bX - \tilde{\bS}$, 
    where $\tilde{\bS}=(\bS_{ij})$ is the data matrix,
    is positive semidefinite.
\end{lemma}
\begin{proof}
    Since $\bS_{ji} = \bU_{ji}\bSigma_{ji}\bV_{ji}^T = \bS_{ij}^T = \bV_{ij}\bSigma_{ij}\bU_{ij}^T$, we can set $\bV_{ij} = \bU_{ji}$ and $\bSigma_{ij} = \bSigma_{ji}$.
    For any $\by = (\by_1^T, \dotsc, \by_m^T)^T$ with $\by_i \in \mathbb{R}^{d_i}$, let $\ba=\bSigma_{ij}^{1/2}\bU_{ij}^T\by_i$ and $\bb = \bSigma_{ij}^{1/2}\bV_{ij}^T\by_j = \bSigma_{ji}^{1/2}\bU_{ji}^T\by_j$.
    Then from the fact $2\ba^T\bb \le \ba^T\ba + \bb^T\bb$, 
$$
    2\by_i^T\bS_{ij}\by_j = 2\by_i^T\bU_{ij}\bSigma_{ij}\bU_{ji}^T\by_j
    \le
    \by_i^T\bU_{ij}\bSigma_{ij}\bU_{ij}^T\by_i + \by_j^T\bU_{ji}\bSigma_{ji}\bU_{ji}^T\by_j
    .
$$
    Thus,
\begin{align*}
    &\by^T(\bX - \tilde{\bS})\by 
    =
    -\sum_{i,j: j\neq i}\by_i^T\bS_{ij}\by_j + \sum_{i=1}^m\by_i^T(\bX_i - \bS_{ii})\by_i \\
    &\ge
    -\frac{1}{2}\sum_{i=1}^m\by_i^T\bU_{ij}\bSigma_{ij}\bU_{ij}^T\by_i 
    -\frac{1}{2}\sum_{j=1}^m\by_j^T\bU_{ji}\bSigma_{ji}\bU_{ji}^T\by_j     
    +\sum_{i=1}^m\by_i^T\left(\sum_{j=1}^m(\bS_{ij}\bS_{ij}^T)^{1/2} - \bS_{ii}\right)\by_i
    \\
    &=
    -\sum_{i=1}^m\sum_{j\neq i}\by_i(\bS_{ij}\bS_{ij}^T)^{1/2}\by_i 
    +\sum_{i=1}^m\sum_{j\neq i}\by_i(\bS_{ij}\bS_{ij}^T)^{1/2}\by_i 
    = 0
    .
\end{align*}
\end{proof}

\begin{proof}[Proof of Corollary \ref{cor:nec_data}]
    If $\bS_{ij}$ has a singular value decomposition $\bS_{ij} = \bV_i\bSigma_{ij}\bV_j^T$ where $\bSigma_{ij} = \bSigma_{ji}$ is $r \times r$ nonnegative diagonal, then 
$$
    \bX_i = \sum_{j=1}^m (\bS_{ij}\bS_{ij}^T)^{1/2} = \bV_i\left(\sum_{j=1}^m\bSigma_{ij}\right)\bV_j^T = \bV_i\bSigma_{i\cdot}\bV_i^T
    ,
    \quad
    \bSigma_{i\cdot} = \sum_{j=1}^m \bSigma_{ij}.
$$
    From Lemma \ref{lem:shapirobotha}, $\bX - \tilde{\bS} = \diag(\bX_1, \dotsc, \bX_m) - \tilde{\bS} \succeq \bzero$.
    Furthermore, 
$$
    \sum_{j=1}^m \bS_{ij}\bV_j = \sum_{j=1}^m \bV_i\bSigma_{ij}\bV_j^T\bV_j = \bV_i\bSigma_{i\cdot}
$$
    and $\bSigma_{i\cdot} = \sum_{j=1}^m \bV_i^T\bS_{ij}\bV_j$.
    Thus we can set $\bar{\bO}_i = \bV_i$ and $\bar{\bLambda}_i = \bSigma_{i\cdot}$ in Theorem \ref{thm:certificate}. 
    Let $\tau_i$ be the smallest diagonal entry of $\bSigma_{i\cdot}$.
    Then,
$$
    \diag\left(\bV_1\bSigma_{1\cdot}\bV_1^T + \tau_1\bV_1^{\perp}\bV_1^{\perp T}, \dotsc, \bV_m\bSigma_{m\cdot}\bV_m^T + \tau_m\bV_m^{\perp}\bV_m^{\perp T}\right)
    \succeq
    \diag(\bX_1, \dotsc, \bX_m) = \bX
    ,
$$
    hence $\bL^{\star} \succeq \bX - \tilde{\bS} \succeq \bzero$.
\end{proof}

\bibliographystyle{chicago}
\bibliography{globalprocrustes_arxivupload.bib}

\newpage
\begin{huge}
\noindent	Supplementary Materials
\end{huge}

\renewcommand{\thesection}{SM\arabic{section}}
\setcounter{section}{0}
\renewcommand{\thetable}{SM\arabic{section}.\arabic{table}}
\setcounter{table}{0}
\renewcommand{\thefigure}{SM\arabic{figure}}
\setcounter{figure}{0}

\section{Relation to eigenvalue optimization}\label{sec:duality}
\subsection{The dual of orthogonal trace-sum maximization}
The certificate results of \S \ref{sec:global} focus on the parallelism between the second equality of the Ky Fan problem \eqref{eqn:kyfan}, i.e., $\max_{\bO\in\mathcal{O}_{d,r}}\tr(\bO^T\bS\bO) = \max_{\bU\in\mathbb{S}^d} \{\tr(\bS\bU): \bzero \preceq \bU \preceq \bI_d, ~\tr(\bU)=r \}$,
and the convex relaxation \eqref{eqn:sdp_primal} of problem \eqref{eqn:tracemax}.
The first equality 
in equation \eqref{eqn:kyfan} suggests that the dual of the SDP $\max_{\bU\in\mathbb{S}^d} \{\tr(\bS\bU): \bzero \preceq \bU \preceq \bI_d, ~\tr(\bU)=r \}$ is the eigenvalue optimization problem 
\[
\sum_{i=1}^r\lambda_i(\bS) =
\min_{\bX\in\mathbb{S}^d:\;\bX\succeq\bS}
	\sum_{i=1}^r \lambda_i(\bX) 
	,
\]
which can be expressed as SDP \eqref{eqn:eigendual}.
Thus it is natural to ask if the dual SDP \eqref{eqn:sdp_dual} has a parallel eigenvalue-optimization formulation.
In this section we show that such a parallelism indeed holds, with two equivalent formulations.

\begin{theorem}\label{thm:duality}
For problem \eqref{eqn:tracemax}
with a data matrix 
$\tilde{\bS}=(\bS_{ij})$,
the following holds.
\begin{subequations}\label{eqn:duality}
\begin{align}
	\max_{\bO_k\in\mathcal{O}_{d_k,r}, k=1,\dotsc,m} \frac{1}{2}&\sum_{i,j=1}^m\tr(\bO_i^T\bS_{ij}\bO_j)  
	\le 
	\max_{\bU\in\;\mathcal{U}} \frac{m}{2}\tr(\tilde{\bS}\bU),
	\label{eqn:weakduality}
	\\	
	& = 
	\min_{\bX\in\mathcal{X}} \frac{m}{2}\left(
	\sum_{i=1}^r \lambda_i(\tilde{\bS}-\bX) + \frac{1}{m}\sum_{k=1}^m\sum_{i=1}^r \lambda_i(\bX_k) \right) 
	\label{eqn:duality_full:dual} \\
	& = 
	\min_{\bX\in\mathcal{X}, \tilde{\bS}-\bX \; \preceq \; \bzero} 
	\frac{1}{2}\sum_{k=1}^m\sum_{i=1}^r \lambda_i(\bX_k) ,
	\label{eqn:duality_full:dual2} 
\end{align}
\end{subequations}
where 
\[
	\mathcal{U}
		=\{\bU\in\mathbb{S}^{D}: \bU \succeq \bzero, 
		~,\bU_{ii} \preceq \bI_{d_k}, ~ \tr(,\bU_{ii})=r, 
		~i=1,\dotsc,m \}, 
\]
and
\[
	\mathcal{X}=\{\diag(\bX_1,\dotsc,\bX_m) \in \mathbb{S}^D: \bX_k\in\mathbb{S}^{d_k},~k=1,\dotsc,m\}
	.
\]
\end{theorem}
\begin{proof}
    See the next subsection.
\end{proof}

Inequality \eqref{eqn:weakduality} is the relation between \eqref{eqn:tracemax} and its SDP relaxation \eqref{eqn:sdp_primal}, which is proved in \S \ref{sec:global} with a sufficient condition for equality in Theorem \ref{thm:certificate}.
Equalities \eqref{eqn:duality_full:dual} and \eqref{eqn:duality_full:dual2} constitute two equivalent weak dual forms of \eqref{eqn:tracemax}, and a strong dualtity of \eqref{eqn:sdp_primal}.
Clearly relation \eqref{eqn:duality} generalizes Ky Fan's classical result \eqref{eqn:kyfan} on the sum of $r$ largest eigenvalues of a symmetric matrix, and also shows that
\[
\sum_{i=1}^r\lambda_i(\bS) =
\min_{\bX\in\mathbb{S}^d}\left(
	\sum_{i=1}^r \lambda_i(\bS-\bX) + \sum_{i=1}^r \lambda_i(\bX) \right) 
	.
\]

\paragraph{Historical notes}
It is interesting to note that, in the fully orthogonal case ($d_1=\dotsb=d_m=d=r$),
the right-hand side of equation \eqref{eqn:duality_full:dual} reduces to the convex dual for orthogonal Procrustes problem \eqref{eqn:procrustes} studied by 
\ifx\citealt\undefined 
	Shapiro and Botha \cite{shapiro1988dual}:
\else
	\citet{shapiro1988dual}:
\fi
\[
	\max_{\bO_k\in\mathcal{O}_{d,d}, k=1,\dotsc,m} \sum_{i<j}\tr(\bO_i^T\bS_{ij}\bO_j) 
	\le \min_{\bX\in\mathcal{X}} \frac{m}{2}\left(
			\sum_{i=1}^d \lambda_i(\tilde{\bS}-\bX)+\frac{1}{m}\tr(\bX) \right)
    ,			
\]
and the right-hand side of inequality \eqref{eqn:weakduality} (i.e., \ref{eqn:sdp_primal}) reduces to the SDP relaxation of the orthogonal least squares problem \eqref{eqn:aug-LScriterion} studied by 
\ifx\citealt\undefined 
	Zhang and Singer \cite{zhang2017disentangling}:
\else
	\citet{zhang2017disentangling}:
\fi
\begin{equation}\label{eqn:zhangsinger}
	\begin{array}{ll}
	\text{minimize} & \frac{1}{2}(\tr(-\bS\bW) - \sum_{k=1}^m\|\bA_k\|_{\mathrm{F}}^2 )\\
	\text{subject~to} & \bW \succeq 0, 
	                ~
					\bW_{ii} = \bI_d, ~ i=1,\dotsc,m,
	\end{array}
\end{equation}
where 
$\bS = -[\bA_1, \dotsc, \bA_m]^T[\bA_1^T, \dotsc, \bA_m^T]\preceq \bzero$.
The optimization variable is $\bW\in\mathbb{S}^{md}$; $\bW_{ii}$ is the $i$th $d\times d$ diagonal block of $\bS$.
Since for problem \eqref{eqn:aug-LScriterion} $d_1=\dotsb=d_m=d=r$ and the diagonal blocks $-\bA_i^T\bA_i$ of $\bS$ are irrelevant to optimization, by
changing the variable to $\bU=\frac{1}{m}\bW$, we see that \eqref{eqn:zhangsinger} is equivalent to
\begin{equation*}
	\begin{array}{ll}
	\text{maximize} & \frac{m}{2}\tr(\tilde{\bS}\bU) \\
	\text{subject~to} & \bU \succeq 0, 
	                ~
	                m\bU_{ii} = \bI_d, ~ i=1,\dotsc,m,
	\end{array}
\end{equation*}
by setting $\bS_{ij}=-\bA_i^T\bA_j$ if $i\neq j$ and $\bS_{ii}=\bA_i^T\bA_i$.
This in turn is equivalent to \eqref{eqn:sdp_primal}, since under
$d_1=\dotsb=d_m=d=r$ and $\bU\succeq\bzero$, $m\bU_{ii} \preceq \bI_d$ and $\tr(m\bU_{kk})=d$ if and only if $m\bU_{ii}=\bI_d$.

Thus Theorem \ref{thm:duality} generalizes these two results.
In particular,
the strong duality between these three-decade apart results 
appears to be novel.

\begin{remark}\label{rem:dualitygap}
    In the fully orthogonal case,
    a sufficient condition for the duality gap in inequality \eqref{eqn:weakduality} to vanish is 
    $\lambda_d(\tilde\bS-\bX^\star)>\lambda_{d+1}(\tilde\bS-\bX^\star)$ 
    for solution $\bX^\star$ to the problem in the right-hand side of equality \eqref{eqn:duality_full:dual}
    \citep[Theorem 3]{shapiro1988dual}.
    Checking this condition requires to solve the corresponding eigenvalue optimization problem, which is nonsmooth, or equivalently \eqref{eqn:sdp_dual}, which requires lifting the dimension.
    On the other hand, condition \eqref{eqn:certificate} in Theorem \ref{thm:certificate} only requires a 
    stationary point, which is computationally cheap to obtain (see \S\ref{sec:blockascent} for a concrete algorithm).
    Furthermore, Theorem \ref{thm:certificate} implies that the gap between the $d$th and the $(d+1)$st eigenvalues is not essential for the tightness of the upper bound.
\end{remark}

\subsection{Proof of Theorem \ref{thm:duality}}
Inequality \eqref{eqn:weakduality} is already proved in \S\ref{sec:eigen}. Since the right-hand side of this inequality is merely  \eqref{eqn:sdp_primal} and its strong duality with \eqref{eqn:sdp_dual} is also discussed in \S\ref{sec:eigen}, it suffices to show that the right-hand sides of equations \eqref{eqn:duality_full:dual} and \eqref{eqn:duality_full:dual2} are both equivalent to \eqref{eqn:sdp_dual}.

To show equality \eqref{eqn:duality_full:dual}, apply the convex form of the Ky Fan problem (second equality of equation \eqref{eqn:kyfan})
to see that the objective function of the right-hand side of equation \eqref{eqn:duality_full:dual} can be written
\begin{align*}
	& \max_{\bU\in\Phi_{md,r}}\tr((\tilde{\bS}-\bX)\bU) 
	+ \frac{1}{m}\sum_{k=1}^m\max_{\bV_k\in\Phi_{d,r}}\tr(\bX_k \bV_k) \\
	&= \max_{(\bU,\bV_1,\dotsc,\bV_m)\in\;\bar{\mathcal{U}}}
	\tr(\tilde{\bS}\bU) + \sum_{k=1}^m \left( -\tr(\bX_k\bU_{kk}) + \frac{1}{m}\tr(\bX_k\bV_k) \right) \\
	&= \max_{(\bU,\bV_1,\dotsc,\bV_m)\in\;\bar{\mathcal{U}}}
	\tr(\tilde{\bS}\bU) - \sum_{k=1}^m \tr\left(\bX_k(\bU_{kk} - \frac{1}{m}\bV_k) \right),
\end{align*}
where 
$\bar{\mathcal{U}}=\{(\bU,\bV_1,\dotsc,\bV_m): \bU\in\Phi_{D,r},\bV_k\in\Phi_{d_k,r},~k=1,\dotsc,m\}$
with
	$\Phi_{d,r}=\{\bB\in\mathbb{S}^d: \bzero \preceq \bB \preceq \bI_d, ~ \tr(\bB)=r \}$
	,
the convex hull of $\mathcal{O}_{d,r}$ \citep{overton1993optimality,hiriart1995sensitivity}.
Thus the dual problem \eqref{eqn:duality_full:dual} is equivalent to
\[
	\min_{\bX\in\mathcal{X}} \; 
	\max_{(\bU,\bV_1,\dotsc,\bV_m)\in\;\bar{\mathcal{U}}}
	\tr(\tilde{\bS}\bU) - \sum_{k=1}^m \tr\left(\bX_k(\bU_{kk} - \frac{1}{m}\bV_k) \right).
\]
Define a saddle function
\[
	\mathcal{L}(\bX;\bU,\bV) = 
	\tr(\tilde{\bS}\bU) - \sum_{k=1}^m \tr\left(\bX_k(\bU_{kk} - \frac{1}{m}\bV_k) \right),
\]
where $\bV=\diag(\bV_1,\dotsc,\bV_m)$. Because
the objective function of the right-hand side of equation \eqref{eqn:duality_full:dual} is continuous and has finite values on $\mathcal{X}$,
the strong duality
\[
	\max_{(\bU,\bV_1,\dotsc,\bV_m)\in\;\bar{\mathcal{U}}} \;
	\min_{\bX\in\mathcal{X}} \; 
	\mathcal{L}(\bX;\bU,\bV)  
	=
	\min_{\bX\in\mathcal{X}} \; 
	\max_{(\bU,\bV_1,\dotsc,\bV_m)\in\;\bar{\mathcal{U}}} \;
	\mathcal{L}(\bX;\bU,\bV)  
\]
holds \citep[Propositions 5.5.1, 5.5.2]{Bertsekas:ConvexOptimizationTheory:2009}.
The saddle function $\mathcal{L}(\bX;\bU,\bV)$ is linear in $\bX$, hence the left-hand side of the above equality is equivalent to the following SDP problem
\begin{equation}\label{eqn:sdpdual}
	\begin{array}{ll}
	\text{maximize} & \tr(\tilde{\bS}\bU) \\
	\text{subject~to} & \bzero \preceq \bU \preceq \bI_{D}, ~ \tr(\bU) = r, \\
					& \bzero \preceq \bV_i \preceq \bI_{d_i}, ~ \tr(\bV_i)=r, 
					~
					\bU_{ii} = \frac{1}{m}\bV_i,~ i=1,\dotsc,m. 
	\end{array}
\end{equation}
This problem looks almost like \eqref{eqn:sdp_primal}. To establish the complete equivalence, 
we need to remove redundant constraints from \eqref{eqn:sdpdual}. The constraints $\bV_k=m\bU_{kk} \succeq \bzero$ and $\tr(\bU)=r$ can be trivially removed.
To see the constraint $\bU \preceq \bI_{D}$ is redundant, we use the following lemma:
\begin{lemma}[
\ifx\citealt\undefined 
	Bourin and Lee \cite{BourinLee2012unitary}
\else
	\citet{BourinLee2012unitary}
\fi
]\label{lem:decomposition}
For a positive semidefinite block matrix  $\bR \in \mathbb{S}^{p+q}$, there holds
\begin{equation*}
    \bR =
	\begin{bmatrix} \bA & \bC \\ \bC^T & \bB \end{bmatrix}
	= \bP\begin{bmatrix} \bA & \bzero \\ \bzero & \bzero \end{bmatrix}\bP^T
	+ \bQ\begin{bmatrix} \bzero & \bzero \\ \bzero & \bB\end{bmatrix}\bQ^T,
\end{equation*}
for some orthogonal matrices $\bP$, $\bQ \in \mathcal{O}_{p+q,p+q}$,
where $\bA\in\mathbb{S}^p$ and $\bB\in\mathbb{S}^q$.
\end{lemma}
\noindent
Then it follows that $\bU \succeq \bzero$ and $\bU_{ii}=(1/m)\bV_i$ for $i=1,\dots,m$ imply $\|\bU\| \le \sum_{i=1}^m\|(1/m)\bV_i\| \le (1/m)\sum_{i=1}^m\|\bI_{d_i}\| = 1$, where $\|\cdot\|$ is the spectral norm. The latter is equivalent to $\bU \preceq \bI_{D}$.

To establish equality \eqref{eqn:duality_full:dual2},
we show that the eigenvalue minimization problem in the right-hand side of \eqref{eqn:duality_full:dual2}
\begin{equation}\label{eqn:evc}
	\begin{array}{ll}
	\text{minimize}	& \frac{1}{2}\sum_{k=1}^m\sum_{i=1}^r \lambda_i(\bX_k) \\
	\text{subject to} & \tilde\bS - \bX \preceq \bzero, \\
	\end{array}
\end{equation}
where the optimization variable is $\bX = \diag(\bX_1,\dotsc,\bX_m)$, $\bX_k \in \mathbb{S}^{d_k}$, $k=1,\dotsc,m$,
is equivalent to \eqref{eqn:sdp_dual}.
The inner sum of the objective of problem \eqref{eqn:evc} can be written as a minimization form using the dual SDP \eqref{eqn:eigendual} of the convex form of the Ky Fan problem \eqref{eqn:kyfan}. Thus, 
we obtain an equivalent SDP formulation of problem \eqref{eqn:evc}:
\begin{equation}\label{eqn:sdpdual3}
\begin{array}{ll}
\text{minimize} & \frac{m}{2}\sum_{k=1}^m (rz_k + \tr(\bM_k) ) \\
\text{subject to} 
  & \bZ+\bM-\bN-\bL = \tilde\bS, 
  \quad
  \bL \succeq \bzero, \\
  & \bM_k \succeq \bzero, 
  ~
  \bN_k \succeq \bzero, ~ k=1,\dotsc,m, 
\end{array}
\end{equation}
where 
$\bZ=\diag(mz_1\bI_{d_1},\dotsc,mz_m\bI_{d_m})$,
$\bM=\diag(m\bM_1,\dotsc,m\bM_m)$, and
$\bN\allowbreak=\allowbreak\diag(\allowbreak m\bN_1,\allowbreak\dotsc\allowbreak,m\bN_m)$;
$\bX_k=mz_k\bI_{d_k} + m\bM_k - m\bN_k$, $k=1, \dotsc, m$, are eliminated.
Compared to \eqref{eqn:sdp_dual}, SDP \eqref{eqn:sdpdual3} includes an additional block diagonal matrix $\bN$ as an optimization variable. In fact, the dual of SDP \eqref{eqn:sdpdual3} is SDP
\begin{equation}\label{eqn:sdpprimal2}
\begin{array}{ll}
\text{minimize} & \tr(\tilde\bS\bU)  \\
\text{subject to} 
  & \bU \succeq \bzero, 
    \quad
  \bzero \preceq m\bU_{ii} \preceq \bI_d,
    ~
    \tr(m\bU_{ii})=r,~i=1,\dotsc,m,
\end{array}
\end{equation}
which is \eqref{eqn:sdp_primal} with redundant constraints $m\bU_{ii} \succeq \bzero$ added;
the $\bN_k$ are the Lagrange multipliers corresponding to these constraints.
Thus SDP \eqref{eqn:sdpprimal2} is equivalent to \eqref{eqn:sdp_primal}. 
By the strong duality, 
SDP \eqref{eqn:evc} is equivalent to \eqref{eqn:sdp_dual} as desired. 
\qed

\begin{remark}\label{rem:unconstrained}
With regard to \eqref{eqn:sdp_primal}, observe that \eqref{eqn:KKT:Ulower}--\eqref{eqn:KKT:Utrace} in \S\ref{sec:global} imply that $\rank(\bU^\star)\ge r$, which, together with \eqref{eqn:KKT:Lslack} imply that the dimension of the nullspace of $\bL^\star$ in Theorem \ref{thm:certificate} is at least $r$.
Since $-\bL^\star=\tilde\bS-\bX^\star$ (see SDP \eqref{eqn:sdpdual3}) if $\bX^\star$ solves the problem in the right-hand side of equality \eqref{eqn:duality_full:dual},
it follows that $\sum_{i=1}^r\lambda_i(\tilde\bS-\bX^\star)=0$. 
It is now clear that the constrained problem in the right-hand side of \eqref{eqn:duality_full:dual2} also solves the unconstrained problem 
in equality \eqref{eqn:duality_full:dual}.
That the constraint $\tilde{\bS}-\bX \preceq \bzero$ gets rid of the first term in the right-hand side of equation \eqref{eqn:duality_full:dual} is not obvious without the reduction of the constraint $\bU\preceq\bI_{d}$ in SDP \eqref{eqn:sdpdual},
which is possible due to 
Lemma \ref{lem:decomposition}.
\end{remark}

\section{Proof of Theorem \ref{thm:convergence}}\label{app:pf_convergence}
In preparation for the proof, we need the following definition, which can be found in \citet{rockafellar2009variational}:
\begin{definition}[Fr\'{e}chet subdifferentials]
	Vector $\bg$ is a Fr\'{e}chet subgradient of a lower semicontinuous function $\psi:\mathbb{R}^n\to\mathbb{R}\cup\{+\infty\}$ at $\bx \in \dom(\psi)=\{x: \psi(x)<+\infty\}$ if	
\[
	\lim_{\by\to \bx}\inf_{\by\neq\bx} \frac{\psi(\by)-\psi(\bx)-\langle \bg, \by-\bx \rangle}{\|\by-\bx\|} \ge 0,
\]
where $\langle \cdot,\cdot\rangle$ and $\|\cdot\|$ are the standard Euclidean inner product and norm, respectively.
The set of Fr\'{e}chet subgradients of $\psi$ at $\bx$ is called the Fr\'{e}chet subdifferential, and denoted by $\hat{\partial}f(\bx)$. If $\bx\not\in\dom(\psi)$, then $\hat{\partial}f(\bx)=\emptyset$.
The \emph{limiting Fr\'{e}chet subdifferential}, or simply \emph{subdifferential} for short, is defined and denoted by
\[
	\partial \psi(\bx)=\{\bg\in\mathbb{R}^n: \exists \bx_m~\text{such that}~\bx_m\to \bx,~\psi(\bx_m)\to \psi(\bx),~\bg_m \in \hat{\partial}\psi(\bx_m),~\bg_m\to \bg \}.
\]
\end{definition}
\noindent The set $\partial\psi(\bx)$ is closed, convex, and possibly empty.  If $\psi(x)$ is convex, then $\partial\psi(\bx)$ reduces to its convex subdifferential. If $\psi(x)$ is differentiable, then $\partial\psi(\bx)$ reduces to its ordinary differential.

\begin{definition}
	A lower semicontinuous function $\psi$ with $\dom(\psi)\neq\emptyset$ is said to possess the \emph{Kurdyka-\L{}ojasiewicz (KL) property} at a point $\bar{\bx} \in \dom(\partial\psi)$ if there exist $\eta > 0$, a neighborhood $\mathcal{B}_{\rho}(\bar{\bx})\triangleq\{\bx : \|\bx - \bar{\bx}\| < \rho\}$, 
$c>0$, and $\theta\in[0,1)$ such that for any $\bx \in \mathcal{B}_{\rho}(\bar{\bx}) \cap \dom(\partial\psi)$ and $\psi(\bar{\bx}) < \psi(\bx) < \psi(\bar{\bx}) + \eta$, it holds
\begin{equation}\label{eqn:KL}
	|\psi(\bx) - \psi(\bar{\bx})|^{\theta} \le c \dist[\bzero, \partial\psi(\bx)],
\end{equation}
where $\dom(\partial\psi) = \{\bx : \partial\psi(\bx)\neq \emptyset\}$ 
and $\dist[\bzero, \partial\psi(\bx)] = \min\{\|\by\| : \by \in \partial\psi(\bx)\}$.
The quantity $\theta$ is called the \emph{\L{}ojasiewicz exponent}.
The tuple $(\eta,\rho,c,\theta)$ may depend on $\bar{\bx}$.
\end{definition}
\begin{definition}[Closed map]
A set-valued map $M$ from a point in $X$ to a subset of $Y$ is said to be closed at $\bx \in X$ if $\bx^k\to\bx$, $\bx^k\in X$, and $\by^k\to\by$, $\by^k\in M(\bx^k)$ imply $\by \in M(\bx)$.
The set-valued map $M$ is said to be closed on $X$ if it is closed at each point of $X$.
\end{definition}

Algorithm \ref{alg:blockrelax} defines a set-valued map 
$M: \mathcal{O}_{d_1,r}\times\dotsb\times\mathcal{O}_{d_m,r} \to \mathcal{O}_{d_1,r}\times\dotsb\times\mathcal{O}_{d_m,r}: \bTheta=(\bTheta_1, \dotsc, \bTheta_m) \mapsto \bO=(\bO_1, \dotsc, \bO_m)$
recursively with blocks
\begin{equation}\label{eqn:argmin}
        \bO_i \in \argmax_{\bX_i\in\mathcal{O}_{d_i,r}} \tr[(\sum_{j=1}^m\bS_{ij}\bTheta_j)^T(\bX_i-\bTheta_i)] - \frac{1}{2\alpha}\|\bX_i-\bTheta_i\|_{\rm F}^2
,   
\end{equation}
Let us denote the maximand by $g_i(\bX_i\mid\bTheta)$.
Since this function is continuous in both arguments,
if $\bTheta^{k}=(\bTheta_1^{k},\dotsc,\bTheta_m^{k})$ converges to $\bTheta=(\bTheta_1,\dotsc,\bTheta_m)$ and $\bO^{k}=(\bO_1^k,\dotsc,\bO_m^k)$ converges to $\bO=(\bO_1,\dotsc,\bO_m)$, then taking limits in
\[
g_i(\bO_i^k\mid\bTheta^k) \ge g_i(\bX_i\mid\bTheta^k)
\]
for arbitrary 
$\bX_i$ yields
\[
g_i(\bO_i\mid\bTheta) \ge g_i(\bX_i\mid\bTheta)
.
\]
Since the composition of two closed maps with compact domains is closed \citep[Chapter 7]{LuenbergerYe08Book},
it follows that $\bO \in M(\bTheta)$, and the map  is everywhere closed.

Problem \eqref{eqn:tracemax} can be considered an unconstrained minimization problem of a lower semicontinuous function
\begin{equation}\label{eqn:psi}
	\psi(\bO_1,\dotsc,\bO_m) = -f(\bO_1,\dotsc,\bO_m) + \sum_{i=1}^m\delta_i(\bO_i)
\end{equation}
in the vector space $\mathcal{V}=\mathbb{R}^{d_1\times r}\times\dotsb\times\mathbb{R}^{d_m\times r}$,
where $\delta_i$ is the indicator function for the Stiefel manifold $\mathcal{O}_{d_i,r}$,  i.e, $\delta_i(\bA)=0$ if $\bA\in\mathcal{O}_{d_i,r}$ and $\delta_i(\bA)=+\infty$ otherwise. 
Observe that 
$\dom(\psi)\neq\emptyset$ and
$f(\bO)$ is a polynomial in $\mathcal{V}$. It is known that polynomials and indicator functions of Stiefel manifolds possess the KL property; so does a sum of these functions \citep{attouch2010proximal}. Hence $\psi$ possesses the KL property at each point of $\partial\psi(\bO)$ for some \L{}ojasiewicz exponent $\theta$.

In the objective function \eqref{eqn:psi}, since $\partial\delta_i(\bO_i)$ is the orthogonal complement of the tangent space of $\mathcal{O}_{d_i,r}$ at $\bO_i$, 
and a tangent vector has a direct-sum representation $\bO_i\bA_i + \bO_i^{\perp}K_i$ where $\bA_i$ is skew-symmetric and $\bO_i^{\perp}\in\mathcal{O}_{d_i,r - d_i}$ with $\bO_i\bO_i^T+\bO_i^{\perp}\bO_i^{\perp I} = \bI_{d_i}$,
it follows that 
$\partial\delta_i(\bO_i) = \{\bO_i\bLambda_i :\bLambda_i\in\mathbb{S}^r\}$
\citep{Boothby86DiffGeometryBook}. 
Thus
if we denote the subdifferential of $\psi$ with respect to $\bO_i$ by  $\partial_i\psi(\bO_i,\bO_{-i})$, where $\bO_{-i}=(\bO_1,\dotsc,\bO_{i-1},\bO_{i+1},\dotsc,\bO_m)$,
then
$\partial_i\psi(\bO_i,\bO_{-i})=\{-\sum_{j=1}^m \bS_{ij}\bO_j+\bO_i\bLambda_i:\bLambda_i\in\mathbb{S}^r\}$, and
\begin{equation}\label{eqn:Frechet}
	\partial\psi(\bO) 
	= \partial_1\psi(\bO_1,\bO_{-1})\times \dotsb \times \partial_m\psi(\bO_m,\bO_{-m}).
\end{equation}
It follows
that the condition $\bzero\in\partial\psi(\bO)$ is equivalent to the 
stationarity condition \eqref{eqn:firstorder2}.

Now we are ready to prove the claimed result.

\begin{proof}[\Pf of Theorem \ref{thm:convergence}]
	We first show that the sequence of differences $\{\|\bO^k-\bO^{k+1}\|_{\rm F}\}$ is square summable.
	Since each $\mathcal{O}_{d_i,r}$ is compact, $\psi$ is bounded below.
	Also since $\nabla_1 f_i(\bO_i,\bO)=\bS_{ii}\bO_i + \sum_{j\neq i}\bS_{ij}\bO_j$ is the partial derivative of $f(\bO)$ with respect to $\bO_i$ and it is Lipschitz continuous with modulus $\|\bS_{ii}\|_2$, we have
	\begin{equation}\label{eqn:lipshitz}
	    -f(\bO_i) \le -f(\bO^{k}) - \tr[\nabla_1 f_i(\bO_i^k, \bO^{\text{prev}})^T(\bO_i-\bO_i^k)] + \frac{L}{2}\|\bO_i-\bO_i^{k}\|_{\rm F}, 
	    \quad \forall \bO_i \in \mathcal{O}_{d_i,r},
	\end{equation}
	where $L=\max_{i=1,\dotsc,m}\|\bS_{ii}\|_2$.
	At the $i$th block update in the $k+1$st cycle of Algorithm \ref{alg:blockrelax}, there holds
	\[
	\bO_i^{k} \in \argmax_{\bO_i\in\mathcal{O}_{d_i,r}}\tr[\nabla_1 f_i(\bO_i^k, \bO^{\text{prev}})^T(\bO_i - \bO_i^k)] - \frac{1}{2\alpha}\|\bO_i-\bO_i^k\|_{\rm F}^2 + \delta_i(\bO_i)
	.
	\]
	Thus
	\begin{equation}\label{eqn:minimizer}
	\begin{split}
	    -\tr&[\nabla_1 f_i(\bO_i^k, \bO^{\text{prev}})^T (\bO_i^{k+1} - \bO_i^k)]
	    + \frac{1}{2\alpha}\|\bO_i^{k+1}-\bO_i^k\|_{\rm F}^2 
	    + \delta_i(\bO_i^k)
	    \\
	    &\le
	    -\tr[\nabla_1 f_i(\bO_i^k, \bO^{\text{prev}})^T(\bO_i^{k} - \bO_i^k)]
	    + \frac{1}{2\alpha}\|\bO_i^{k}-\bO_i^k\|_{\rm F}^2 
	    + \delta_i(\bO_i^k)	    
	    .
	\end{split}
	\end{equation}
	Substituting $\bO_i$ in inequality \eqref{eqn:lipshitz} with $\bO_i^k$ and adding to inequality \eqref{eqn:minimizer} yields
	\[
	    \left(\frac{1}{2\alpha}-\frac{L}{2}\right)\|\bO_i^{k+1}-\bO^k\|_{\rm F}^2
	    \le
	    [-f(\bO^k) + \delta_i(\bO_i^k)] - [-f(\bO^{k+1}) + \delta_i(\bO_i^k)]
	    .
	\]
	(Recall that $1/(2\alpha)-L/2 > 0$.)
	Summing the above inequality from $i=1$ to $m$, we obtain
	\begin{equation}\label{eqn:squaresummable}
        \left(\frac{1}{2\alpha}-\frac{L}{2}\right)\|\bO^{k+1}-\bO^{k}\|_{\rm F}^2
        =
        \left(\frac{1}{2\alpha}-\frac{L}{2}\right)\sum_{i=1}^m\|\bO_i^{k+1}-\bO_i^{k}\|_{\rm F}^2
        \le
        \psi(\bO^k) - \psi(\bO^{k+1})
	\end{equation}
	for each $k$. 
	This implies
	$\sum_{k=0}^\infty\|\bO^k-\bO^{k+1}\|_{\rm F}^2 \le 
	\left(\frac{1}{2\alpha}-\frac{L}{2}\right)^{-1}[\psi(\bO^0)-\bar{\psi}] < \infty$,
	where $\bar{\psi}=\lim_{k\to\infty}\psi(\bO^k) > -\infty$.

	We now show that all the limit points of the sequence $\{\bO^k=(\bO_1^k,\dotsc,\bO_m^k)\}$ are stationary points.
	The square summability above implies that $\lim_{k\to\infty}\|\bO^k-\bO^{k+1}\|_{\rm F}=0$ and by Ostrowski's theorem 
\ifx\citealt\undefined 
	(see, e.g., \cite[Proposition 12.4.1]{lange2013}), 
\else
	\citep[see, e.g.,][Proposition 12.4.1]{lange2013}, 
\fi
	the set $W$ of limit points of the $\{\bO^k\}$ is compact and connected.
	If a subsequence $\{\bO^{k_j}\}$ of $\{\bO^k\}$ converges to $\hat\bO=(\hat{\bO}_1,\dotsc,\hat{\bO}_m)$, then it follows $\bO^{k_j+1}\to\hat\bO$. From inequality \eqref{eqn:minimizer}, we have
	\begin{align*}
	    -\tr&[\nabla_1 f_i(\bO_i^{k_j}, \bO^{\text{prev}})^T (\bO_i^{k_j+1} - \bO_i^{k_j})]
	    + \frac{1}{2\alpha}\|\bO_i^{k_j+1}-\bO_i^{k_j}\|_{\rm F}^2 
	    + \delta_i(\bO_i^k)
	    \\
	    &\le
	    -\tr[\nabla_1 f_i({\hat\bO}_i, \bO^{\text{prev}})^T(\hat{\bO}_i - \bO_i^k)]
	    + \frac{1}{2\alpha}\|\bO_i^{k}-\bO_i^k\|_{\rm F}^2 
	    + \delta_i(\hat{\bO}_i)	   
	\end{align*}
	with $\bO^{\text{prev}}=(\bO_1^{k_j+1},\dotsc,\bO_{i-1}^{k_j+1},\bO_{i}^{k_j},\bO_{i+1}^{k_j},\dotsc,\bO_m^{k_j})$.
	Take limit superior on both sides of this 
	inequality. 
	Since $\nabla_1 f_i$ is continuous, this yields
	$\limsup_{j\to\infty} \allowbreak \delta_i(\bO_i^{k_j+1}) \allowbreak \le \allowbreak \delta_i(\hat{\bO}_i)$.
	However, $\delta_i$ is lower semicontinuous. Thus 
	$\lim_{j\to\infty} \delta_i(\bO_i^{k_j+1}) = \delta_i(\hat{\bO}_i)$.
	Together with the continuity of $f$, this implies that $\lim_{j\to\infty}\psi(\bO^{k_j+1})=\psi(\hat{\bO})=\lim_{k\to\infty}\psi(\bO^k)\allowbreak=\bar{\psi}$. 
	Thus $\psi$ takes the constant value $\bar{\psi}$ on $W$.
	Finally, closedness of the map $\bO^{k+1} \in M(\bO^k)$ implies $\hat{\bO} \in M(\hat{\bO})$. This in turn implies
	$\bzero \in \partial \psi(\hat{\bO})$, because
		$\hat{\bO}_i \in \argmin_{\bO_i\in\mathbb{R}^{d_i\times r}}
			\left\{ \psi(\bO_i,\hat{\bO}_{-i}) + \frac{1}{2\alpha}\|\bO_i-\hat{\bO}_i\|_{\rm F}^2 \right\}$
		implies
		$\bzero \in \partial_i \psi(\hat{\bO}_i,\hat{\bO}_{-i})$ 
	for $i=1,\dotsc,m$.

	To see that the whole sequence converges, observe that $\{\bO^k\}$ is bounded and thus has a finite limit point $\bO^\star$.
	From the representation of the subdifferential \eqref{eqn:Frechet}, 
	$-\sum_{j=1}^m\bS_{ij}\bO_j^k + \bO_i^k\bLambda_i^{k+1} \in \partial_i\psi(\bO_i^k,\bO_{-i}^k)$,
	where $\bLambda_i^{k+1}\in\mathbb{S}^r$ is the Lagrange multiplier associated with 
	$\bO_i^{k+1}$, the $k+1$st update of $\bO_i$.
	An argument similar to \S\ref{sec:global:local} reveals that $\bLambda^{k+1}$ satisfies
\begin{align}\label{eqn:lambda_k+1}
	\bO_i^{k+1}\bLambda_i^{k+1} &= \sum_{j<i}\bS_{ij}\bO_j^{k+1}+\sum_{j=i}^m\bS_{ij}\bO_j^{k} + {\alpha}^{-1}(\bO_i^{k}-\bO_i^{k+1}), 
\end{align}
which implies
	\begin{align*}
		\sum_{j<i}\bS_{ij}(\bO_j^{k+1}-\bO_j^k) + (\bO_i^k-\bO_i^{k+1})(\bLambda_i^{k+1}+\alpha^{-1}\bI_r) 
		\in \partial_i\psi(\bO_i^k,\bO_{-i}^k)
		.
	\end{align*}
	Let $\hat{\bS}$ be the strictly upper block triangular matrix obtained from 
	$\tilde{\bS}=(\bS_{ij})$.
	It follows
	\[
	\hat{\bS}(\bO^{k+1}-\bO^k)  - (\bO^{k+1}-\bO^k) \diag(\bLambda_1^{k+1}+\alpha^{-1}\bI_r,\dotsc,\bLambda_m^{k+1}+\alpha^{-1}\bI_r)
		\in \partial\psi(\bO^k).
	\]
	Thus
	\begin{align}\label{eqn:distupperbound}
	\dist[\bzero,& \partial\psi(\bO^k)]
	 \\
	&\le 
	\|\hat{\bS} - \diag(\bLambda_1^{k+1}+\alpha^{-1}\bI_r,\dotsc,\bLambda_m^{k+1}+\alpha^{-1}\bI_r)\|_2 \|\bO^{k+1}-\bO^k \|_{\rm F} \nonumber \\
	&\le L \|\bO^{k+1}-\bO^k \|_{\rm F},
	\nonumber
	\end{align}
where $L=\|\hat{\bS}\|_2 + \sqrt{m} \max_{i=1,\dotsc,m}\|(\tilde{\bS}+\alpha^{-1}\bI_D)_{i\cdot}\|_2$;
$(\tilde{\bS}+\alpha^{-1}\bI_D)_{i\cdot}$ denotes the $i$th $d_i\times D$ row block of $\tilde{\bS}+\alpha^{-1}\bI_D$.
To see the last inequality,
first observe that
	\begin{align*}
	\|\hat{\bS} - \diag(\bLambda_1^{k+1}+\alpha^{-1}\bI_r,\dotsc,\bLambda_m^{k+1}+\alpha^{-1}\bI_r)\|_2
	&\le
	\|\hat{\bS}\|_2 + \max_{i=1,\dotsc,m}\|\bLambda_i^{k+1}+\alpha^{-1}\bI_r\|_2  
	.
	\end{align*}
Let 
$\bTheta=\frac{1}{\sqrt{m}}\big[(\bO_1^{k+1})^T,\dotsc,(\bO_{i-1}^{k+1})^T,(\bO_i^{k+1})^T,(\bO_{i+1}^k)^T,\dotsc,(\bO_m^{k})^T\big]^T \in \mathcal{O}_{D,r}$.
Then from equation \eqref{eqn:lambda_k+1},
	\begin{align*}
	\|\bLambda_i^{k+1}&+\alpha^{-1}\bI_r\|_2
	\\
	&= \max_{\bx\neq\bzero}
		{\|(\bO_i^{k+1})^T[\sum_{j\le i}\bS_{ij}\bO_j^{k+1}+\sum_{j>i}\bS_{ij}\bO_j^{k} + {\alpha}^{-1}(\bO_i^{k}-\bO_i^{k+1})]\bx+\alpha^{-1}\bx\|}/{\|\bx\|} \\
	&\le \max_{\bx\neq\bzero}
		{\|(\sum_{j\le i}\bS_{ij}\bO_j^{k+1}+\sum_{j>i}\bS_{ij}\bO_j^{k} + {\alpha}^{-1}\bO_i^{k})\bx\|}/{\|\bx\|} 
	\\
	&= \max_{\bx\neq\bzero}
		{\|\sqrt{m}[\bS_{i1},\dotsc,\bS_{i,i-1},\bS_{ii}+\alpha^{-1}\bI_{d_i},\bS_{i,i+1},\dotsc,\bS_{im}]\bTheta\bx\|}/{\|\bTheta\bx\|} 
	\\
	&\le \sqrt{m}\max_{\by\neq\bzero}
		{\|[\bS_{i1},\dotsc,\bS_{i,i-1},\bS_{ii}+\alpha^{-1}\bI_{d_i},\bS_{i,i+1},\dotsc,\bS_{im}]\by\|}/{\|\by\|} 
	\\
	&= \sqrt{m}\|(\tilde{\bS}+\alpha^{-1}\bI_D)_{i\cdot}\|_2
	,
	\end{align*}
where the first inequality holds because $\bO_i^{k+1} \in \mathcal{O}_{d_i,r}$, and the second inequality is due to that $\|\bTheta\bx\|=\|\bx\|$.

	Now we can exploit the KL property of the problem.
	From inequality \eqref{eqn:distupperbound} and the KL inequality \eqref{eqn:KL}, we have
	\[
		|\psi(\bO^k)-\bar{\psi}|^{\tilde{\theta}}
		= 
		|\psi(\bO^k)-\bar{\psi}-\psi(\bO^\star)+\bar{\psi}|^{\tilde{\theta}}
		=
		|\psi(\bO^k)-\psi(\bO^\star)|^{\tilde{\theta}}
		\le \tilde{c} L \|\bO^{k+1}-\bO^k\|_{\rm F}
	\]
	in the neighborhood of $\bO^\star$ 
	such that $\bO\in\mathcal{B}_{\tilde{\rho}}\cap \dom(\partial\psi)$
	and $\psi(\bO^\star) < \psi(\bO) < \psi(\bO^\star) + \tilde{\eta}$.
	Note $\tilde\eta$, $\tilde\rho$, $\tilde{c}$, and $\tilde{\theta}$ depend on $\bO^\star$.
	In order to get rid of dependency on $\bO^\star$, cover $W$ by a finite number of balls $\mathcal{B}_{\tilde\rho(\bO^{\star j})}$, $\bO^{\star j}\in W$, $j=1,\dotsc,J$, and
	take $\theta=\max_j\tilde\theta(\bO^{\star j})$,
	$c=\max_j\tilde{c}(\bO^{\star j})$, and
	$\eta=\min_j\tilde{\eta}(\bO^{\star j})$.
	For a sufficiently large $K$, every $\bO^k$ with $k\ge K$ falls within these balls and satisfies $|\psi(\bO^k)-\bar{\psi}|<1$.
	Without loss of generality assume $K=0$. 
	The KL inequality now entails
	\begin{equation}\label{eqn:KLbound}
		|\psi(\bO^k)-\bar{\psi}|^{\theta}
		\le c L \|\bO^{k+1}-\bO^k\|_{\rm F}.
	\end{equation}
	If $\psi(\bO^k) = \bar{\psi}$ for some $k$, then since the sequence $\{\psi(\bO^k)\}$ is nonincreasing,
	$\bar{\psi}=\psi(\bO^k)=\psi(\bO^{k+1})=\dotsb$. From inequality \eqref{eqn:squaresummable}, $0=\psi(\bO^k)-\psi(\bO^{k+1}) \ge (2\alpha)^{-1}\|\bO^k-\bO^{k+1}\|_{\rm F}^2$.
	Thus $\bO^k=\bO^{k+1}=\dotsb=\bO^\star$.
	Otherwise we have $\psi(\bO^k) > \bar{\psi}$ for all $k$.
	In combination with the concavity of the function $a^{1-\theta}$ on $[0, \infty)$, inequalities \eqref{eqn:squaresummable} and \eqref{eqn:KLbound} imply
	\begin{align*}
	[\psi(\bO^k)-\bar{\psi}]^{1-\theta} -[\psi(\bO^{k+1})-\bar{\psi}]^{1-\theta}
	& \ge \frac{(1-\theta)[\psi(\bO^k)-\psi(\bO^{k+1})]}{[\psi(\bO^k)-\bar{\psi}]^{\theta}} 
	\\
	&
	\ge \frac{1-\theta}{2\alpha cL}\|\bO^{k+1}-\bO^{k}\|_{\rm F}.
	\end{align*}
	Rearranging this inequality and summing over $k$ yield
	\[
		\sum_{k=0}^{\infty}\|\bO^{k+1}-\bO^k\|_{\rm F} \le
		\frac{2\alpha c L}{1-\theta}[\psi(\bO^0)-\bar{\psi}]^{1-\theta}.
	\]
	Thus, the sequence $\{\bO^k\}$ is Cauchy and hence converges to a unique limit $\bO^\star$.
	
    Finally, the rate of convergence of Algorithm \ref{alg:blockrelax} follows directly from Theorem 3 in \citet{xu2017globally}:
\begin{proposition}
Under the assumptions of Theorem \ref{thm:convergence}, suppose $\lim_{k\to\infty} \allowbreak \bO^k \allowbreak = \allowbreak \bO^\star$. 
Let the \L{}ojasiewicz exponent of \eqref{eqn:tracemax} near $\bO^\star$ be $\theta$.
Then the following holds.
\begin{enumerate}
	\item If $\theta\in[0,\frac{1}{2}]$, then $\|\bO^k-\bO^\star\|_{\rm F}\le Cr^k$, $\forall k$, for some $C>0$ and $r\in[0,1)$;
	\item If $\theta\in(\frac{1}{2},1)$, then $\|\bO^k-\bO^\star\|_{\rm F}\le Ck^{-(1-\theta)/(2\theta-1)}$, $\forall k$, for some $C > 0$.
\end{enumerate}
\end{proposition}
\end{proof}
 \section{Additional experiments}\label{sec:additionalnumerical}
\begin{example}[Example 5.1, \citet{liu2015maximization}]\label{ex:liu5_1}
    For MAXBET,
    the data matrix is
\begin{equation*}
    \tilde{\bS} = \begin{bmatrix} 
     4.3299 & 2.3230 & -1.3711  & -0.0084  & -0.7414 \\
     2.3230 & 3.1181 &  1.0959  &  0.1285  &  0.0727 \\
    -1.3711 & 1.0959 &  6.4920  & -1.9883  & -0.1878 \\
    -0.0084 & 0.1285 & -1.9883  &  2.4591  & 1.8463  \\
    -0.7414 & 0.0727 & -0.1878  &  1.8463  & 5.8875 
    \end{bmatrix}
\end{equation*}
    with $m=2$ and $d_1 = 2$, $d_2 = 3$. 
    Ranks $r=1$ and $2$ were tested. For the former, the globally optimal value is known to be $7.365$ \citep[p. 1500]{liu2015maximization}.
    For MAXDIFF, we set the corresponding diagonal blocks to zero.
    Tables \ref{tbl:liu5_1_maxdiff} and \ref{tbl:liu5_1_maxbet} summarize the results.
    Global optima were found for all cases, which appear to be new to the literature except for MAXBET, $r=1$.
    Again strategies ``sb'' and ``tb'' were more effective than ``eye'' and ``lww1,'' 
    and Algorithm \ref{alg:blockrelax} (PBA) was orders of magnitudes faster than Manopt.
\end{example}

\begin{example}[Example 5.2, \citet{liu2015maximization}]\label{ex:liu5_2}
    For MAXBET,
    the data matrix is
    $$
    \tilde{\bS} = \begin{bmatrix} 
     45 & -20 &   5 &   6 & 16 &   3 \\
    -20 &  77 & -20 & -25 & -8 & -21 \\
      5 & -20 &  74 &  47 & 18 & -32 \\
      6 & -25 &  47 &  54 &  7 & -11 \\
     16 &  -8 &  18 &   7 & 21 &  -7 \\
      3 & -21 & -32 & -11 & -7 &  70
    \end{bmatrix}
    $$
    with $m=3$ and $d_1 = d_2 = d_3 = 2$.
    Ranks $r=1$ and $2$ were tested. 
    For MAXDIFF, we set the corresponding diagonal blocks to zero.
    Tables \ref{tbl:liu5_2_maxdiff} and \ref{tbl:liu5_2_maxbet} summarize the results.
    Global optima were confirmed for all cases except for MAXBET with $r=1$. 
    For the latter, the converged objective value ($189.48$) is believed to be globally optimal \citep{hanafi2003global}.
    For the others, the found global optima appear to be novel.
    There was a case (MAXDIFF, $r=2$) that Manopt converged to a  stationary point that is not a local maximizer if started with the ``eye'' strategy, while Algorithm \ref{alg:blockrelax} correctly found a global optimum.
    In all cases, Algorithm \ref{alg:blockrelax} was orders of magnitudes faster than Manopt.
\end{example}

\begin{figure}[t!]
\begin{center}
\begin{tabular}{ll}
\includegraphics[width=0.45\textwidth]{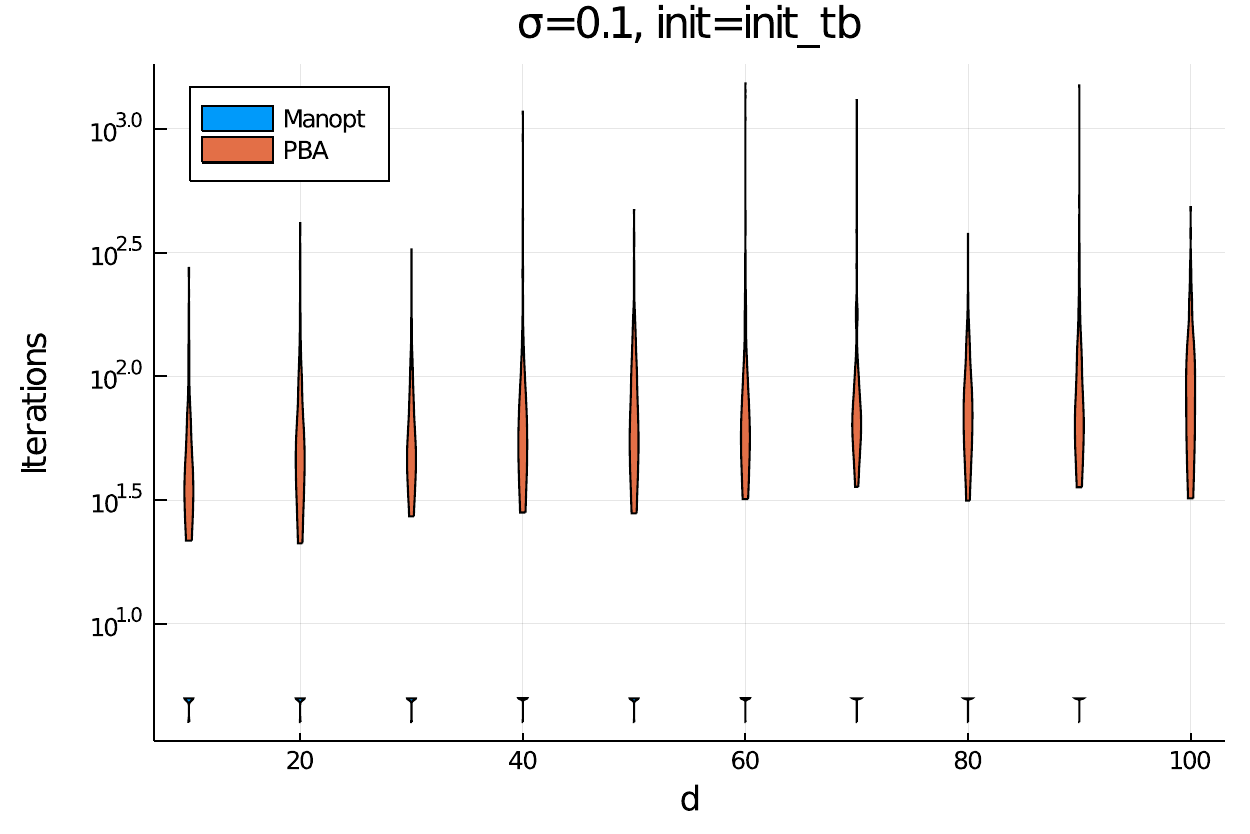}
&
\includegraphics[width=0.45\textwidth]{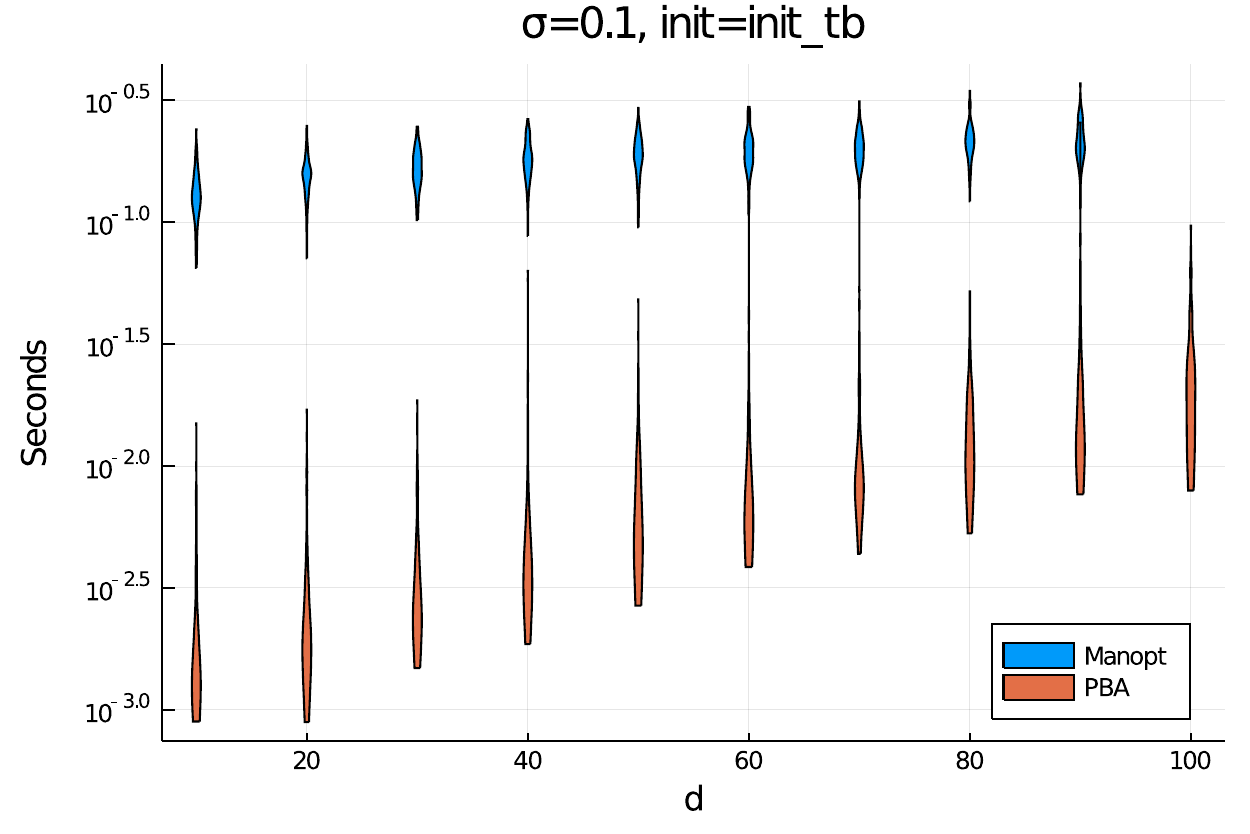}
\\
\includegraphics[width=0.45\textwidth]{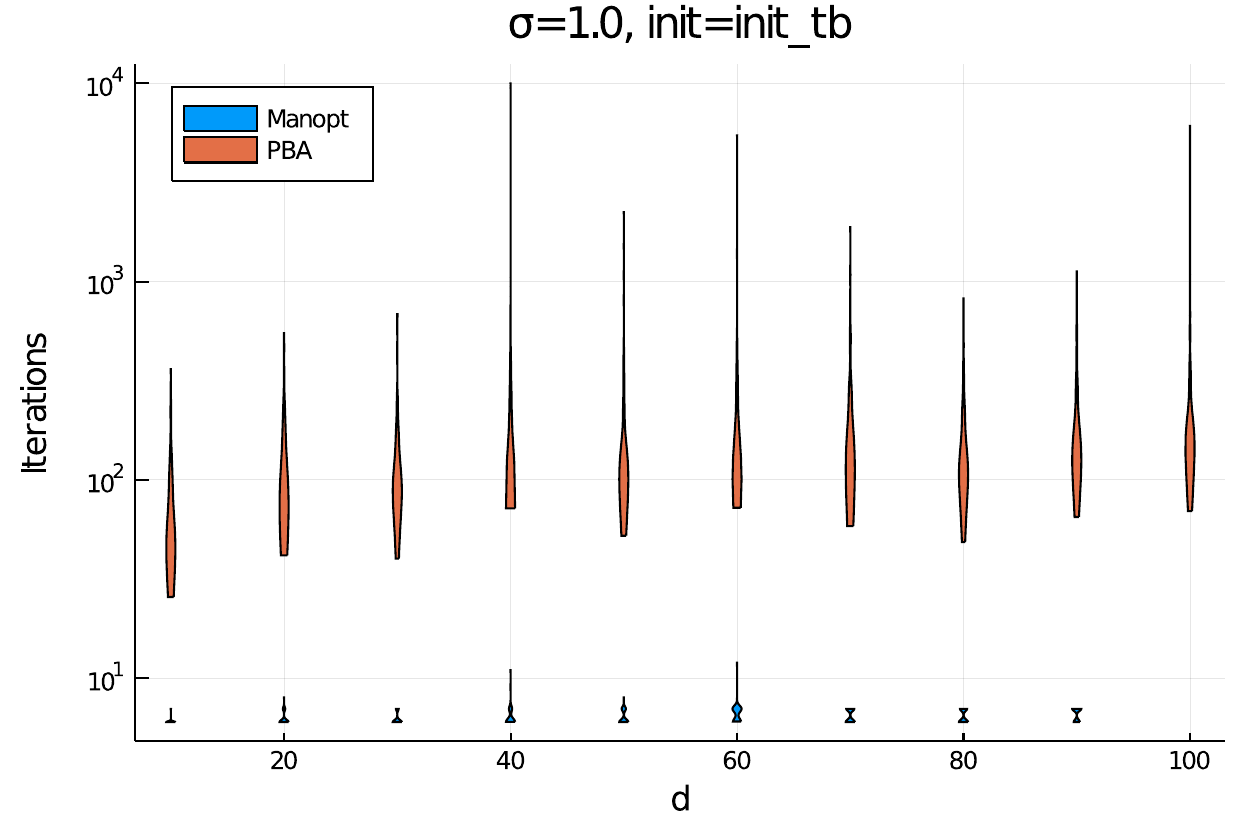}
&
\includegraphics[width=0.45\textwidth]{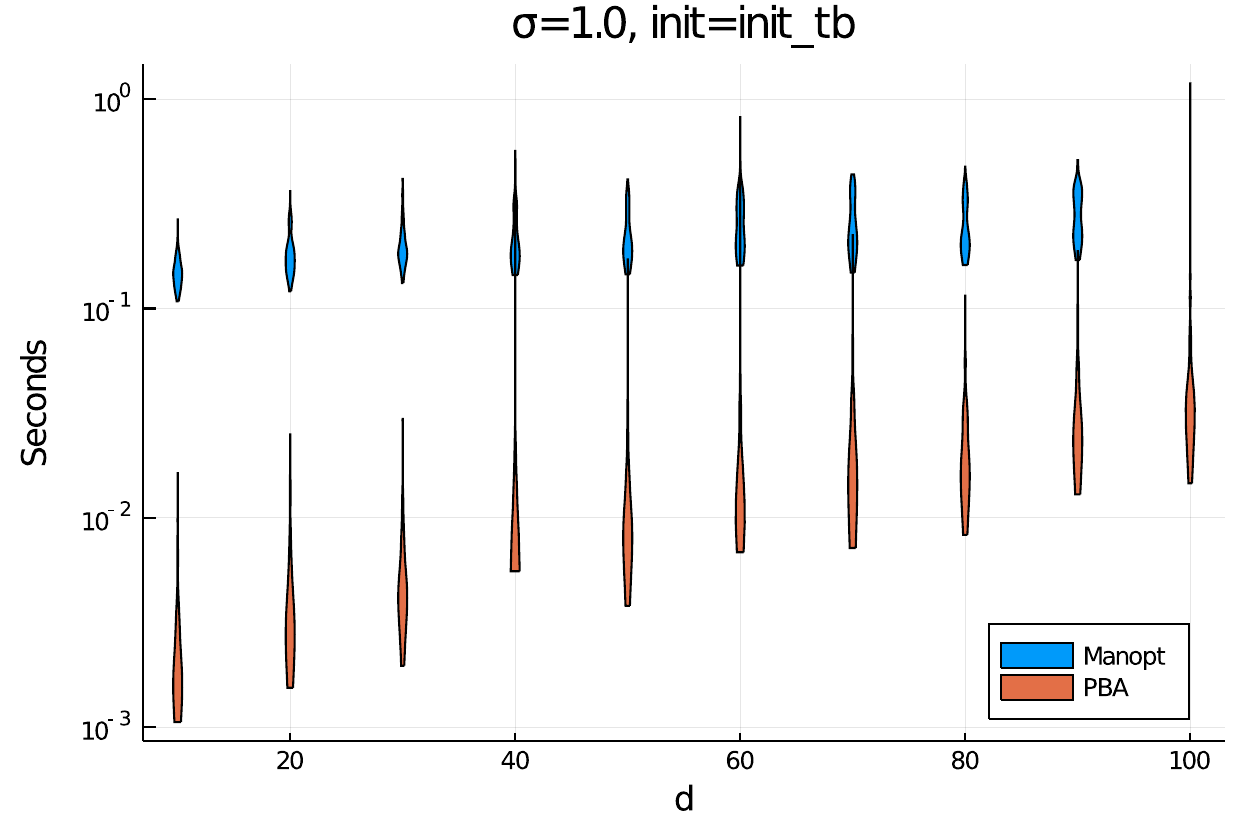}
\\
\includegraphics[width=0.45\textwidth]{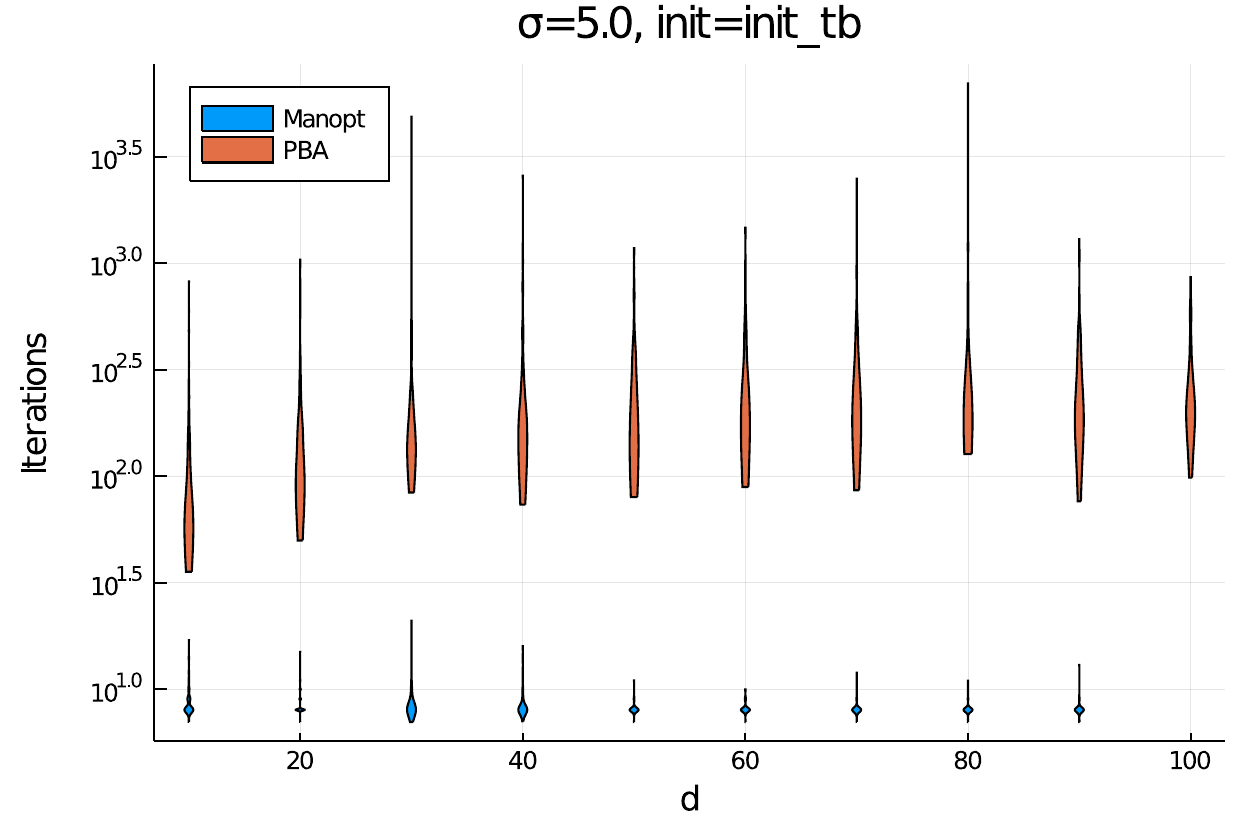}
&
\includegraphics[width=0.45\textwidth]{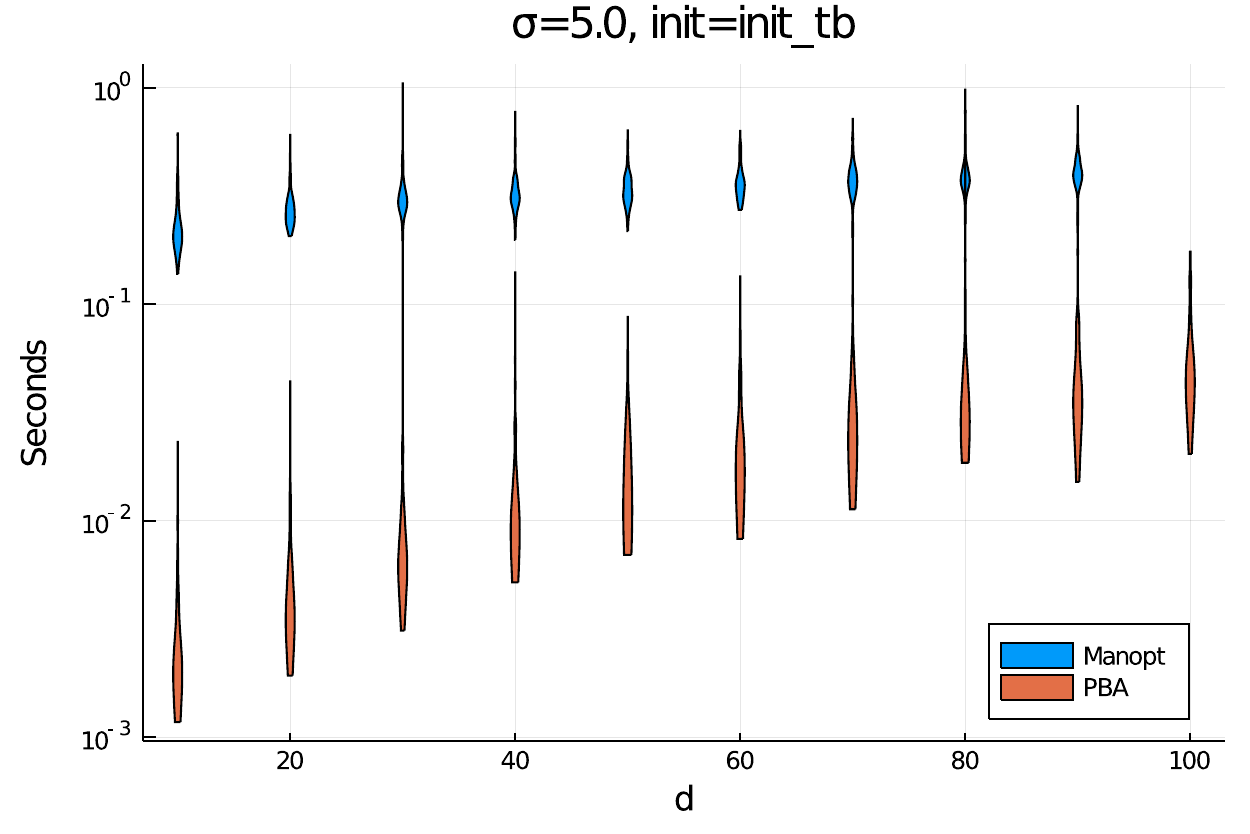}
\\
\includegraphics[width=0.45\textwidth]{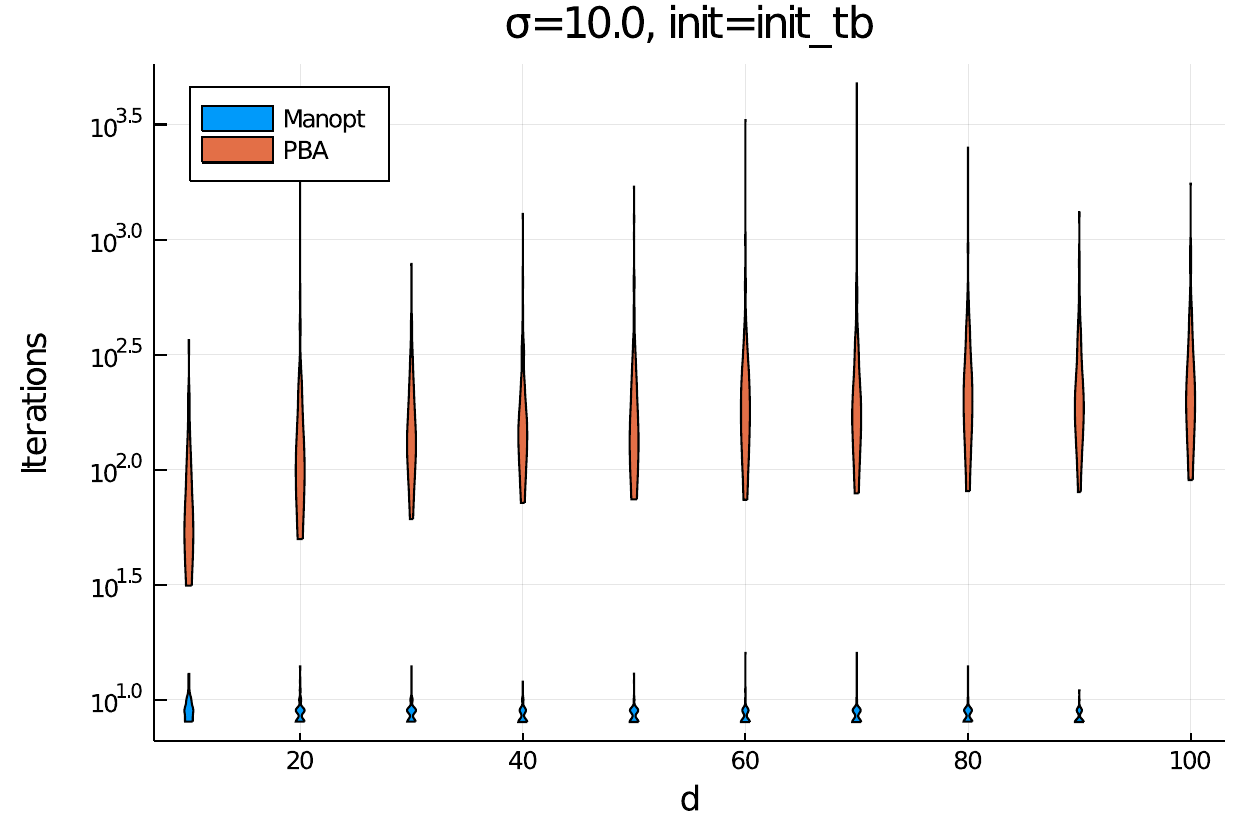}
&
\includegraphics[width=0.45\textwidth]{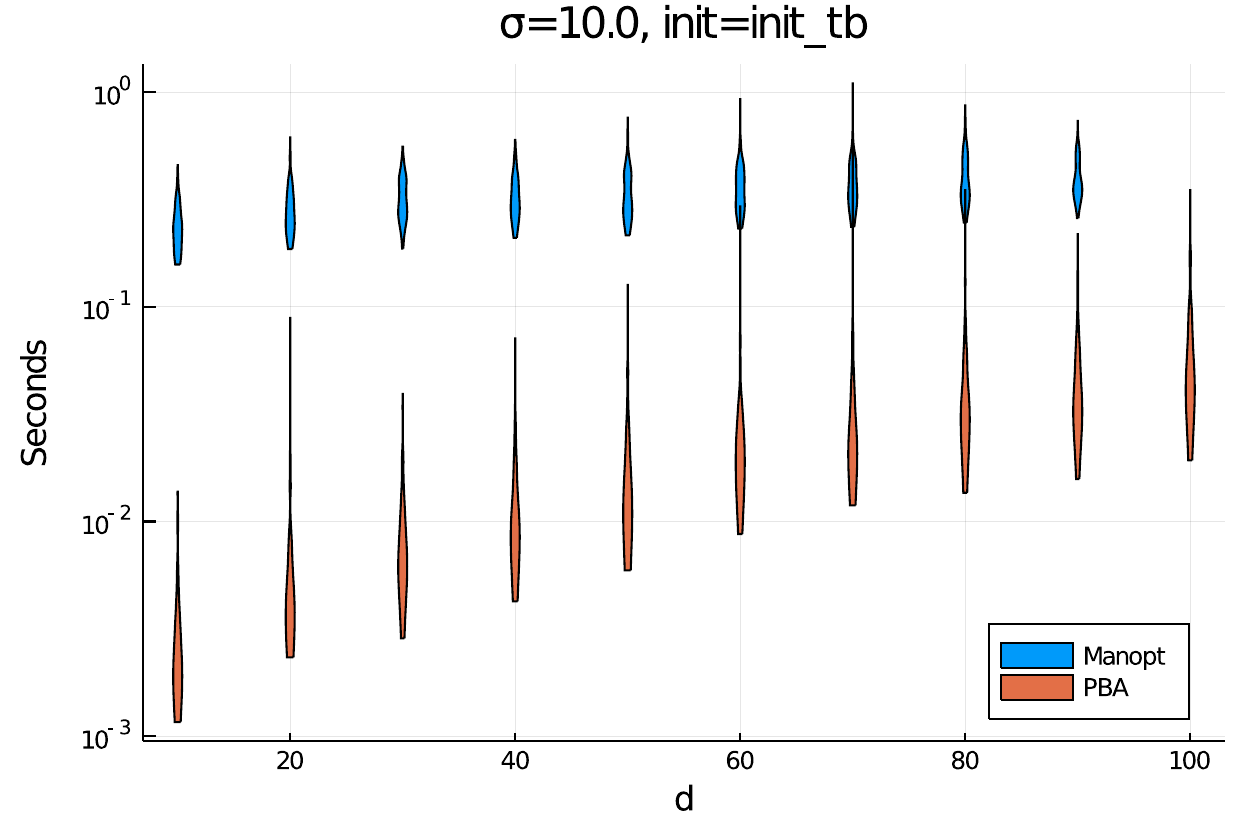}
\end{tabular}
\end{center}
\caption{Simulation studies. Left: number of iterations until convergence for each $d$. Right: time (in seconds) until convergence. Initialization strategy ``tb'' was used.}\label{fig:summary_r3m5supp}
\end{figure}

\paragraph{Additional tables and plots}
Table \ref{tbl:portwine_maxbet}, generalized CCA of Port Wine Data (Example \ref{ex:portwine} for the MAXBET criterion);
Fig. \ref{fig:summary_r3m5supp}, violin plots of number of iterations, and timing vs. dimension $d$ for the large scale examples in \S \ref{sec:numerical:large}.
 \begin{table*}
\begin{center}
\caption{Example \ref{ex:liu5_1}: MAXDIFF}\label{tbl:liu5_1_maxdiff}
\medskip
\footnotesize
\begin{tabular}{lllrlllll}
\toprule
$r$ & init & method & iter & time (sec) & obj & classification & $\lambda {\min}(\bL^{\star})$ & $\lambda {\min}(\mathcal{L}^{\star})$ \\
\midrule 
\multirow{8}{*}{1} & \multirow{2}{*}{eye} & PBA & 8 & 9.759e-5 & 1.870 & global opt. & -2.987e-16 & -2.987e-16 \\
 &  & Manopt & 9 & 0.02981 & 1.870 & global opt. & -4.652e-16 & -4.652e-16 \\
 & \multirow{2}{*}{sb} & PBA & 9 & 9.784e-5 & 1.870 & global opt. & -6.317e-16 & -6.317e-16 \\
 &  & Manopt & 7 & 0.02202 & 1.870 & global opt. & -4.018e-15 & -4.018e-15 \\
 & \multirow{2}{*}{tb} & PBA & 2 & 5.466e-5 & 1.870 & global opt. & -1.171e-15 & -1.171e-15 \\
 &  & Manopt & 1 & 0.0007613 & 1.870 & global opt. & \;4.989e-17 & \;4.989e-17 \\
 & \multirow{2}{*}{lww1} & PBA & 9 & 0.0001105 & 1.870 & global opt. & -4.652e-16 & -4.652e-16 \\
 &  & Manopt & 7 & 0.01505 & 1.870 & global opt. & -5.052e-16 & -5.052e-16 \\
\cmidrule{2-9}
\multirow{8}{*}{2} & \multirow{2}{*}{eye} & PBA & 5 & 0.0001356 & 2.265 & global opt. & -1.350e-15 & -1.475 \\
 &  & Manopt & 8 & 0.02020 & 2.265 & global opt. & -9.940e-15 & -1.475 \\
 & \multirow{2}{*}{sb} & PBA & 5 & 9.158e-5 & 2.265 & global opt. & -1.525e-16 & -1.475 \\
 &  & Manopt & 7 & 0.01931 & 2.265 & global opt. & -2.487e-15 & -1.475 \\
 & \multirow{2}{*}{tb} & PBA & 2 & 6.768e-5 & 2.265 & global opt. & -4.513e-16 & -1.475 \\
 &  & Manopt & 1 & 0.0007638 & 2.265 & global opt. & \;6.422e-17 & -1.475 \\
 & \multirow{2}{*}{lww1} & PBA & 5 & 0.0001288 & 2.265 & global opt. & -3.695e-16 & -1.475 \\
 &  & Manopt & 10 & 0.03455 & 2.265 & global opt. & -8.954e-16 & -1.475 \\
\bottomrule
\end{tabular}
\end{center}
\end{table*}

\begin{table*}
\begin{center}
\caption{Example \ref{ex:liu5_1}, MAXBET}\label{tbl:liu5_1_maxbet}
\medskip
\footnotesize
\begin{tabular}{lllrlllll}
\toprule
$r$ & init & method & iter & time (sec) & obj & classification & $\lambda {\min}(\bL^{\star})$ & $\lambda {\min}(\mathcal{L}^{\star})$ \\
\midrule 
\multirow{8}{*}{1} & \multirow{2}{*}{eye} & PBA & 101 & 0.0004694 & 7.051 & not loc. opt. & \text{--} & \text{--} \\
 &  & Manopt & 8 & 0.03417 & 7.051 & not loc. opt. & \text{--} & \text{--} \\
 & \multirow{2}{*}{sb} & PBA & 80 & 0.0004606 & 7.365 & global opt. & -5.231e-15 & -5.231e-15 \\
 &  & Manopt & 8 & 0.03018 & 7.365 & global opt. & \;3.573e-15 & \;3.573e-15 \\
 & \multirow{2}{*}{tb} & PBA & 64 & 0.0003227 & 7.365 & global opt. & -5.054e-15 & -5.054e-15 \\
 &  & Manopt & 7 & 0.02155 & 7.365 & global opt. & -3.141e-15 & -3.141e-15 \\
 & \multirow{2}{*}{lww1} & PBA & 67 & 0.0003648 & 7.365 & global opt. & -1.688e-15 & -1.688e-15 \\
 &  & Manopt & 6 & 0.02293 & 7.365 & global opt. & \;5.923e-16 & \;5.923e-16 \\
\cmidrule{2-9} 
\multirow{8}{*}{2} & \multirow{2}{*}{eye} & PBA & 49 & 0.0003298 & 12.75 & global opt. & -1.153e-15 & -3.269 \\
 &  & Manopt & 8 & 0.02223 & 12.75 & global opt. & -7.428e-16 & -3.269 \\
 & \multirow{2}{*}{sb} & PBA & 44 & 0.0003020 & 12.75 & global opt. & -2.203e-15 & -3.269 \\
 &  & Manopt & 6 & 0.01648 & 12.75 & global opt. & -6.327e-15 & -3.269 \\
 & \multirow{2}{*}{tb} & PBA & 44 & 0.0002937 & 12.75 & global opt. & \;3.348e-16 & -3.269 \\
 &  & Manopt & 5 & 0.01177 & 12.75 & global opt. & -1.458e-15 & -3.269 \\
 & \multirow{2}{*}{lww1} & PBA & 62 & 0.0004145 & 12.15 & stationary & -0.599 & -3.438 \\
 &  & Manopt & 6 & 0.01761 & 12.15 & stationary & -0.599 & -3.438 \\
\bottomrule
\end{tabular}
\end{center}
\end{table*}

 \begin{table*}
\begin{center}
\caption{Example \ref{ex:liu5_2}, MAXDIFF}\label{tbl:liu5_2_maxdiff}
\medskip
\footnotesize
\begin{tabular}{lllrlllll}
\toprule
$r$ & init & method & iter & time (sec) & obj & classification & $\lambda {\min}(\bL^{\star})$ & $\lambda {\min}(\mathcal{L}^{\star})$ \\
\hline
\multirow{8}{*}{1} & \multirow{2}{*}{eye} & PBA & 9 & 0.0001254 & 66.57 & global opt. & -1.744e-14 & -1.744e-14 \\
 &  & Manopt & 8 & 0.04351 & 66.57 & global opt. & -1.899e-15 & -1.899e-15 \\
 & \multirow{2}{*}{sb} & PBA & 10 & 0.0001336 & 66.57 & global opt. & \;2.135e-14 & \;2.135e-14 \\
 &  & Manopt & 8 & 0.03467 & 66.57 & global opt. & -4.456e-16 & -4.456e-16 \\
 & \multirow{2}{*}{tb} & PBA & 9 & 0.0001278 & 66.57 & global opt. & -2.231e-15 & -2.231e-15 \\
 &  & Manopt & 4 & 0.01699 & 66.57 & global opt. & -5.917e-15 & -5.917e-15 \\
 & \multirow{2}{*}{lww1} & PBA & 13 & 0.0001677 & 66.57 & global opt. & -5.44e-16 & -5.44e-16 \\
 &  & Manopt & 10 & 0.03902 & 66.57 & global opt. & \;1.797e-14 & \;1.797e-14 \\
\cmidrule{2-9} 
\multirow{8}{*}{2} & \multirow{2}{*}{eye} & PBA & 11 & 0.0001296 & 93.05 & global opt. & \;2.225e-15 & -28.54 \\
 &  & Manopt & 8 & 0.02821 & 79.75 & not loc. opt. & \text{--} & \text{--} \\
 & \multirow{2}{*}{sb} & PBA & 11 & 0.0001183 & 93.05 & global opt. & -1.855e-15 & -28.54 \\
 &  & Manopt & 6 & 0.02081 & 93.05 & global opt. & -1.13e-15 & -28.54 \\
 & \multirow{2}{*}{tb} & PBA & 10 & 0.0001076 & 93.05 & global opt. & -9.259e-15 & -28.54 \\
 &  & Manopt & 5 & 0.01888 & 93.05 & global opt. & -1.727e-14 & -28.54 \\
 & \multirow{2}{*}{lww1} & PBA & 12 & 0.0001502 & 93.05 & global opt. & -3.618e-14 & -28.54 \\
 &  & Manopt & 9 & 0.02662 & 93.05 & global opt. & -1.846e-14 & -28.54 \\
\bottomrule
\end{tabular}
\end{center}
\end{table*}

\begin{table*}
\begin{center}
\caption{Example \ref{ex:liu5_2}, MAXBET}\label{tbl:liu5_2_maxbet}
\medskip
\footnotesize
\begin{tabular}{lllrlllll}
\toprule
$r$ & init & method & iter & time (sec) & obj & classification & $\lambda {\min}(\bL^{\star})$ & $\lambda {\min}(\mathcal{L}^{\star})$ \\
\midrule
\multirow{8}{*}{1} & \multirow{2}{*}{eye} & PBA & 20 & 0.0001750 & 189.5 & stationary & -0.4819 & -0.4819 \\
 &  & Manopt & 9 & 0.05333 & 189.5 & stationary & -0.4819 & -0.4819 \\
 & \multirow{2}{*}{sb} & PBA & 21 & 0.0001798 & 189.5 & stationary & -0.4819 & -0.4819 \\
 &  & Manopt & 9 & 0.04469 & 189.5 & stationary & -0.4819 & -0.4819 \\
 & \multirow{2}{*}{tb} & PBA & 16 & 0.0001653 & 189.5 & stationary & -0.4819 & -0.4819 \\
 &  & Manopt & 4 & 0.01736 & 189.5 & stationary & -0.4819 & -0.4819 \\
 & \multirow{2}{*}{lww1} & PBA & 18 & 0.0002036 & 189.5 & stationary & -0.4819 & -0.4819 \\
 &  & Manopt & 6 & 0.02912 & 189.5 & stationary & -0.4819 & -0.4819 \\
\cmidrule{2-9} 
\multirow{8}{*}{2} & \multirow{2}{*}{eye} & PBA & 45 & 0.0003748 & 250.2 & stationary & -12.65 & -109.4 \\
 &  & Manopt & 8 & 0.02823 & 250.2 & stationary & -12.65 & -109.4 \\
 & \multirow{2}{*}{sb} & PBA & 30 & 0.0002464 & 263.6 & global opt. & -2.863e-13 & -102.9 \\
 &  & Manopt & 6 & 0.01766 & 263.6 & global opt. & -7.745e-15 & -102.9 \\
 & \multirow{2}{*}{tb} & PBA & 31 & 0.0002772 & 263.6 & global opt. & -8.992e-15 & -102.9 \\
 &  & Manopt & 6 & 0.01924 & 263.6 & global opt. & -6.802e-15 & -102.9 \\
 & \multirow{2}{*}{lww1} & PBA & 36 & 0.0003419 & 263.6 & global opt. & -1.158e-13 & -102.9 \\
 &  & Manopt & 9 & 0.02834 & 263.6 & global opt. & -2.375e-14 & -102.9 \\
\bottomrule
\end{tabular}
\end{center}
\end{table*}

 \begin{table*}
\begin{center}
\caption{Port Wine Data, MAXBET}\label{tbl:portwine_maxbet}
\medskip
\footnotesize
\begin{tabular}{lllrlllll}
\toprule
$r$ & init & method & iter & time (sec) & obj & classification & $\lambda {\min}(\bL^{\star})$ & $\lambda {\min}(\mathcal{L}^{\star})$ \\
\midrule 
\multirow{8}{*}{1} & \multirow{2}{*}{eye} & PBA & 13 & 0.0002722 & 302.8 & global opt. & -3.351e-15 & -3.351e-15 \\
 &  & Manopt & 9 & 0.1007 & 302.8 & global opt. & \;1.581e-14 & \;1.581e-14 \\
 & \multirow{2}{*}{sb} & PBA & 13 & 0.0002102 & 302.8 & global opt. & -4.491e-14 & -4.491e-14 \\
 &  & Manopt & 7 & 0.07019 & 302.8 & global opt. & \;1.581e-14 & \;1.581e-14 \\
 & \multirow{2}{*}{tb} & PBA & 11 & 0.0001744 & 302.8 & global opt. & -8.152e-15 & -8.152e-15 \\
 &  & Manopt & 5 & 0.05140 & 302.8 & global opt. & \;5.021e-15 & \;5.021e-15 \\
 & \multirow{2}{*}{lww1} & PBA & 14 & 0.0002237 & 302.8 & global opt. & \;2.326e-14 & \;2.326e-14 \\
 &  & Manopt & 10 & 0.07510 & 302.8 & global opt. & -1.779e-14 & -1.779e-14 \\
\cmidrule{2-9} 
\multirow{8}{*}{2} & \multirow{2}{*}{eye} & PBA & 28 & 0.0005198 & 389.9 & global opt. & \;7.302e-15 & -121.9 \\
 &  & Manopt & 11 & 0.1146 & 389.9 & global opt. & \;1.492e-14 & -121.9 \\
 & \multirow{2}{*}{sb} & PBA & 28 & 0.0005463 & 389.9 & global opt. & -4.004e-14 & -121.9 \\
 &  & Manopt & 7 & 0.1310 & 389.9 & global opt. & -1.103e-13 & -121.9 \\
 & \multirow{2}{*}{tb} & PBA & 26 & 0.0004592 & 389.9 & global opt. & -1.189e-14 & -121.9 \\
 &  & Manopt & 7 & 0.1233 & 389.9 & global opt. & -5.900e-14 & -121.9 \\
 & \multirow{2}{*}{lww1} & PBA & 27 & 0.0004891 & 389.9 & global opt. & -4.200e-14 & -121.9 \\
 & & Manopt & 11 & 0.1486 & 389.9 & global opt. & -1.976e-14 & -121.9 \\
\cmidrule{2-9} 
\multirow{8}{*}{3} & \multirow{2}{*}{eye} & PBA & 41 & 0.001651 & 409.4 & stationary & -3.332 & -162.7 \\
 &  & Manopt & 15 & 0.3253 & 409.4 & stationary & -3.332 & -162.7 \\
 & \multirow{2}{*}{sb} & PBA & 15 & 0.0005503 & 414.9 & global opt. & -1.575e-13 & -159.3 \\
 &  & Manopt & 7 & 0.1370 & 414.9 & global opt. & -3.376e-14 & -159.3 \\
 & \multirow{2}{*}{tb} & PBA & 17 & 0.0006829 & 414.9 & global opt. & -1.226e-13 & -159.3 \\
 &  & Manopt & 7 & 0.1306 & 414.9 & global opt. & -3.034e-14 & -159.3 \\
 & \multirow{2}{*}{lww1} & PBA & 38 & 0.001179 & 409.4 & stationary & -3.332 & -162.7 \\
 &  & Manopt & 17 & 0.2742 & 409.4 & stationary & -4.288 & -162.8 \\
\bottomrule
\end{tabular}
\end{center}
\end{table*}

\end{document}